\relax
\documentclass[]{article} 
\usepackage{times}  % DO NOT CHANGE THIS
\usepackage{helvet}  % DO NOT CHANGE THIS
\usepackage{courier}  % DO NOT CHANGE THIS
\usepackage[hyphens]{url}  % DO NOT CHANGE THIS
\usepackage{graphicx} % DO NOT CHANGE THIS
\usepackage{tikz} % nice language for creating drawings and diagrams
\usetikzlibrary{shapes.geometric, arrows,automata,petri}

\tikzset{
    square/.style={%
        draw=none,
        circle,
        append after command={%
            \pgfextra \draw[black] (\tikzlastnode.north-|\tikzlastnode.west) rectangle 
                (\tikzlastnode.south-|\tikzlastnode.east);\endpgfextra}
    }
}

\usetikzlibrary{positioning}
\usepackage{caption} % DO NOT CHANGE THIS AND DO NOT ADD ANY OPTIONS TO IT
\DeclareCaptionStyle{ruled}{labelfont=normalfont,labelsep=colon,strut=off} % DO NOT CHANGE THIS
\newcommand{\YD}[1]{{ \color{Blue}{\bf [YD:} #1]}}
\newcommand{\YDT}[1]{{ \color{Blue}{\bf [YD:]} #1}}

\newcommand{\ES}[1]{{ \color{Purple}{\bf [ES:} #1]}}

\newcommand{\remove}[1]{}

\usepackage{amssymb}
\usepackage{amsthm}
\usepackage{mathtools}
\usepackage{algorithm} 
\usepackage{algpseudocode} 

\usepackage{color}
\usepackage{graphicx}
\usepackage{relsize}
\usepackage{mathtools}  
\usepackage{amsmath}
\usepackage{amssymb}
\usepackage{tabulary}
\usepackage{booktabs}
\definecolor{brown}{cmyk}{0,0.5,1,0.2}
\definecolor{myGreen}{rgb}{0,0.8,0.4}
\definecolor{Green}{rgb}{0,1,0}
\definecolor{Red}{rgb}{1,0,0}
\definecolor{myBlue}{rgb}{0,0.6,1}
\definecolor{Blue}{rgb}{0,0,1}
\definecolor{Purple}{rgb}{0.7,0.0,1.0}

\newtheorem{definition}{Definition}%[section]
\newtheorem{theorem}{Theorem}[section]
\newtheorem{lemma}{Lemma}[section]
\newtheorem{observation}[theorem]{Observation}
\newtheorem{proposition}[theorem]{Proposition}
\newtheorem{property}[lemma]{Property}

\title{On Existence of Must-Include Paths and Cycles\\ in Undirected Graphs
\footnote{Manuscript submitted for publication in SICOMP}}

\date{}

\author{Yefim Dinitz, dinitz@bgu.ac.il\\ Solomon Eyal Shimony, shimony@cs.bgu.ac.il\\
Dept. of Computer Science,
Ben-Gurion University}

\begin{document}

%\begin{keyword}
%Biconnected graphs,
%triconnected graphs,
%SPQR trees,
%must-include simple paths and cycles
%\end{keyword}

\maketitle

\begin{abstract}
Given an undirected graph $G=(V,E)$ and
vertices $s,t,w_1,w_2\in V$, we study finding whether there exists a simple path $P$ from $s$ to $t$
such that $w_1,w_2 \in P$.
As a sub-problem, we study the question: given an undirected graph and three of its edges, does there exist a simple cycle containing all those edges?
We provide necessary and sufficient conditions for the existence of such paths and cycles,
and develop efficient algorithms to solve this and related problems.
%\YDT {Additionally, we present a new description of the SPQR tree of a graph, which eliminates the post-processing for generating its S nodes.}\ES{I would avoid placing it like this in the abstract like the plague! Mentioning it in weak language in contributions "as far as we know" is one thing, and is OK. Placing it like this in the abstract is very risky. Suppose Roberto Tammasia, for example, has a similar description in his new version of his textbook? If he gets this to review this almost asks him to reject the paper,,,}

$~$

\end{abstract}

\section{Introduction}

We define and examine the problem of determining 
whether a simple path from a source to a target
through two other given vertices exists
in an undirected graph.
This problem, as well as
its extension to finding all vertex pairs where such
paths do not exist,
arose from the need to quickly determine must-include vertex pairs in an admissible
heuristic in a combinatorial search for longest constrained paths \cite{ShimonyEtAl2022SOCS}.
Finding simple paths that include a given
set of vertices is also of interest
in communication networks  \cite{5366808}.
In addition,
this problem, as well as the 
related question on a must-include simple cycle, are of theoretical interest. 
%\YD {Delete: in combinatorics. Maybe, the graph theory is not less relevant than combinatorics.}

%*****************************************
\subsection{Problem statements}
\label{ss:problem statements}

Formally, we examine the following 
problems on simple paths and cycles in undirected graphs, starting with the main problem:
%\YD {Maybe, delete?: , starting with the main problem: }

\begin{definition}[Path Existence Problem (PEP)]
Given an undirected graph $G=(V,E)$ and
vertices $s,t,w_1,w_2\in V$, does there exist
a simple path $P$ from $s$ to $t$
such that $w_1,w_2 \in P$?
\end{definition}

%\YD {I suggest eliminating the word "Definition", retaining just the problem names, in bold. Decide as you wish, Eyal.}

We henceforth call a such a path a
\emph{must-include} $\{ w_1,w_2\}$-path from $s$ to $t$. In triconnected graphs, the PEP always has an affirmative answer \cite{ShimonyEtAl2022SOCS}.
A more challenging alternate of this problem
that is also of interest is:

\begin{definition}[Excluded Pairs Enumeration
(EPE)]
Given an undirected graph $G=(V,E)$ and
vertices $s,t\in V$,
find all pairs of vertices $w_1,w_2$
such that there is no simple
path $P$ from $s$ to $t$ such that $\{w_1,w_2\} \in P$.
\end{definition}

For conciseness, we call the vertices and edges of a graph its \emph{elements}.
Another problem similar to PEP is:

%the \emph{elements} of a graph.
%As a sub-problem having an independent interest, we also study the following:

\begin{definition}[Cycle Existence Problem (CEP)]
Given an undirected graph and three of its elements, does there exist a simple cycle containing all these elements?
\end{definition}

Given a graph $G$ with distinguished vertices $s$ and $t$ (either with or without edge $(s,t)$), 
we denote $G \cup (s,t)$ by $G^+$.\footnote
{Henceforth, we omit the braces $\{ ~\}$ for singleton sets in expressions involving set operators such as $\setminus$ and $\cup$, when unambiguous.
Additionally, we apply set operators to a graph as shorthand for the operation applied to its vertex set or edge set, as appropriate from the context. E.g., for a graph $G=(V,E)$, we denote by $G \cup (s,t)$ the graph $(V, E\cup \{ (s,t)\})$.} 
Note that PEP$(G,s,t,w_1,w_2)$ is equivalent to CEP$(G^+,(s,t),w_1,w_2)$: their solutions differ just by adding/removing edge $(s,t)$.
Due to this reduction, any result for the CEP straightforwardly implies 
a similar result for the PEP, including algorithms and their run-times.
Hence, in what follows, we provide theoretical results and algorithmic solutions only for the CEP, 
thereby immediately covering the PEP as well.

For the CEP on a triconnected graph, the answer is known (as a folklore), except for the case when three edges are given: in all but the latter case, such a cycle always exists. Where all given elements are vertices, this is covered by Dirac's theorem \cite{Dirac1960}, while the case where at most two elements are edges was proved in informal discussions \cite{k-connected, k-connect}.
%\footnote{
%We found it at https://math.stackexchange.com/questions/156032/every-k-vertices-in-an-k-connected-graph-are-contained-in-a-cycle?rq=1 \cite{k-connected} and https://math.stackexchange.com/questions/3599622/k-connected-graph-there-exists-a-cycle-that-contains-any-2-edges-and-any-k-2?noredirect=1\&lq=1 \cite{k-connected-two-edges}.} Misha Lavrov
%https://math.stackexchange.com/questions/3600673/g-a-k-connected-graph-show-that-a-set-of-k-2-vertices-and-a-set-of-two-edges-li?rq=1
%Nick Matteo
%In Section \ref{s:background}, we provide an outline of the proof.(?????)

\remove{
Our results are as follows:
\begin{itemize}
    \item For the CEP on a triconnected graph with three distinguished edges, we provide the necessary and sufficient condition for the cycle existence. 
    \item For the general PEP and CEP, %on an arbitrary graph, 
    we provide the necessary and sufficient condition for the cycle existence.
    \item For the general PEP, CEP, and EPE,
    we provide efficient (linear time) algorithms deciding them.
\end{itemize}
}

%*******************************************
\subsection{Related work}

Here we briefly discuss the work on
on problems related to the PEP".
and related problems.
Issues directly related to our methods and used in developing our results are examined in more
detail in the background section.

Much of the literature on existence of simple paths
centers on disjoint paths with certain
endpoint requirements \cite{berczi_et_al:LIPIcs:2017:7824,EILAMTZOREFF1998113}.
Out of these, the most relevant to the
PEP is: given $k$ pairs of vertices
$(v_{j_1},v_{j_2})$,  $k\geq j\geq 1$,
find $k$  (edge-wise or vertex-wise)
disjoint paths such that each path $P_j$ connects the pairs $v_{j_1},v_{j_2}$.
These problems are NP-hard  for general $k$ \cite{10.1145/1061425.1061430}, but for 
constant $k$, Robertson and Seymour
\cite{DBLP:journals/jct/RobertsonS95b}
provide an $O(|V|^3)$ solution to the problem
of $k$ vertex-disjoint paths. However, as the authors admit,
the algorithm stated therein is impractical. 
Given an algorithm for computing these
$k$ disjoint paths, one can solve the PEP
by finding a path from $s$ to $w_1$,
a path from $w_2$ to $t$, and a path from
a neighbor of $w_1$ to a neighbor of $w_2$,
all disjoint. Merging these paths is a solution to the PEP. 
Making this scheme complete may require iterating over
all pairs of neighbors, and also examining
the case where the paths are $s$ to $w_2$
and $t$ to $w_1$. Still, the runtime is polynomial. Finding $k$ disjoint trees
can be done in linear time (for fixed $k$), in graphs restricted to be planar \cite{DisjointTrees}.
%\YD {Delete this sentence? I do not think that the mentioned problem is related to our ones.} \ES{Actually it is related, as it has k disjoint paths as a special case, and as stated above, from that you can answer the PEP. But this result only works for planar graphs}

The direct generalization of the PEP,
where the must-include set consists
of more than two vertices, has also been examined. 
For unbounded $k$, it is easily
shown that this problem is NP-hard, by reduction from Hamiltonian Cycle
\cite{garey1979computers}.
An algorithm using path splicing involving small $k$ was shown in \cite{5366808},
but that work provides no explicit runtime guarantees and
their algorithm does not guarantee finding
a path if one exists.

Our PEP falls under the conditions
for existence of an integral multi-flow
in a network with four terminals
stated by Seymour \cite{DBLP:journals/networks/Seymour80}.
However, it is not explicitly stated therein
 how that translates into an efficient algorithm for finding such paths or deciding
their existence.
%\YD {Delete?: In triconnected graphs, the PEP always has an affirmative answer \cite{ShimonyEtAl2022SOCS}.
%For the CEP, the solution is known (as a folklore), except for the case of three edges: % That is,
%in all but the latter case, such a cycle always exists.\footnote{
%We found it at 
%https://math.stackexchange.com/questions/156032/every-k-vertices-in-an-k-connected-graph-are-contained-in-a-cycle?rq=1 \cite{k-connected} and 
%https://math.stackexchange.com/questions/3599622/k-connected-graph-there-exists-a-cycle-that-contains-any-2-edges-and-any-k-2?noredirect=1\&lq=1 \cite{k-connected-two-edges}} In Section \ref{s:background}, we provide an outline of the proof.}
%This solution allows an effective verifying algorithm.

%*****************************************
\subsection{Structure and contributions of this paper}

%\YD {I suggest switching this subsection with the Related Work one, since it is central and more important.}

We begin with background on graph connectivity (Section \ref{s:background}): Menger's theorem and the decomposition of a graph into its 2- and 3-connected components and the respective properties of simple paths within them;
in passing, we show a simple reduction of our problems to the case of a biconnected graph.
In particular, we describe SPQR trees (Section \ref{sec:SPQR}), which are heavily used in this paper.
For the reader's convenience, we provide an SPQR tree description not coinciding with any previously known one, though essentially equivalent to that of \cite{SPQRtrees}.

%As far as we know, our description is  novel, as it eliminates the (post-processing) merge operations required in related work \cite{SPQRtrees} to define the S nodes in the SPQR tree.

%, examining in detail SPQR trees which are heavily used in this paper, and other necessary properties of simple paths 
For CEP, we provide (Section \ref{sec:triconnected}, the first contribution of this paper)
a structural solution to the case of three edges in triconnected graphs: 
simple necessary and sufficient conditions for cycle existence.
%\YD {not in Section 3:  and a decision algorithm with a linear complexity.
%}
(Note that Section~\ref{sec:triconnected} can be read before Section \ref{sec:SPQR}, as it does not require knowing SPQR trees.)
%only Menger's theorem, which we revisit in Section \ref{ss:basics}.

We then provide necessary and sufficient conditions for the existence of a cycle as in CEP (and thus of the path as in PEP)
in general biconnected undirected graphs
(Section \ref{sec:biconnected}, the second contribution of this paper).
This is done by
%converting the PEP to a CEP, and then 
showing that the CEP reduces
to cycle existence in an appropriately
defined \emph{central component} of the SPQR tree of the graph.
The difficult subproblem here is
in components of type R,
where we use the above-stated conditions
for CEP in triconnected graphs.

Finally, we develop efficient algorithms for both CEP and EPE (Section \ref{sec:algorithms}).
Our CEP (and thus PEP) algorithm has runtime $O(|E|)$,
assuming the input graph is connected.
For the EPE in the case where \emph{implicit} output of the enumerated pairs is allowed (for example, the admissible heuristic in \cite{ShimonyEtAl2022SOCS}
only needs access to various counts of such pairs,
thus can use the implicit form directly), we provide an $O(|E|)$ algorithm as well;
it becomes $O(|V|^2)$ if \emph{explicit} enumeration is required.

%All graphs referred to in this paper are undirected.

%*****************************************
\section{Background on Graph Connectivity}
\label{s:background}

In this paper, we consider only connected undirected graphs $G=(V,E)$.
We begin with some background on vertex connectivity and the supporting flow algorithms;
edge connectivity is also mentioned as needed.
See, e.g., books \cite{even2011graph,Harary1969,10.5555/1051910} for details.
A subset $V' \subseteq V$ of vertices is called a  \emph{$($vertex$)$ cut} if the removal of the vertices in $V'$ makes $G$ not connected and if no subset of $V'$ has this property. 
We call a cut $V'$ a ``$k$-cut'' if $|V'| =k$.
A graph is (vertex) \emph{$k$-connected} if either $|V| \ge k+1$ and no $k-1$ or fewer of its vertices form a cut, or if it is a complete graph on $k$ vertices. 
We use the terms \emph{biconnected} and \emph{triconnected graphs} meaning 2- and 3-connected graphs, respectively. 
Two or more paths are called \emph{internally vertex-disjoint} if they are vertex
%\YD {vertex --$>$ fully} 
%\ES{But this gives us "fully disjoint", does this have an agreed meaning?}
%\YD {Maybe, remove both "vertex" and "fully"?}
disjoint except possibly for their end-vertices.
%disjoint except their end-vertices (also called terminals). We call such paths {\em internally vertex disjoint} paths

%**************************************
\subsection{Basics of Connectivity}
\label{sec:basics}\label{ss:basics}

We begin with some properties and theorems in the graph connectivity theory and algorithms that are used in the rest of the paper.
The main tool for exploring $k$-connected graphs is 
%the \YD {Delete: (generalized) [[since below we write "generalization" and "extension" of Menger's theorem]]}
Menger's theorem:

\begin{theorem}[Menger \cite{Menger1927}]
\label{th:Menger}
In any $k$-connected graph $G=(V,E)$,
$|V| \ge k+1$, for every two vertices $s,t \in V$,
there exist $k$ simple, internally vertex-disjoint paths between $s$ and $t$.
\end{theorem}

Additionally, for two sets of vertices $V_s,V_t$ (``sides''),
each of cardinality at least $k$, there
exist $k$ simple, fully disjoint paths between $V_s$ and $V_t$. % vertex-disjoint at their terminals, as well as internally.
(Note that the sides need not be disjoint.)
Throughout this paper, we assume that \emph{no internal vertex of those paths is in $V_s,V_t$}; indeed, in such a case, the extra paths' prefixes and/or suffixes can be safely truncated.
For conciseness, we henceforth
take paths between two objects (such as paths, cycles) to mean paths between the vertex sets of these objects.

The theorem also naturally generalizes to the case
where either side is a singleton set,
e.g., $V_s=\{s\}$: there exist $k$
paths as above that are vertex-disjoint except at that
singleton set.
We refer to the above generalizations of 
Theorem \ref{th:Menger} as the
\emph{extension to Menger's theorem}.
%Also the cases $s, V_t$ and $V_s, t$ are admissible, with a natural generalization.

%%**************************************
%\subsection{Basics of Connectivity Related Algorithms}
%\label{ss:basics-alg}

%\YD {A temporary place; I intend to think on the right one.}

When we need to actually find the set of
vertex-disjoint paths between $s$ and $t$ as in Menger's theorem, the following efficient algorithms are used.
To find a path from a given source vertex $s$ to a target vertex $t$ (the case $k=1$ of Menger's theorem), efficient schemes execute a labeling algorithm (either BFS or DFS) scanning $G$ from $s$ up to labeling $t$, and then restore the path 
using the backward labels; this takes time $O(|E|)$.
(See any book in basic graph algorithms as
a source, such as \cite{10.5555/1614191}.) 

When we need to find the set of several
vertex-disjoint paths between $s$ and $t$ as in  Menger's theorem, network flow techniques are used.
The given $k$-connected graph $G=(V,E)$, $|V| \ge k+1$, is turned into a (directed) flow network $N^{\mbox {v}} =(\vec G^{\mbox {v}},s,t,c^{\mbox {v}})$ with the flow source $s$ and sink $t$.
Every vertex $v$ of $G$, except for $s$ and $t$, is turned into two vertices $v_1$ and $v_2$ with a directed edge $(v_1,v_2)$ of capacity $c^{\mbox {v}}(v_1,v_2)=1$ from $v_1$ to $v_2$ in $\vec G^{\mbox {v}}$. Every (undirected) edge $(u,v)$ of $G$ is turned into two directed edges $(u_2,v_1)$ and $(v_2,u_1)$ of infinite capacity each in $\vec G^{\mbox {v}}$.
The following max-flow min-cut theorem \cite{10.5555/1942094} guarantees existence of a flow of size $k$ in $N^{\mbox {v}}$. % (and in fact, generalizes Menger's theorem).

\begin{theorem}
\label{th:max-flow min-cut}
In any flow network, the minimal capacity of an $s,t$-cut equals the maximal size of a flow from $s$ to $t$.
\end{theorem}

%states that the minimal size of an $s,t$-cut equals the maximal flow size; hence, there exists a flow of size $k$ in $N^{\mbox {v}}$.
By executing the Ford-Fulkerson algorithm \cite{10.5555/1942094}, we find such a flow $f$ in $N^{\mbox {v}}$; it is guaranteed to be integral, that is the flow $f(e)$ in each edge $e$ of $N^{\mbox {v}}$ is either 0 or 1. The edges assigned with flow 1 constitute $k$ vertex-disjoint paths from $s$ and $t$ in $\vec G^{\mbox {v}}$, which naturally define $k$ vertex-disjoint paths between $s$ and $t$ in $G$, as required. 
(Thus, Theorem~\ref{th:max-flow min-cut} generalizes Menger's theorem.)

The Ford-Fulkerson algorithm (FF) works in iterations, beginning from the zero flow. Each iteration constructs a flow-augmenting path from $s$ to $t$ in the residual network $N^{\mbox {v}}_f$ w.r.t. the current flow $f$ in $N^{\mbox {v}}$, by executing a labeling algorithm (either DFS or BFS) in $N^{\mbox {v}}_f$; after finding such a path, it augments the current flow using that path and updates the residual network. Each iteration of FF increases the flow size in $N^{\mbox {v}}$ by 1. The runtime of each iteration is $O(|E|)$. That is, for any constant $k$, the time of finding a flow of size $k$ in $N^{\mbox {v}}$, and thus of finding $k$ vertex-disjoint paths between $s$ and $t$ in $G$, is also $O(|E|)$.

Handling the generalization 
of Menger's theorem is easy.
If a set $V_s$ is given instead of $s$, then, before constructing $\vec G^{\mbox {v}}$, enhance $G$ by adding an artificial source $\bar s$ with edges to every vertex in $V_s$. Likewise, add edges
from set $V_t$ to new vertex  $\bar t$. After finding $k$ vertex-disjoint paths between $\bar s$ and $\bar t$, we just truncate these paths by removing $\bar s$ and $\bar t$.

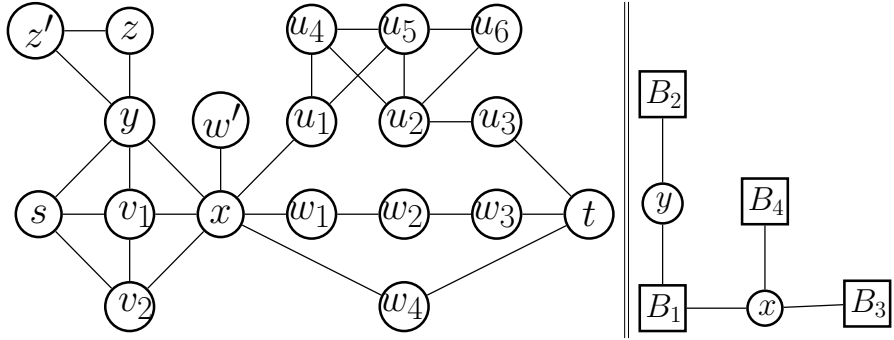
\begin{figure}
    \centering
\resizebox{1.0\columnwidth}{!}{
\begin{tabular}{c||c}
\begin{tikzpicture}[minimum size=10mm,
  node distance=0.9cm and 0.9cm,
  >=stealth,
  bend angle=45,
  auto]
  \tikzstyle{every node}=[font=\Huge]
\node (0) [circle, ultra thick, text centered, text width=0.5cm, draw=black] {$s$};
\node (1) [circle, ultra thick, text centered, text width=0.5cm, draw=black, right=of 0] {$v_1$};
\node (2) [circle, ultra thick, text centered, text width=0.5cm, draw=black, above=of 1] {$y$};
\node (3) [circle, ultra thick, text centered, text width=0.5cm, draw=black, below=of 1] {$v_2$};
\node (4) [circle, ultra thick, text centered, text width=0.5cm, draw=black, right=of 1] {$x$};
\node (100) [circle, ultra thick, text centered, text width=0.5cm, draw=black, above=of 4] {$\! w'$};
\node (5) [circle, ultra thick, text centered, text width=0.5cm, draw=black, above=of 2] {$z$};
\node (16) [circle, ultra thick, text centered, text width=0.5cm, draw=black, left=of 5] {$z'$};
\node (6) [circle, ultra thick, text centered, text width=0.5cm, draw=black, right=of 4] {$\!\!w_1$};
\node (7) [circle, ultra thick, text centered, text width=0.5cm, draw=black, right=of 6] {$\!\!w_2$};
\node (8) [circle, ultra thick, text centered, text width=0.5cm, draw=black, right=of 7] {$\!\!w_3$};
\node (9) [circle, ultra thick, text centered, text width=0.5cm, draw=black, right=of 8] {$t$};
\node (10) [circle, ultra thick, text centered, text width=0.5cm, draw=black, above=of 6] {$\!u_1$};
\node (11) [circle, ultra thick, text centered, text width=0.5cm, draw=black, above=of 7] {$\!u_2$};
\node (15) [circle, ultra thick, text centered, text width=0.5cm, draw=black, above=of 8] {$\!u_3$};
\node (12) [circle, ultra thick, text centered, text width=0.5cm, draw=black, above=of 10] {$\!\! u_4$};
\node (13) [circle, ultra thick, text centered, text width=0.5cm, draw=black, above=of 11] {$\!\! u_5$};
\node (66) [circle, ultra thick, text centered, text width=0.5cm, draw=black, above=of 15] {$\!\! u_6$};
\node (14) [circle, ultra thick, text centered, text width=0.5cm, draw=black, below=of 7] {$\!\!w_4$};
%\node (50) [rectangle,minimum width = 1cm, minimum height= 1cm, ultra thick, text centered, text width=0.75cm, draw=black, below right=of 0, yshift=0.1cm] {$X^*$};

\draw [thick,-,>=stealth] (0) to (1);
\draw [thick,-,>=stealth] (0) to (2);
\draw [thick,-,>=stealth] (0) to (3);
\draw [thick,-,>=stealth] (1) to (2);
\draw [thick,-,>=stealth] (1) to (3);
\draw [thick,-,>=stealth] (1) to (4);
\draw [thick,-,>=stealth] (2) to (4);
\draw [thick,-,>=stealth] (100) to (4);
\draw [thick,-,>=stealth] (2) to (5);
\draw [thick,-,>=stealth] (2) to (16);
\draw [thick,-,>=stealth] (16) to (5);
\draw [thick,-,>=stealth] (3) to (4);
\draw [thick,-,>=stealth] (4) to (6);
\draw [thick,-,>=stealth] (6) to (7);
\draw [thick,-,>=stealth] (66) to (11);
\draw [thick,-,>=stealth] (66) to (13);
\draw [thick,-,>=stealth] (7) to (8);
\draw [thick,-,>=stealth] (8) to (9);
%\draw [thick,-,>=stealth] (10) to (11);
\draw [thick,-,>=stealth] (12) to (11);
\draw [thick,-,>=stealth] (13) to (11);
\draw [thick,-,>=stealth] (12) to (13);
\draw [thick,-,>=stealth] (10) to (12);
\draw [thick,-,>=stealth] (10) to (13);
\draw [thick,-,>=stealth] (4) to (10);
\draw [thick,-,>=stealth] (4) to (14);
\draw [thick,-,>=stealth] (11) to (15);
\draw [thick,-,>=stealth] (9) to (15);
\draw [thick,-,>=stealth] (14) to (9);
\end{tikzpicture}

& 
\begin{tikzpicture}[minimum size=5mm,
  node distance=1.3cm and 1.3cm,
  >=stealth,
  bend angle=45,
  auto]
  \tikzstyle{every node}=[font=\huge]
%\node (0) [circle,draw=green!60, minimum width = 1cm, minimum height= 1cm, ultra thin, text centered, text width=0.01cm, draw=black] {$~$};
\node (2) [rectangle,minimum width = 1cm, minimum height= 1cm, ultra thick, text centered, text width=0.75cm, draw=black] {$B_1$};
\node (1) [circle, ultra thick, text centered, text width=0.3cm, draw=black, above=of 2] {$y$};
\node (3) [rectangle,minimum width = 1cm, minimum height= 1cm, ultra thick, text centered, text width=0.75cm, draw=black, above=of 1, yshift=0.1cm] {$B_2$};
\node (4) [circle, ultra thick, text centered, text width=0.3cm, draw=black, right=of 2] {$x$};
\node (5) [rectangle,minimum width = 1cm, minimum height= 1cm, ultra thick, text centered, text width=0.75cm, draw=black, right=of 4, yshift=0.1cm] {$B_3$};
\node (6) [rectangle,minimum width = 1cm, minimum height= 1cm, ultra thick, text centered, text width=0.75cm, draw=black, above=of 4, yshift=0.1cm] {$B_4$};

\draw [thick,-,>=stealth] (1) to (2);
\draw [thick,-,>=stealth] (1) to (3);
\draw [thick,-,>=stealth] (4) to (2);
\draw [thick,-,>=stealth] (4) to (5);
\draw [thick,-,>=stealth] (4) to (6);
\end{tikzpicture}

\end{tabular}
}

    \caption{A graph (left), and its  block tree (right).}
    \label{fig:example}
\end{figure}

For any non-empty $V' \subseteq V$, we define the \emph{subgraph $G(V')$ induced by $V'$} as $V'$ together with all edges of $G$ between its vertices. % in $G(v)$. 
Vertex $a$ is called an \emph{articulation} (or \emph{separator}) vertex of $G$ if $\{a\}$ is a 1-cut. 
Removing an articulation vertex $a$
from the graph results in two or more induced connected subgraphs $G(V_i)$, $1 \le i \le r$, $r \ge 2$, $\cup_i V_i = V \setminus a$. 
%\YD {Delete: We call the induced connected subgraphs $G(V_i \cup \{a\})$ the \emph{parts} separated by $a$ in $G$.}
Continuing to recursively divide the induced connected subgraphs $G(V_i \cup a)$ by their separating vertices (each is a separating vertex of $G$ as well), the final induced subgraphs are biconnected; they are called \emph{blocks} of $G$.
A block that consists of a single edge
is called \emph{trivial}.
The set of blocks and the \emph{block tree} structure of interleaved articulation vertices and blocks are unique (that is, they do not depend on the order of partitioning of $G$ by the articulation vertices). 
See an example of a graph  and its block tree in Figure~\ref{fig:example}.

%\YD {Let us remove this paragraph, since we never use it in the paper (as far as I recall).
%If you agree, do not remove the text, just comment it, since I plan to use it for another, similar explanation.}
%Every simple path in $G$ can be described as the concatenation of simple paths in the blocks along a chain of blocks in the block tree ``glued'' by the articulation vertices between the consecutive blocks (see Figure~\ref{fig:example}).
%By Menger's theorem, for any three distinct vertices $a$, $b$, $c$ in the same block $B$, there exists a simple path from $a$ to $c$ going via $b$. That is because there exist two disjoint simple paths between $\{a,c\}$ and $b$.
%Therefore, for every chain of blocks and every set of vertices, at most one in each block, there exists a simple path in $G$ containing the chosen vertices.

Observe that any simple cycle in $G$ lies entirely in a single block. Indeed, if a simple path passes from one block to another, via a separator vertex, then it has no way to return.
Therefore, the CEP has a sense only in the case when $s$ and $t$ are in the same block, $B$; in such a case, we can safely discard all other blocks of $G$, thus reducing the problem to the (biconnected) graph $B$.
Due to the equivalence between the PEP and CEP as above, the PEP in a general graph can also be easily reduced to either a PEP or a CEP in a single block of $G^+ = G \cup (s,t)$. 
Henceforth, we thus assume that \emph{ graphs $G$ for the CEP and $G^+$ for the PEP are biconnected}.
%\YD {Based on this assumption, I suggest canceling the notion BGCEP, thus retaining CEP only.}
(\emph{Remark\/}:
Note that $G^+$ is biconnected if and only if the block tree of $G$ has a form of a chain of blocks such that $s$ and $t$ are in the blocks at the two ends of that chain.)

A biconnected graph  can contain
vertex 2-cuts.
Such a vertex pair $\{ a, b\}$ that separates $G$ is henceforth called
a \emph{separation pair} of $G$.
A similar, albeit more complex, sub-division of a biconnected graph $G$
using separation pairs
into triconnected components and auxiliary structural elements is possible %\cite{doi:10.1137/0202012}
\cite{DBLP:journals/algorithmica/BattistaT96}.
%\YD {Let us remove the reference to [12] here, since [12] does not introduce auxiliary structural elements.}
This structure, named an SPQR tree, is described in Section~\ref{sec:SPQR} and 
used extensively in Section~\ref{sec:biconnected}.

%\subsubsection{Edge Cuts and Flows}
%\label{subsec:edgecutsandflows}

In this paper, we also use the concept of an \emph{edge cut},
which is an inclusion-minimal \emph{set of edges} $E'$ of $G$ such that removing the edges in $E'$ disconnects $G$.
Unlike vertex cuts, every edge cut divides $G$ into exactly two connected induced subgraphs; thus, their vertex sets $V',V''$ form a 2-partition of $V$.
Note that  the edges in $E'$ and 
{\em only} these edges connect vertices from $V'$ and $V''$.
Additionally, for any pair of edge cut $E'$ and cycle $L$, the number of edges common to both is even; indeed, traversing $L$ by each edge in $E'$ switches between subgraph $G(V')$ and
subgraph $G(V'')$ or vice versa.
An edge cut is called a $k$-edge cut if $|E'|=k$. To distinguish
such cuts from vertex cuts, we always say ``edge cut'' in full in this paper, leaving ``cut'' as a nickname for a vertex cut.

Any edge cut of the minimal size (cardinality) is called a \emph{minimum} edge cut of $G$.
Two edge cuts 
%$E_1$ and $E_2$ 
of $G$ are called \emph{crossing} 
% if the removal 
if the two corresponding 2-partitions of $V$ subdivide $V$ into four non-empty subsets (otherwise, they subdivide it into three non-empty subsets).
The following theorem is proved in \cite{DKL}:

\begin{theorem}\label{th:crossing}
No minimum edge cuts of odd cardinality cross. 
\end{theorem}

The reader can easily check by hand the special case of this theorem that we use in this paper: if there are two crossing 3-edge cuts, then there exists a
1- or 2-edge cut in the graph.
We also use the following known result on edge and vertex cuts. As we could not find a published proof, we also provide a proof.

\begin{theorem}
\label{l:edge-vertex cuts}
If there is a $k$-edge cut in a graph $G$ with at least $k+2$ vertices, then there exists a vertex cut of cardinality at most $k$ in $G$.
\end{theorem}

\begin{proof}
Let $E'$, $|E'|=k$, be an edge cut partitioning $G$ into $G(V_1)$ and $G(V_2)$. Denote by $\tilde V_i$, $i = 1,2$, the set of end-vertices of edges in $E'$ in $V_i$. Clearly, 
% the sets $V_1$ and $V_2$ are disjoint, and 
$|\tilde V_1|, |\tilde V_2| \le k$.
If for some $i$, the set $V_i \setminus \tilde V_i$ is non-empty, then $\tilde V_i$ is a vertex cut separating $V_i \setminus \tilde V_i$ from $V_{3-i}$, as required. Otherwise, $V_i = \tilde V_i$, $i = 1,2$; we assume so for the rest of the proof.

Suppose that there exist vertices $v_1 \in V_1$ and $v_2 \in V_2$ such that $(v_1,v_2) \not \in E'$.
Then, the following vertex set is a cut separating $v_1$ from $v_2$. For each edge in $E'$ incident on $v_i$, $i = 1,2$, pick its other  end-vertex. For any other edge in $E'$, pick an arbitrary end-vertex thereof. Removing all those vertices implies removing all edges in $E'$, thus separating $v_1$ and $v_2$.
Overall, vertex picks occur $|E'|=k$ times;  since some of the vertices may have been picked more than once, the total vertex cut size is at most $k$, as required.

The only remaining case is where $E' = V_1 \times V_2$. Then, $|V_1| \cdot |V_2| = k$; assume, w.l.o.g., that $|V_1| \le |V_2|$, so that $|V_1| = q \le \sqrt k$. 
Let us show that $|V| = q + k/q \le k+1$, contradicting the condition of the lemma.
If $q=|V_1|=1$, then $|V_2|=k$, thus $|V|=k+1$.
By the assumption $q \le \sqrt k$, the derivative of the function $q+k/q$ is non-positive:
$1 - k/q^2 \le 1 - k/(\sqrt k)^2 = 0$, which suffices.
%
%A simple idea is that either one of the sets of end-vertices of the edges in $E'$ in $V'$ and in $V''$ form a vertex cut as required. However,  . . .
\end{proof}

\begin{figure}
    \centering
\resizebox{1.0\columnwidth}{!}{
\begin{tabular}{c||c}
\begin{tikzpicture}[minimum size=10mm,
  node distance=0.9cm and 0.9cm,
  >=stealth,
  bend angle=45,
  auto]
  \tikzstyle{every node}=[font=\Huge]
\node (4) [circle, ultra thick, text centered, text width=0.5cm, draw=black, right=of 1] {$x$};
\node (6) [circle, ultra thick, text centered, text width=0.5cm, draw=black, right=of 4] {$\!\!w_1$};
\node (7) [circle, ultra thick, text centered, text width=0.5cm, draw=black, right=of 6] {$\!\!w_2$};
\node (8) [circle, ultra thick, text centered, text width=0.5cm, draw=black, right=of 7] {$\!\!w_3$};
\node (9) [circle, ultra thick, text centered, text width=0.5cm, draw=black, right=of 8] {$t$};
\node (10) [circle, ultra thick, text centered, text width=0.5cm, draw=black, above=of 6] {$\!u_1$};
\node (11) [circle, ultra thick, text centered, text width=0.5cm, draw=black, above=of 7] {$\!u_2$};
\node (15) [circle, ultra thick, text centered, text width=0.5cm, draw=black, above=of 8] {$\!u_3$};
\node (12) [circle, ultra thick, text centered, text width=0.5cm, draw=black, above=of 10] {$\!u_4$};
\node (13) [circle, ultra thick, text centered, text width=0.5cm, draw=black, above=of 11] {$\!u_5$};
\node (14) [circle, ultra thick, text centered, text width=0.5cm, draw=black, below=of 7] {$\!\!w_4$};
\node (66) [circle, ultra thick, text centered, text width=0.5cm, draw=black, above=of 15] {$\!\! u_6$};
%\node (50) [rectangle,minimum width = 1cm, minimum height= 1cm, ultra thick, text centered, text width=0.75cm, draw=black, below right=of 0, yshift=0.1cm] {$X^*$};

\draw [thick,-,>=stealth] (4) to (6);
\draw [thick,-,>=stealth] (6) to (7);
\draw [thick,-,>=stealth] (66) to (11);
\draw [thick,-,>=stealth] (66) to (13);
\draw [thick,-,>=stealth] (7) to (8);
\draw [thick,-,>=stealth] (8) to (9);
%\draw [thick,-,>=stealth] (10) to (11);
\draw [thick,-,>=stealth] (12) to (11);
\draw [thick,-,>=stealth] (13) to (11);
\draw [thick,-,>=stealth] (12) to (13);
\draw [thick,-,>=stealth] (10) to (12);
\draw [thick,-,>=stealth] (10) to (13);
\draw [thick,-,>=stealth] (4) to (10);
\draw [thick,-,>=stealth] (4) to (14);
\draw [thick,-,>=stealth] (11) to (15);
\draw [thick,-,>=stealth] (9) to (15);
\draw [thick,-,>=stealth] (14) to (9);
\end{tikzpicture}\\

\begin{tikzpicture}[minimum size=5mm,
  node distance=0.5cm and 0.5cm,
  >=stealth,
  bend angle=45,
  auto]
  \tikzstyle{every node}=[font=\huge]
\node (2) [rectangle,minimum width = 2.5cm, minimum height= 1cm, ultra thick, text centered, text width=0.75cm, draw=black, yshift=0.1cm] {$\!\!\!\! P(x,t)$};
\node (3) [rectangle,minimum width = 5.5cm, minimum height= 1cm, ultra thick, text centered, text width=0.75cm, draw=black, right=of 2, xshift=1.5cm] {$\!\!\!\!\!\!\!\!\!\!\!\!\!\!\!\!\!\! S(x,w_1,w_2,w_3,t)$}
edge [line width=4.0pt] node [anchor=south] {$e_2^{x,t}$} (2);
\node (4) [rectangle,minimum width = 5.5cm, minimum height= 1cm, ultra thick, text centered, text width=0.75cm, draw=black, above=of 2, yshift=1cm] {$\!\!\!\!\!\!\!\!\!\!\!\!\!\!\!\!S(x,u_1,u_2,u_3,t)$}
edge [line width=4.0pt] node [anchor=east] {$e_3^{x,t}$} (2);
\node (7) [rectangle,minimum width = 5.5cm, minimum height= 1cm, ultra thick, text centered, text width=0.75cm, draw=black, above=of 4, yshift=1cm] {$\!\!\!\!\!\!\!\!\!\!\!\!\!\!\!\!R(u_1,u_2,u_4,u_5)$}
edge [line width=4.0pt] node [anchor=east] {$e^{u_1,u_2}$} (4);
\node (6) [rectangle,minimum width = 4.5cm, minimum height= 1cm, ultra thick, text width=0.75cm, draw=black, below=of 2, yshift=-1.1cm] {$\!\!\!\!\!\!\!\!\!\!\!S(x,w_4,t)$}
edge [line width=4.0pt] node [anchor=east] {$e_1^{x,t}$} (2);
\node (5) [rectangle,minimum width = 2.7cm, minimum height= 1cm, ultra thick, text centered, text width=0.75cm, draw=black, right=of 7, xshift=0.3cm] {$\!\!\!\!\!\!\!\!\!P(u_2,u_5)$}
edge [line width=4.0pt] node [anchor=south,yshift=0.5cm] {$e_1^{u_2,u_5}$} (7);
\node (8) [rectangle,minimum width = 3.8cm, minimum height= 1cm, ultra thick, text width=0.75cm, draw=black, right=of 5, xshift=0.3cm] {$\!\!\!\!\!\!\!\!\!\!\!\!\!S(u_2,u_5,u_6)$}
edge [line width=4.0pt] node [anchor=south,yshift=0.5cm] {$e_2^{u_2,u_5}$} (5);

\end{tikzpicture}
& 
\begin{tikzpicture}[minimum size=3mm,
  node distance=0.5cm and 0.5cm,
  >=stealth,
  bend angle=45,
  auto]
  \tikzstyle{every node}=[font=\huge]
\node (0) [circle, ultra thick, text centered, text width=0.5cm, draw=black] {$x$};
\node (1) [circle, ultra thick, text centered, text width=0.5cm, draw=black, right=of 0] {$t$}
  edge [thick, dashed]
    node [anchor=south] {\small $Vir_2$}(0)
  edge [thick, dashed, bend left] 
    node {\small $Vir_1$} (0)
  edge [thick, dashed, bend right] 
  node [anchor=south] {\small $Vir_3$} (0);
  
\node (2) [circle, ultra thick, text centered, text width=0.5cm, draw=black,below=of 0] {$x$};
\node (3) [circle, ultra thick, text centered, text width=0.5cm, draw=black,below=of 1] {$t$}
   edge [thick, dashed]
    node [anchor=south] {\small $Vir_1$}(2); 
\node (4) [circle, ultra thick, text centered, text width=0.5cm, draw=black,below=of 2] {$w_4$}
   edge [thick] (2)
   edge [thick] (3);

\node (5) [circle, ultra thick, text centered, text width=0.5cm, draw=black,above=of 0] {$x$};
\node (6) [circle, ultra thick, text centered, text width=0.5cm, draw=black,right=of 5] {$t$}
   edge [thick, dashed]
    node [anchor=north] {\small $Vir_3$}(5); 
\node (8) [circle, ultra thick, text centered, text width=0.5cm, draw=black,above=of 6] {$u_3$}
   edge [thick] (6);

\node (9) [circle, ultra thick, text centered, text width=0.5cm, draw=black,above=of 8] {$u_2$}
   edge [thick] (8)
; 
\node (7) [circle, ultra thick, text centered, text width=0.5cm, draw=black,left=of 9] {$u_1$}
edge[thick] (5)
   edge [thick, dashed]
    node [anchor=south,yshift=0.2cm] {\small $Vir_6$}(9);

\node (12) [circle, ultra thick, text centered, text width=0.5cm, draw=black, above=of 7] {$u_1$};
\node (13) [circle, ultra thick, text centered, text width=0.5cm, draw=black, above=of 9] {$u_2$}
  edge [thick, dashed] 
  node [anchor=north,yshift=-0.2cm] {\small $Vir_6$} (12);
\node (14) [circle, ultra thick, text centered, text width=0.5cm, draw=black, above=of 12] {$u_4$}
  edge [thick] (12)
  edge [thick] (13);
\node (15) [circle, ultra thick, text centered, text width=0.5cm, draw=black, above=of 13] {$u_5$}
  edge [thick] (12)
  edge [thick] (14)
  edge [thick, dashed] node [yshift=-0.1cm] {\small $Vir_5$} (13);    
  
\node (11) [circle, ultra thick, text centered, text width=0.5cm, draw=black, right=of 13] {$u_2$};
\node (10) [circle, ultra thick, text centered, text width=0.5cm, draw=black, above=of 11] {$u_5$}
  edge [thick] (11)
  edge [thick, dashed, bend left] 
    node [yshift=-0.1cm] {\small $Vir_4$} (11)
  edge [thick, dashed, bend right] 
  node [anchor=east,xshift=0.1cm,yshift=0.2cm] {\small $Vir_5$} (11);

\node (21) [circle, ultra thick, text centered, text width=0.5cm, draw=black, right=of 11] {$u_2$};
\node (20) [circle, ultra thick, text centered, text width=0.5cm, draw=black, above=of 21] {$u_5$}
  edge [thick, dashed] node [anchor=east, yshift=0.2cm] {\small $Vir_4$} (21);
\node (22) [circle, ultra thick, text centered, text width=0.5cm, draw=black, right=of 20] {$u_6$}
  edge [thick] (21)
  edge [thick] (20);
    
\node (16) [circle, ultra thick, text centered, text width=0.5cm, draw=black,right=of 3] {$x$};
\node (17) [circle, ultra thick, text centered, text width=0.5cm, draw=black,right=of 16] {$t$}
   edge [thick, dashed]
    node [anchor=north] {\small $Vir_2$}(16); 
\node (18) [circle, ultra thick, text centered, text width=0.5cm, draw=black,above=of 17] {$\!w_3$}
   edge [thick] (17);

\node (19) [circle, ultra thick, text centered, text width=0.5cm, draw=black,above=of 18] {$\!w_2$}
   edge [thick] (18)
; 
\node (20) [circle, ultra thick, text centered, text width=0.5cm, draw=black,left=of 19] {$\!w_1$}
   edge[thick] (16)
   edge [thick] (19);
\end{tikzpicture}

\end{tabular}
}
    \caption{
    A biconnected graph (top left), its SPQR tree structure (bottom left), and details of
    its components (right).
    On the right,
    for each structural edge (bottom left, thick line),  the two copies of the corresponding virtual edge (right, dashed) are denoted by $Vir_j$ with the same $j$.}
    %(These indices are assigned arbitrarily for clarity,
    %with no relation to their indexing in the definitions.)
    %\YD {Now. with the heavy indexing of virtual edges in the definitions, I see no need of the text in the parentheses. Especially, of its second half, but also for its first half.}}
    
    \label{fig:exampleSPQR}
\end{figure}
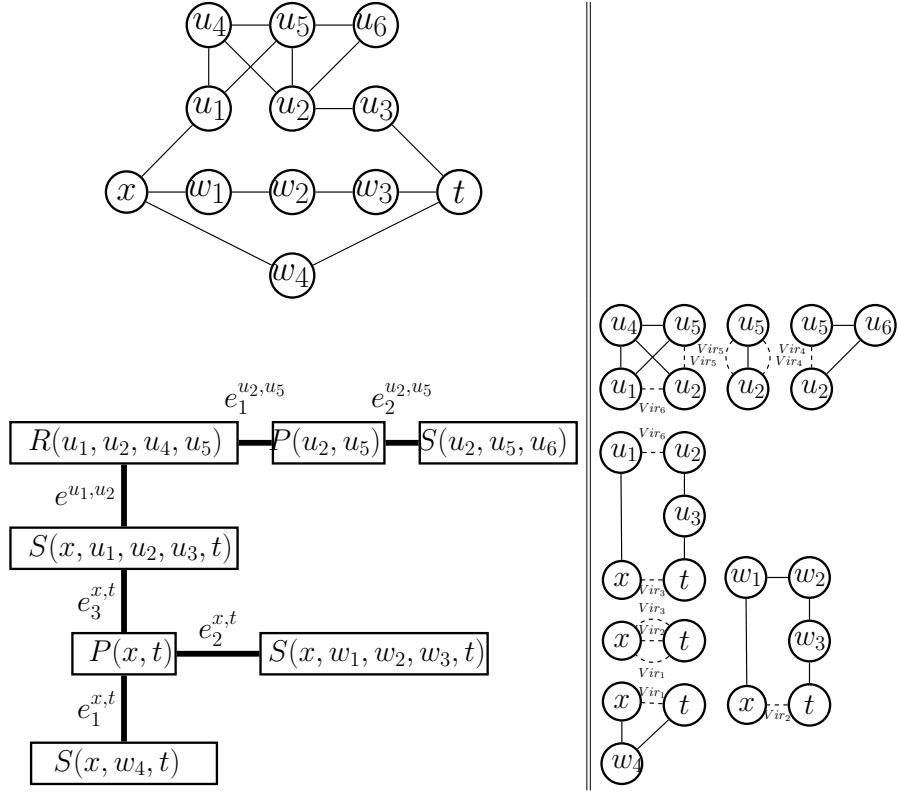

%***************************************
\subsection{Triconnected Components and SPQR trees}
\label{sec:SPQR}

%A similar, albeit more complex, sub-division of a \emph{2-connected graph $G$} into 3-connected components and auxiliary structural elements is possible \cite{doi:10.1137/0202012}.
A neat representation of the sub-division of a biconnected graph $G$ into triconnected components and auxiliary structural elements 
%\YD {Delete: introduced in  \cite{doi:10.1137/0202012}}
%\YD {My reason is that we refer to the "neat representation", not just to the division of a graph into its triconnected components introduced in \cite{doi:10.1137/0202012} (BTW, that division was introduced before by Mac Lane, Saunders (1937)). See Wikipedia on SPQR trees.}
%\YDT {for planar graphs}
%\YD {I do not so like this preference of Hopcroft and Tarjan before Tamassia and Di Battista. So, if you like to keep it, I suggest the following: a) the addition on planar graphs, see above, and 2) addition of Mac Lane, Saunders (1937) (see Wikipedia on SPQR trees) together with Hopcroft and Tarjan.} 
is called the SPQR tree of $G$ \cite{DBLP:journals/algorithmica/BattistaT96}; such a tree can be constructed in a linear time $O(|E|)$ \cite{SPQRtrees}.
The SPQR tree of $G$, ${\cal T} = {\cal T}(G)$ 
consists of nodes, each being a graph called a ``triconnected component''\footnote{This is the standard historically used term, even though not all these components are triconnected graphs.}
 of $G$ 
%or ``skeleton graph'' in \cite{DBLP:journals/algorithmica/BattistaT96,SPQRtrees}
%\YD {Where is the notion "skeleton graph" from? I do not recall it.}
%\ES{skeleton (mu): appears in Battista and Tamassia.}
%\YD {Maybe, we would omit this notion? I do not see what it adds, and we do not use it.}
(henceforth called simply a \emph{component}) of
one of 4 types: S (series),
P (parallel), R (rigid), and Q (representing a single edge). Henceforth, we use ``X node/component'' as shorthand for ``node/component of type X''.
%\YD {Let us choose what we write *everywhere*: either X-component/node, or X component/node. As for now, there are more cases with than without "-". Though, we can decide as we wish.}

Each component is a part of graph $G$ augmented with auxiliary {\em virtual edges}.
The components (nodes) are connected in tree $\cal T$ by {\em structural edges}.
Figure \ref{fig:exampleSPQR} (bottom left) shows
the gross structure of the SPQR tree of the example graph in Figure \ref{fig:exampleSPQR} (top left): rectangles represent components, with structural edges between them;
the components are shown in  Figure \ref{fig:exampleSPQR} (right). 
As we use SPQR trees extensively in this paper, we further delve into their details below.
There exist different
definitions of SPQR trees, which vary in details; we adopt
the definition of \cite{SPQRtrees} that omits the
Q nodes.
%\YDT {%We mainly follow the 
The main difference of our definition from that of \cite{SPQRtrees} is that our SPQR tree is \emph{unrooted}.\footnote{After \cite{DBLP:journals/algorithmica/BattistaT96,SPQRtrees} were published, R. Tamassia agreed that the unrooted form is preferable for SPQR trees (a personal communication to the first author).
}
%\emph{Remark\/}:
As we are concerned here with intuitive understanding rather than with how to construct the SPQR tree, we provide a description %of the tree %and its properties 
that, unlike prior work, is not
meant to lead to an efficient algorithm of creating the tree.
%\YD{I suggest this change, since, to my mind, our definition DOES desrcibe a recursive process. Argue, if you feel that I am not right.}
%See \cite{SPQRtrees} on how to construct the SPQR tree of a graph  efficiently.
%\YD {The last sentence doubles the information given above.}

Consider the graph shown in Figure \ref{fig:exampleSPQR} (top left). 
There is a must-include path from $w_4$ to $x$ via $u_1$ and $u_6$: $(w_4,t,u_3,u_2,u_6,u_5,u_1,x)$, but must-include paths from $w_4$ to $x$ through $u_1$ and $w_2$
and from $u_2$ to $u_4$ through $w_4$ and $u_6$
do not exist.
%There is a must-include path from $u_1$ to $u_4$ via $w_4$ and $u_6$: $(u_1,x,w_4,t,u_3,u_2,u_6,u_5,u_4)$, but a must-include path from $u_2$ to $u_4$ through the same $w_4$ and $u_6$ does not exist.
We show the reasons for these query answers in terms of the SPQR tree of the graph, at the end of Section~\ref{sec:biconnected}.

To define the SPQR tree, we begin with auxiliary definitions.
%We introduce them for multi-graphs (where there may be two or more parallel edges between a pair of vertices), since multi-graphs can be produced as intermediate objects.
%If there are at least two edges between two vertices, then the bunch of all edges between them is called a \emph{bond}.
We assume that %there are at least three vertices in 
the given graph $G$ is not a single edge. % and not a pair of parallel edges.
%A \emph{split pair} of $G$ is either a separation pair or the pair of end-vertices of a bond if there are at least three vertices in the graph.
For every separator (2-cut) $\{a,b\}$ of $G$,
the graph $G \setminus \{ a, b\}$ consists of $r^{a,b} \ge 2$ 
maximal connected induced subgraphs, which we denote by $G(V_i^{a,b})$, $1 \leq i \leq r^{a,b}$.
%(????), or a single such graph $G(V_1^{a,b})$ and the bunch $G(V_0^{a,b})$ of parallel edges between $a$ and $b$, $V_0^{a,b}=\{a,b\}$. 
The \emph{split subgraphs w.r.t. $\{a,b\}$} are %$G_i^{a,b} = 
$G(V_i^{a,b}\cup \{a,b\}) \setminus (a,b)$,
%with all edges between $a$ and $b$ removed,
$1 \leq i \leq r^{a,b}$. %, and $G(\{a,b\})$, if there is at least one edge $(a,b)$ in $G$.
%(Note that, by definition, there are at least two split graphs w.r.t. any separation pair.)
%Denote by $G^{a,b}_i$ the \emph{extended} induced subgraph $G(V_i^{a,b}\cup \{a,b\})\setminus\{ (a,b)\}$.
%\YD {Why do we omit the superscript \{a,b\} at $G_i$ and $G'_i$? These objects are not unique in $G$, and I suppose that we might miss those superscripts in our text.}
    \remove{
If there are exactly two split subgraphs w.r.t. $\{a,b\}$,
we call $\{a,b\}$ a \emph{binary} separator, otherwise,
for reasons made obvious below, it is called a \emph{P-separator}. 
}
For example, in the graph of Figure \ref{fig:exampleSPQR} (top left), there are two split subgraphs w.r.t. 2-cuts $\{u_1, u_2\}$, $\{u_2,u_5\}$, and $\{x, w_3\}$
and three split subgraphs w.r.t. $\{x, t\}$.
    \remove{
We call separator pair $\{a,b\}$ of $G$ a
\emph{proper} separator (of $G$)
if no separator pair $\{c, d\}$ separates $a$ from $b$ in $G$.
In Figure \ref{fig:exampleSPQR}
(top left), $\{u_1, u_2\}$ and $\{x, t\}$
are proper separators. $\{x, w_3\}$
is not a proper separator of $G$, because
$\{ w_1, t\}$ separates $x$ from $w_3$.
Note that every P-separator $\{a,b\}$ is also proper, because
by definition there exist three internally disjoint paths between $a$ and $b$ (one in each $G_i$ and $(a,b)$, if it exists), which cannot all be cut by only two vertices $c, d$.
}

%We define the following separation pair \emph{partition operators} on a biconnected graph, used to construct the SPQR tree, as follows:

We define the SPQR tree ${\cal T}(G)$ of any given biconnected graph $G$ recursively,
with the base cases as follows:

\begin{description}
 %   \item \YD {deleted on bond}
% If $G$ is a bond of at least three edges between two vertices, ${\cal T}(G)$ consists of a single P component $G$.
    \item[Cycle (S node):] If $G$ is a cycle, ${\cal T}(G)$ consists of a single S component $G$.
    \item[Triconnected (R node):] If there is no 2-cut in $G$, and $G$ is not a triangle, then $G$ is a triconnected graph on at least four vertices. ${\cal T}(G)$ consists of a single R component $G$.
\end{description}

\begin{figure}[h]
\centering
\includegraphics[scale=0.22]{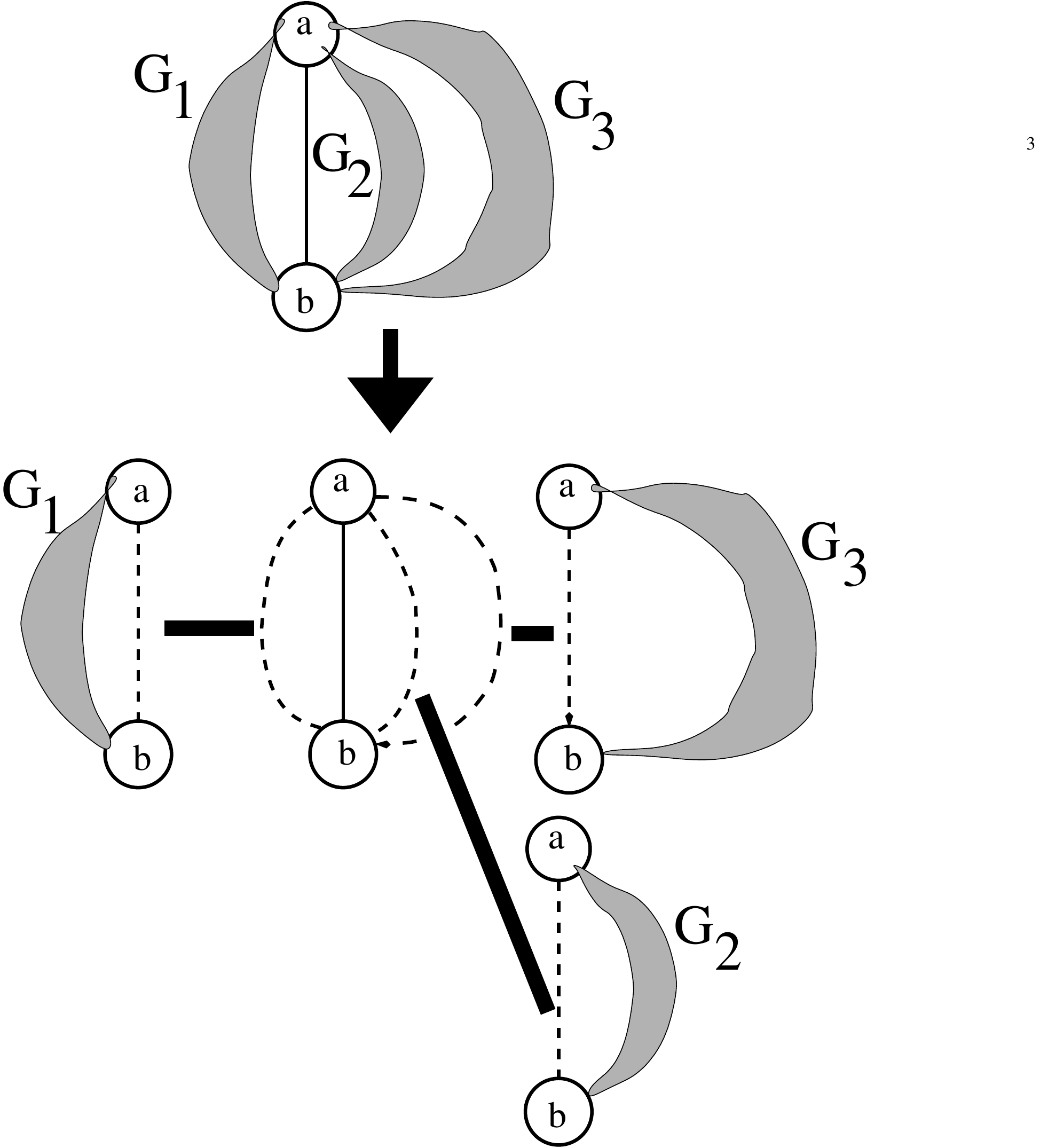}
\includegraphics[scale=0.3]{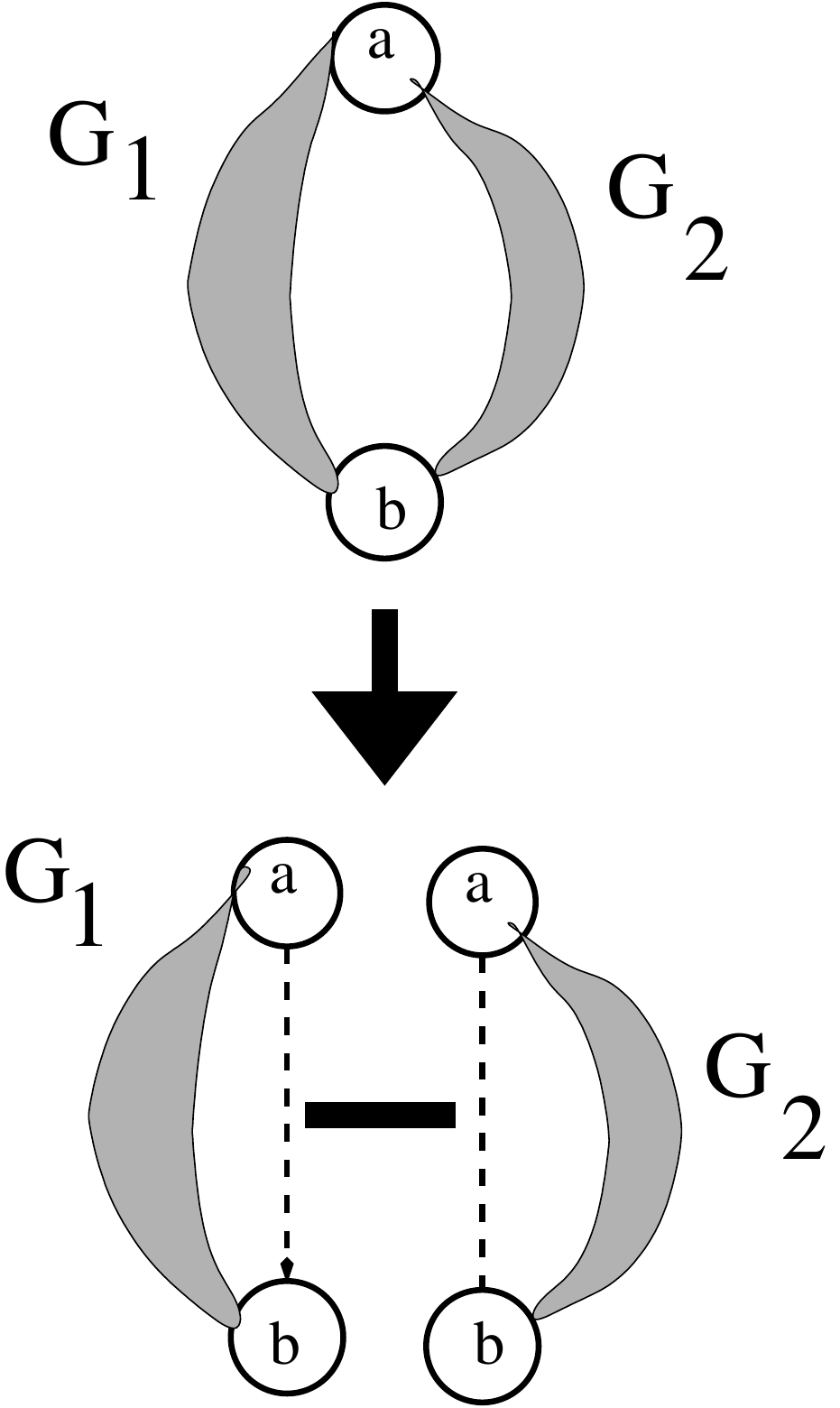}
\hspace{1cm}
\includegraphics[scale=0.23]{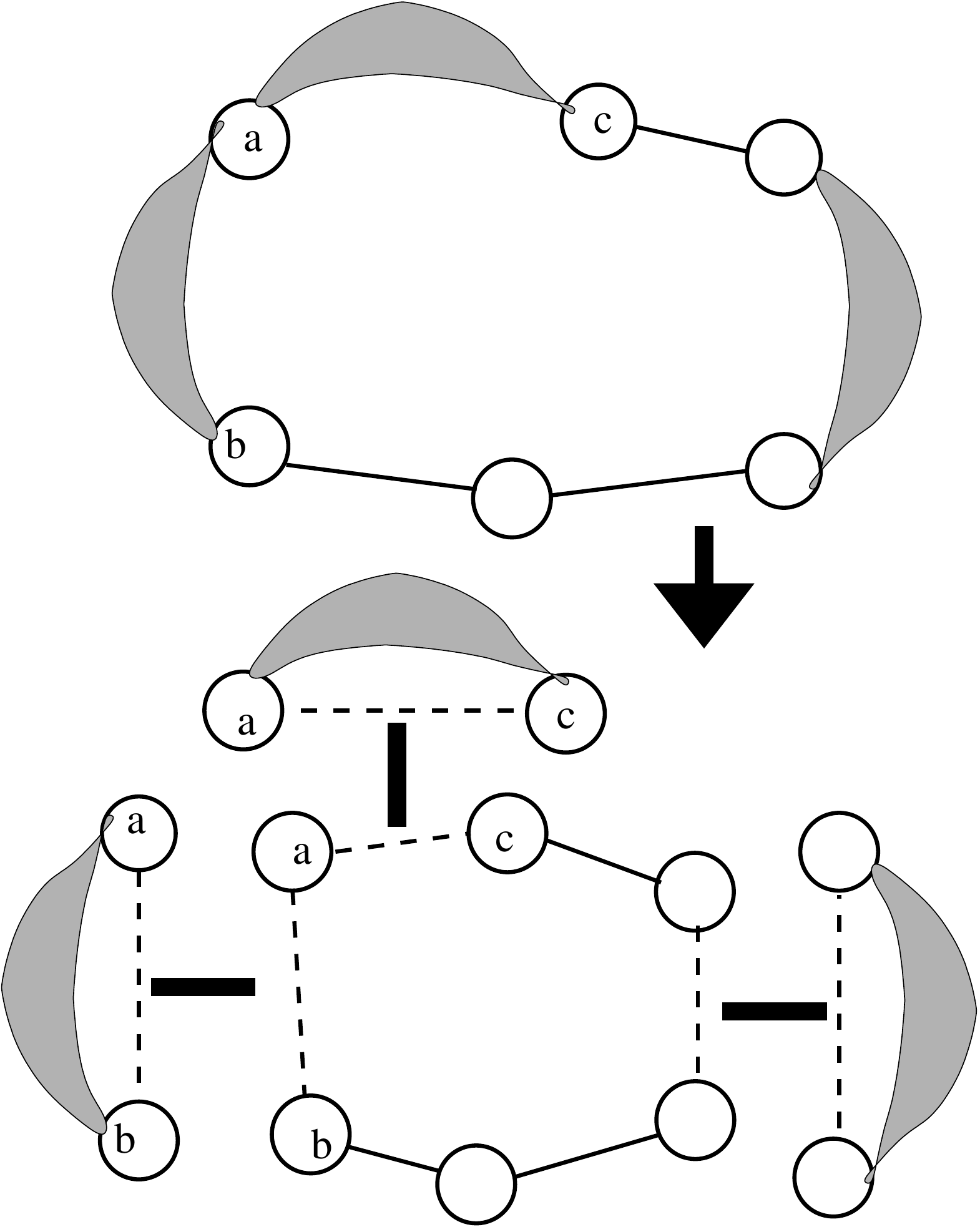}
\caption{Parallel case (left), binary case (middle), and S node construction (right)}
\label{fig:SPQRcases}
\end{figure}

Otherwise, if $G$ has 2-cuts and is not a cycle, % $\{a,b\}$, 
then ${\cal T}(G)$ is defined by the following recursion. See illustration in Figure~\ref{fig:SPQRcases} (left and middle).

\begin{description}
\item[Parallel case] If $\{a,b\}$ is a 2-cut with either $r^{a,b} \ge 3$ or with $(a,b) \in E$ (or both), we create new 
\emph{virtual edges} $e^{a,b}_i=(a,b)$, 
%denoted $Vir_i^{a,b}$, 
    %\YDT {and add it to}
%for each split subgraph $G(V_i^{a,b})$, 
$1 \leq i \leq r^{a,b}$.
We also create a \emph{P component} $C^P$ %, i.e. an auxiliary part $P(a,b)$:
as the multi-graph %consisting of the two vertices 
on two vertices $a, b$ with (at least three) edges between them: all virtual edges $e^{a,b}_i$ as well as edge $(a,b)$, if $(a,b)\in E$.
The tree ${\cal T}(G)$ consists of node $C^P$ and $r^{a,b}$ trees ${\cal T}(G(V_i^{a,b}) \cup e^{a,b}_i)$, each connected to $C^P$ by a \emph{structural edge} from the component in ${\cal T}(G)$ that contains $e^{a,b}_i$.
%(See Figure~\ref{recursion} (left).
    \remove{
The resulting parts are 
$G^\prime_i=G_i \cup Vir_i^{a, b}$ for $r\geq i\geq 1$. 
Each part w.r.t.\ $\{a,b\}$ is connected to $P(a,b)$ by a structural edge, corresponding
to the respective virtual edge $Vir_i^{a,b}$.
As above, we say that $Vir_i^{a,b}$ \emph{represents}
$G_i$ and all of its elements except for $a, b$ in $P(a,b)$, and also
\emph{represents} $G \setminus G_i$ and all of its elements %except for $a, b$ 
in $G^\prime_i$.
}
%\ES{I think $G \setminus G_i$ does not include $a,b$, so we don't need to say "except for a,b" at the end (though still technically correct)}

\item[Binary case] Otherwise, if $\{a,b\}$ is a 2-cut with exactly two split subgraphs w.r.t. $\{a,b\}$ and at least one of them \emph{is biconnected},
%(ES: Simply change this to: at least one is biconnected, and SEE BELOW)
then we create a new \emph{virtual edge} $e^{a,b}=(a,b)$. The tree ${\cal T}(G)$ consists of trees ${\cal T}(G(V_1^{a,b}) \cup e^{a,b})$ and ${\cal T}(G(V_2^{a,b}) \cup e^{a,b})$ connected by a \emph{structural edge} between the components in those trees containing $e^{a,b}$.
\remove{
\item [Series case] (ES: now can drop this entire case!!!
That is because if there are no split pairs
that are parallel or the (modified) binary then the graph is ALREADY a bond, triconnected, or a cycle!)

It holds if there are exactly two split subgraphs $G_1$ and $G_2$ w.r.t. $\{a,b\}$ and at least one of them is not biconnected.
Then, if $G_j$, $j=1,2$, is not biconnected, then it is a chain ${\cal B}_j$ of at least two blocks between $a$ and $b$.
Summarizing, we have a circular chain of at least three blocks of $G_1$ and $G_2$ separated by $a$, $b$, and the articulation vertices of ${\cal B}_1$ and ${\cal B}_2$.
Let us denote by $B_i$, $1 \leq i \leq \ell^{a, b}$ the \emph{non-trivial} blocks among them; note that there is at least one such block, since $G$ is not a cycle.
We create a new \emph{virtual edge} $e_i=(c,d)$ 
%denoted $Vir_i^{a,b}$, 
    %\YDT {and add it to}
for each block $B_i$, $1 \leq i \leq \ell^{a,b}$, where $c$ and $d$ are the two vertices separating $B_i$ in the above circular chain.
We also create an \emph{S component} $C^S$ %, i.e. an auxiliary part $P(a,b)$:
as the cycle %consisting of the two vertices 
consisting of all the above virtual edges $e_i$ and all edges that are the trivial blocks in ${\cal B}_1$ and ${\cal B}_2$, in the order of the above circular chain.
The tree ${\cal T}(G)$ consists of node $C^S$ and $\ell^{a,b}$ trees ${\cal T}(B_i \cup e_i)$, each connected to $C^S$ by a \emph{structural edge} from the component in it containing $e_i$.
}
\end{description}

For example, consider a recursive composition of the SPQR tree of the graph $G$ in Figure \ref{fig:exampleSPQR}
(top left). %The first one is as follows.
For illustration, follow the other parts of the figure.
Henceforth, we call edges of $G$ in components ``real edges''.
By splitting $G$ w.r.t. $\{u_1, u_2\}$, we get two split subgraphs (above and below $\{u_1, u_2\}$ in the figure).
Since the top subgraph is biconnected, %but the bottom is not: its blocks are three edges and the non-trivial bottom block on six vertices separated by vertices $x$ and $t$. 
we are in the binary case
of the recursion.
We create a new virtual edge $(u_1,u_2)$ (denoted $Vir_6$ in the figure) and 
%According to the description of the case, we create S component $S(x,u_1,u_2,u_3,t)$ containing three real and two virtual edges: $(u_1,u_2)$ (marked by $Vir_6$ in the figure) and $(x,t)$ (marked by $Vir_3$ there).
%We then 
continue to construct the SPQR trees of the top and bottom split subgraphs with virtual edge $(u_1,u_2)$ added to each (henceforth, the ``top'' and ``bottom'' graphs, resp.). %and of the bottom block on six vertices with added virtual edge $(x,t)$ (henceforth, ``the bottom graph'').
%(Note that splitting the original graph w.r.t. any other pair of vertices among $x,u_1,u_2,u_3,t$ (e.g., w.r.t. $\{x, u_3\}$) instead of $\{u_1, u_2\}$, would lead to the same circular chain of blocks and thus to exactly the same recursive composition of ${\cal T}(G)$.)

By splitting the top graph by
$\{u_2,u_5\}$ (parallel case), we get two split subgraphs: 
the ``left'' and ``right'' ones.
The arising P component $P(u_2,u_5)$ consists of three edges $(u_2,u_5)$: the left virtual (denoted $Vir_5$), the real, and the right virtual (denoted $Vir_4$) ones.
The SPQR trees of the left and right graphs with added virtual edges are of a single node each (base cases of the recursion): R component $R(u_1,u_2,u_4,u_5)$ (a complete graph on four vertices) and S component $S(u_2,u_5,u_6)$ (a cycle on three vertices), respectively. 
These nodes are connected to $P(u_2,u_5)$ each by the structural 
edges $e_1^{u_2,u_5}$ and $e_2^{u_2,u_5}$ corresponding to the left and right virtual edges $(u_2,u_5)$, resp., thus forming the entire SPQR tree of the top graph. % w.r.t. $\{u_1, u_2\}$.

By splitting the bottom graph w.r.t. $\{x, t\}$ (parallel case), we get P component $P(x,t)$ with three virtual edges (denoted $Vir_3$, $Vir_2$, and $Vir_1$), and three split subgraphs. %, which form a cycle each after adding the related virtual edge. 
These three subgraphs with added respective virtual edges $(x,t)$ have a single-node SPQR tree each: S components $S(x,u_1,u_2,u_3,t)$, $S(x,w_1,w_2,w_3,t)$, and $S(x,w_4,t)$. Each of these S nodes is connected to $P(x,t)$ by the structural 
edges $e_3^{x,t}$, $e_2^{x,t}$, and $e_1^{x,t}$ corresponding to the virtual edges $Vir_3$, $Vir_2$, and $Vir_1$, resp., thus forming the SPQR tree of the bottom graph. 
Finally, the SPQR trees of the top and bottom graphs are connected by the structural edge 
$e^{u_1,u_2}$ between their components $R(u_1,u_2,u_4,u_5)$ and $S(x,u_1,u_2,u_3,t)$ containing virtual edge $(u_1,u_2)$ (denoted $Vir_6$ in the figure), thus arriving at the entire SPQR tree ${\cal T}(G)$.

\remove{
Note that different ways of recursion can be chosen.
For example, let us look on the recursive composition of the SPQR tree of the \emoh{bottom graph} as above which begins from splitting it w.r.t. $\{w_1, w_3\}$ (series case). It separates S component $S(x,w_1,w_2,w_3,t)$ and refers to the SPQR tree of the graph on vertices $x,w_4,t$ with two \emph{parallel} edges $(x,t)$: existing $Vir_3$ and newly created $Vir_2$.
The only way of its recursive composition is w.r.t. $\{x, t\}$ (series case), which creates S component $S(x,w_4,t)$.
The next recursion step is the recursion leaf on the bond on $(x,t)$, creating P component $P(x,t)$.
Finally, $S(x,w_4,t)$ and $P(x,t)$ are connected by a structural edge, resulting in the same SPQR tree of the bottom graph.
Overall, the same SPQR tree ${\cal T}(G)$ is built.
}

We now show that our definition of the SPQR tree is essentially the same as in \cite{SPQRtrees}, with only technical differences. %, let us delve a bit into the history of the topic.
The classic papers of McLane \cite{10.1215/S0012-7094-37-00336-3} and Hopcroft and Tarjan \cite{doi:10.1137/0202012} subdivided the given biconnected graph into its triconnected components via small steps,  splitting by 2-cuts \emph{into two} in all cases. 
Moreover, some of the splits made are canceled in a post-processing phase.
The founders of SPQR trees \cite{DBLP:journals/algorithmica/BattistaT96,SPQRtrees} introduced the tree structure on the components, while making the construction steps bigger. Each step adds a new component in its final form, thus eliminating the need for post-processing.
The entire construction propagates from the (arbitrarily chosen) root to the leaves of the SPQR tree.
We essentially retained their definition, with two technical changes (partially ``reviving'' the approach of \cite{10.1215/S0012-7094-37-00336-3, doi:10.1137/0202012}).
First, we canceled rooting and propagation, replacing it by a recursion with an arbitrary order of splitting by (carefully chosen) 2-cuts. Second, we replaced their %(rather heavy) 
one-step constructions of S and R components by smaller steps, %(similar to those of \cite{McLane} and Hopcroft and Tarjan \cite{doi:10.1137/0202012}
so that such components get their final form as the base cases of the recursion. % (with no need of post-processing).

Let us compare our approaches using the example in Figure~\ref{fig:SPQRcases} (right). When the algorithm of \cite{SPQRtrees} arrives at 2-cut $\{a, b\}$ when propagating from the grey area with $a$ and $b$ at its boundary, it reveals the entire chain of blocks between $a$ and $b$ in the rest of the graph. Then, it creates the S component along the circular chain of blocks, replacing each one of the non-trivial blocks
therein (the grey areas) by a virtual edge. The further propagation is to those blocks, except for that between $a$ and $b$, with the respective virtual edge added to each.
We do the same by separating these blocks \emph{one by one}, using the 2-cuts on their boundaries, so that the S component as above is finally revealed as a base case of the recursion. Thus, we arrive at the same SPQR tree structure.
A similar technical difference occurs for R components.

As our definition and that of \cite{SPQRtrees}
are equivalent,
we can use the known properties of SPQR trees, listed in what follows.
For any biconnected graph $G=(V,E)$, the tree ${\cal T}(G)$ defined as above is unique (independent of choosing its recursive composition) and coincides with that defined in \cite{SPQRtrees}.
The total number of elements in ${\cal T}(G)$, including all structural edges, components, and all the elements within components: vertices, edges of $G$, and virtual edges,
is $O(|E|)$ \cite{SPQRtrees}.
By construction, SPQR trees have no adjacent pairs of S nodes 
%\YD {Note that this "no adjacent pair of S-nodes" is not a trivial consequence of the construction, as far as I see now. Do we need to base this? (preferably, no)} \ES{Easily follows from stuff in the appendix, we can decide how to say this later.  It was actually trivial from Yefim's alternate definition of basic cuts.}, 
and no adjacent pairs of P nodes.
%\YD {Maybe, the right place of these two properties is among Properties listed below?}
%For illustration, see Figure~\ref{fig:exampleSPQR}.
%In this example, all virtual edges are due to P-nodes,
%except for $Vir_4$ between
%$u_1$ and $u_2$ which is shared directly
%between the R-node and the S-node having the cycle $x,u_1,u_2,u_3,x$.
Note that any component that is a triangle (such as that on $u_2,u_5, u_6$ in our example) is an S-node, rather than an R-node (which should contain at least four vertices, by definition).

\remove{
\begin{theorem}[\cite{SPQRtrees}]
\label{pr:SPQR-size-time}
Given a biconected graph $G=(V,E)$,
the total number of elements in $\cal T$,
the SPQR tree of $G$
(including all structural edges, components, and all
the elements within components: vertices, real edges, and virtual edges)
is $O(|E|)$. The tree $\cal T$ can be constructed in time $O(|E|)$.
\end{theorem}
}

Recall that each element of any component either is real (taken from $G$) % (in which case we call it is ``real''), 
or is a virtual edge.
%By the construction of $\cal T$, 
Any element of $G$ is an element of at least one component.
Every edge of $G$ is always a real edge of exactly one component, while any vertex of $G$ is a vertex of more than one component if and only if it is a member of a separation pair of $G$.
By construction, any virtual edge has exactly two copies, in the components connected by the structural edge corresponding to that virtual edge.

Let $e=(a,b)$ be a virtual edge in a component $C$ and %$e_{\cal T}(e) = 
$(C,C')$ be the structural edge in $\cal T$ corresponding to $e$.
We denote by $B(C,e)$ the sub-tree of $\cal T$ hanging on $(C,C')$ from $C$.
Let $G(C,e)$ be the subgraph of $G$ consisting 
of all vertices and edges of $G$ appearing in the nodes (components) of $B(C,e)$ as their real elements.
By construction, SPQR trees have the following property.

\begin{property}
\label{prop:sub-tree}
For any virtual edge $e=(a,b)$ and the structural edge $(C,C')$ corresponding to it,
the subgraphs $G(C,e)$ and $G(C',e)$ cover together the entire $G$ and have only vertices $a,b$ and no edges in common.
\end{property}

%The essential property of 
Tree ${\cal T}={\cal T}(G)$ \emph{models  all 2-cuts of $G$} in the following way.
Consider any P node $C^P$, which is a bond on vertices $a,b$. The pair $\{a, b\}$ is a 2-cut of $G$, so that 
the $r^{a,b}$ maximal connected induced subgraphs $G(V_i^{a,b})$ are the graphs $G(C^P,e^{a,b}_i)$ for all virtual edges $e^{a,b}_i$ in $C^P$.
Consider any structural edge $(C,C')$, corresponding to virtual edge $e=(a,b)$, such that neither $C$ nor $C'$ is a P node. The pair $\{a, b\}$ is a 2-cut of $G$, so that there are exactly two split graphs $G(C,e)$ and $G(C',e)$ w.r.t. it.
Consider any two non-adjacent vertices $a,b$ of an S component $C^S$. %, such that $(a,b) \notin E$. 
Let $a$ and $b$ break cycle $C^S$ into paths $P_1$ and $P_2$. The pair $\{a, b\}$ is a 2-cut of $G$, so that the two split graphs w.r.t. it are the graphs $\cup_{e \in P_i} G(C^S,e) \cup (P_i \cap G)$, $i=1,2$.
There are no other 2-cuts in $G$.

\begin{property}
\label{prop:biconnected}
$G(C,(a,b))\cup (a,b)$ is biconnected.
\end{property}

Hence, there exists a simple path between $a$ and $b$ in $G(C,(a,b))$.

\begin{property}\label{prop:pathWvertex}
Given any element $x \in G(C,e)$, there exists a simple path $P$ between $a$ and $b$ in $G(C,e)$ such that $x\in P$.
\end{property}

The latter property follows from Property~\ref{prop:biconnected} due to the extension of Menger's theorem applied to $\{a,b\}$ and either $x$, if $x$ is a vertex, or $\{u,v\}$, if $x$ is an edge $(u,v)$, in $G(C,(a,b))\cup (a,b)$.
%\YD {I suggest replacing the parentheses by commas in the previous sentence.}
Note that by Property~\ref{prop:sub-tree}, the paths between $a$ and $b$ %in Properties~\ref{prop:biconnected} and \ref{prop:pathWvertex} 
as above
do not use any elements of $C$ other than %their end-vertices. 
$a$ and $b$.

In this paper, we use the following notion of \emph{representation} of elements of $G$ in components.
Let $C$ be an arbitrary component of $\cal T$.
%We define the notion of \emph{representation} of the elements of $G$ in $C$ as follows.
%Every element of $G$ is \emph{represented} in $C$
We say that any real element $x$ of $C$ represents itself in $C$. %; formally: $r_C(x)=x$.
%All other elements of $G$ are represented in $C$ as follows.
For any virtual edge $e=(a,b) \in C$,
%and $e_{\cal T}(e)$ be the structural edge in $\cal T$ corresponding to $e$.
%We denote by $B(C,e)$ the sub-tree of $\cal T$ hanging on $e_{\cal T}(e)$ from $C$.
%Let $G(C,e)$ be the subgraph of $G$ consisting of all vertices and edges of $G$ appearing in the nodes (components) of $B(C,e)$ as their real elements.
all elements of $G(C,e)$ except for $a$ and $b$ are represented in $C$ by $e$. %; formally, $r_C(x)=e$.
Thus, every element $x$ of $G$ is represented exactly once in $C$; we denote the unique object representing it---either $x$ itself or a virtual edge---by $r_C(x)$.

%*****************************************
\subsection{PEP and CEP: Known Properties}\label{sec:PEPknown}

If the graph 
%(or a component thereof in which all vertices of interest reside) 
happens to be triconnected, the answer to
the PEP is always positive.
%\YD {I see the complicated previous sentence be extra: (a) the reader would not be able to understand/make use of it, (b) in the general case, we do not have PEP in components, just CEP (PEP is a special boundary case).}
That is due to the following theorem from
\cite{ShimonyEtAl2022SOCS}:

\begin{theorem}
\label{th:triconnected}
Let $G$ be a triconnected graph. Then, for every $s,w_1,w_2,t \in V$,
there exists a simple path in $G$ from $s$ to $t$ that includes $w_1$ and $w_2$. \end{theorem}

For the CEP, observe that by Dirac's theorem \cite{Dirac1960}, in every
$k$-connected graph, $k\geq2$, given any set of $k$ vertices, there exists
a simple cycle containing these vertices.
This theorem can be extended
to $m\in\{1,2\}$ edges and $k-m$ vertices. This extension is known as "folklore" of internet discussions \cite{k-connect,k-connected}, and
as we could not find it formally published, we state a proof
outline, as follows.
For $k=2$, given edges $(u,u')$ and $(v,v')$,
existence of a cycle follows immediately from the extension to Menger's theorem. 
Indeed, since there exist two vertex-disjoint paths from $u,u'$ to $v,v'$, 
they can be spliced with the given edges, resulting with the cycle as required.
The case where an edge and a vertex are given is similar.
For $k>2$, this is proved by induction:
there exists a cycle $L$ with the first $k-1$ elements,
including the one or two given edges. If the last vertex $w$ is not
already on $L$, there are $k$ vertex-disjoint (other than at $w$) simple paths $P_i$ from $w$ to vertices in $L$.
Since $L$ has only $k-1$ path segments between its given $k-1$ elements, at least two of the $P_i$, w.l.o.g. $P_1,P_2$, end at some vertices $w_1,w_2$ lying at one of these segments, $P$. 
Then, the required cycle is the same as $L$ where path $P$ is spliced to include the detour $P_1,w,P_2$. Therefore,
as mentioned in the introduction,
in a triconnected graph, the
answer to the CEP is positive
for any set of three elements
in which at most two are edges. Formally:
%\YD{I suggest making this statement a theorem, instead of lemma.}
%(a lemma is usually an auxiliary statement).}

\begin{theorem}\label{th:CEP}
In a triconnected graph $G$, given a set $S$ of three elements of $G$
where $S$ contains at most two edges,
there exists a simple cycle in $G$ that
traverses all the elements of $S$.
\end{theorem}

%This leaves the case where the three given 
%elements are all edges
%as the only non-trivial CEP problem
%of this type 
%in triconnected graphs.

%***********************************
\section{Cycle Existence Problem in Triconnected Graphs}\label{sec:triconnected}

We have seen that the only nontrivial case for the CEP in triconnected graphs is when the three given elements are edges.
%\YD {Is the style with so close by doubling of almost the same sentence OK? Maybe to  remove the last sentence of the previous section?}
%\ES{Nope, a short while ago these sentences were far away from each other. Former deleted as you asked.}
In this section, we prove the following main theorem, thus completing the solution to the CEP in triconnected graphs:
% \YDT{In this section, we show necessary and sufficient conditions for cycle existence in triconnected graphs. In this section, we prove the following main theorem:}

\begin{theorem}\label{th:cycle}
Let $G=(V,E)$ be a triconnected graph, and
$e_1, e_2,e_3\in E$.
Then, there exists a simple cycle $L$ such that
$e_1, e_2,e_3\in L$ if and only if
neither of the following conditions occur:
\begin{enumerate}
%    \item There is a 
%    vertex $v\in e_1,e_2,e_3$.\\
    \item $e_1,e_2,e_3$ all share a common end-vertex.
    \item $\{ e_1, e_2,e_3\}$ is an edge cut of $G$.
\end{enumerate}
\end{theorem}

%\YD {Maybe, we will number the statements (theorems, lemmas, etc.) by sections? For example, Theorem 3.1 seems sounding better than Theorem 9. And Theorem 2.x would clearly sound as a basic one.}\ES{OK, but later please.}

Note that the \emph{only if} part of the theorem statement is straightforward: exception 1
trivially implies that the edges cannot be on a simple cycle; for exception 2, this is due to every edge cut intersecting any cycle by an even number of edges (see Section~\ref{sec:basics}).
The rest of this section proves the \emph{if} part. Since
the latter is rather complex, we begin with proving some auxiliary statements.
First, we consider the \emph{main case} where all six end-vertices of the given edges are distinct, postponing the other cases to the very end.
%proving a simpler case and some lemmas required to show the main result.
Henceforth, we denote the \emph{end-vertices of $e_3$} by $s$ and $t$.

Let us introduce some notation.
For any path $P$, % from $x$ to $y$, 
we denote by $b(P)$ its starting vertex, by $f(P)$ its ending vertex, and by $\bar P$ the path reverse to $P$.
%When some set of paths is fixed, ... from that set 
We denote by $P_{xy}$ the path from a given collection with $b(P)=x$ and $f(P)=y$, if this does not lead to ambiguity,
and the reverse path $\bar P_{xy}$ by$P_{yx}$.
%the path reverse to $P_{xy}$.
For two vertices $x'$ and $y'$ lying on path $P$ in that order, we denote by $P[x'..y']$ the sub-path of $P$ from $x'$ to $y'$.
%\YD {Delete (since this was already said in Section 2.1):
%In what follows, when we assume
%or require the existence of a path, $P$, between two objects (vertices, paths, cycles), we assume, w.l.o.g., that only the end-vertices of $P$ are in those objects; indeed, otherwise, some sub-path of $P$ satisfies this property.}

%We begin with the following useful lemma:

\begin{lemma}\label{lem:path-splice}
Let $G=(V,E)$ be a 
%[[ YD: deleted: (vertex) 3-connected ]]
graph, $s,t\in V$
and $e_1=(u,u')$, $e_2=(v,v')$ edges in $E$, where all above vertices are distinct in $G$.
If there are a simple cycle 
$L=((u,u'),P_{u'v'},(v',v),P_{vu})$ in $G$, such that 
$s,t \not \in L$, and mutually vertex-disjoint paths
$P_1$ from $s$ to $L$ and $P_2$ from $t$ to $L$,
either both ending at $P_{vu}$ or both ending  at $P_{u'v'}$,
then there exists a simple path from $s$ to $t$ containing $e_1$ and $e_2$.
\end{lemma}

\begin{proof}
%Denote the final vertex on a path $P_i$ by $f(P_i)$, and the respective beginning vertex by $b(P_i)$.
Assume w.l.o.g.\ that $f(P_1),f(P_2)\in P_{vu}$
and that $f(P_1)$ is on the part of $P_{vu}$ between $v$ and $f(P_2)$. 
Then, the required $s,t$-path is:

\[
(P_1, {P_{uv}}[f(P_1)..v], (v,v'), P_{v'u'}, (u',u), {P_{uv}}[u..f(P_2)],\bar P_2).
\]
\end{proof}

We now prove the main case of Theorem~\ref{th:cycle} ($e_1, e_2, e_3$ do not have a common end-vertex) for a certain special case, thus providing a template for the general case.
We now assume that $G^- = G\setminus (s,t)$ is still triconnected.
%\YD {Why do we need the next two sentences, Eyal? We prove the if part, so they are not relevant.}
%\ES{Not sure any more. At some point I thought these were an explanation I added because you said something was not clear. But I am not even sure of that,}
%As a consequence, 
%exception 2 in Theorem \ref{th:cycle} implies
%a 2-edge cut in a triconnected graph $G^-$, so cannot
%occur due to Theorem \ref{l:edge-vertex cuts},
%, because $G^-$ is triconnected.
%In this special case we show that the theorem holds, as follows.
%By the assumption of the main case, exception 1 in Theorem \ref{th:cycle} also cannot occur. 
%Now, we show our special case result:

\begin{proposition}\label{prop:tri+cycle}
Let $G^-=(V,E^-)$ be a 
triconnected graph, $s,t$ vertices
and $e_1=(u,u')$, $e_2=(v,v')$ edges in $G^-$ (all above vertices distinct in $G^-$), such
that $(s, t) \not \in E^-$.
%Let $G' =G \cup (s,t)$.
Then, there exists a simple cycle
containing $e_1, e_2, (s,t)$ in $G = G^- \cup (s,t)$.
\end{proposition}

\begin{proof}
By the extension to Menger's theorem, there exist three vertex-disjoint paths
from $\{s,u,u'\}$ to
$\{t,v,v'\}$. Fix such a set of paths
${\cal P}$. Since both above vertex sets have cardinality 3, 
let $M$ be the symmetric mapping between these vertex sets, with $M(x)=y$ just when the path $(x...y) \in {\cal P}$.
There are two possible cases,
and in both cases below we construct a simple %vertex-disjoint 
path from $s$ to $t$ through $e_1, e_2$ in $G^-$, to which we then can add $(s,t)$ to complete the required cycle in $G$.

{\bf Case 1}:
If $M(t)=u$ or $M(t)=u'$, 
then  $M(s)=v$ or $M(s)=v'$.
Assume the former in each case, w.l.o.g.\
(due to symmetry).
%
%\YD {1) stop here, as it was previously (preferred, as I explained to you)}
%\YD {2) otherwise, we can flip the vertices' names}
%we can just rename $u$ to be $M(t)$,
%and $u'$ as the other end of $(u,u')$, and %likewise for $(v,v')$.)
Therefore, we have vertex-disjoint paths $P_{tu}, P_{sv},  P_{v'u'} \in {\cal P}$ with names denoting their respective endpoints.
Then, $P=(P_{sv}, (v,v'), P_{v'u'}, (u',u), P_{ut})$ is a simple path from $s$ to $t$ as required.
%; from here and on, let $P_{yx}$ denote the reverse path to $P_{xy}$.

{\bf Case 2}: If $M(t)=s$, then, w.l.o.g., $M(v)=u$ and $M(v')=u'$.
Thus, we have $P_{st}, P_{vu}, P_{v'u'} \in {\cal P}$, all vertex-disjoint. 
Therefore, $(P_{vu}, (u,u'), P_{u'v'}, (v',v))$ is a simple cycle,
which we denote by $L$.
Since $(s,t) \not \in E$,
path $P_{st}$ must have at least three vertices.
Since $L$ has at least four vertices, by the extension of Menger's theorem,
there exist three vertex-disjoint paths $P_1, P_2, P_3$ from some vertices in $P_{st}$ to $L$.

Then, at least two of $f(P_1),f(P_2),
f(P_3)$ must be on either $P_{vu}$ or $P_{v'u'}$.
Assume, w.l.o.g., (the other cases are symmetrical)
that $f(P_1),f(P_2)\in P_{vu}$
and that $b(P_1)$ is on the part of $P_{st}$
between $s$ and $b(P_2)$.
Then, $P'_1=(P_{st}[s..b(P_1)]\cdot P_1)$ and $P'_2=(P_{ts}[t..b(P_2)]\cdot P_2)$
are vertex-disjoint paths from $s$ and $t$ to
$P_{vu}$.
Note that $L$ is a simple cycle intersecting vertex-disjoint paths $P'_1$ and $P'_2$ only at $f(P'_1)$ and $f(P'_2)$.
Now, we have the required simple path
from $s$ to $t$ due to Lemma \ref{lem:path-splice}.
\end{proof}

We now move back to the general case, where $G$
is triconnected but $G\setminus (s,t)$ may be not
triconnected. We use the same proof
outline as in Proposition \ref{prop:tri+cycle}.
That is, there exist three vertex-disjoint paths from $\{s,u,u'\}$ to
$\{t,v,v'\}$.
We have the same two cases w.r.t.\ mapping $M$. 
The construction in case 1 does not require edge $(s,t)$,
so edge $(s,t)$ can be added to $P$
to complete the required cycle as before, and we are done.
In case 2, %however,
if $P_{st}$ consists  only of edge $(s,t)$
then the proof fails, because the cardinality of the vertex set of $P_{st}$ is only 2.
Note, however, that if there exists any path in $G$ 
from $s$ to $t$ other than  $(s,t)$
that is vertex-disjoint with cycle $L=((v,v'), P_{v'u'},(u',u),P_{uv})$, we can still
apply the method of case 2.
So, we %\YDT {only(?)} 
need to prove the remaining
case, stated as the following Lemma, in which the only $s,t$-path that is vertex-disjoint  from $L$ is $(s,t)$. 
%\ES{What did you mean above? [T]}
%\YD {Delete (since is assumed immediately after the statement of Th. 10): This is stated as the
%following Lemma, assuming for now that
%all the vertices in the in $e_1,e_2,e_3$
%are distinct.}
Note that in this case,  
$L$ contains neither $s$ nor $t$,
since by the construction,
the paths $P_{uv}$, $P_{u'v'}$,
$P_{st}=(s,t)$ are all mutually vertex-disjoint.

\begin{lemma}
\label{lem:distinct}
Let $G=(V,E)$ be a 
triconnected graph, $e_1=(u,u')$, $e_2=(v,v'),(s, t)$ edges in $G$ (all above vertices distinct in $G$).
Let $L$ be a cycle in $G$ such that
$e_1,e_2 \in L$, $s,t \not \in L$,
and there is no path from $s$ to $t$ in 
$G\setminus L\setminus (s,t)$.
Then, there exists a simple cycle in $G$ containing $e_1,e_2, (s,t)$ unless
$e_1,e_2, (s,t)$ form an edge cut of $G$.
\end{lemma}

\begin{proof}

W.l.o.g. let $L=(e_1,P_{uv},e_2,P_{v'u'})$.
We call a simple path with end-vertices in $L\cup \{ s,t\}$ but no other vertices in this set a {\em connecting path}, or a {\em connector}
for short. 
We call a connector from $s$ or $t$
to $L$ an {\em external} connector,
and specifically an $s$-connector
or $t$-connector, respectively.
A connector from $P_{uv}$ to $P_{u'v'}$ that 
is not
either $e_1$ or $e_2$ we call 
a {\em bridge}, 
and one either from
$P_{uv}$ to $P_{uv}$
or from $P_{u'v'}$ to $P_{u'v'}$ a {\em bypass}.

As $G$ is triconnected, there exist three vertex-disjoint (other than at $s$) paths from $s$ to distinct vertices on $L$. 
Since at most one of  such three paths %the set of such external connectors
can be through $t$, then
at least two of them are $s$-connectors, and if $t$ is on the third path, then its suffix from $t$ is a $t$-connector.
We call such three connectors an \emph{$s$-triple}.
%, and thus at least two external connectors for $s$.
Symmetrically, there exist three vertex-disjoint (other than at $t$) paths from $t$ to distinct vertices on $L$, so that
%Since at most one of the set of such connectors can be though $s$,
at least two of them are $t$-connectors, and if $s$ is on the third path, then its suffix from $s$ is an $s$-connector;
we call them  a \emph{$t$-triple}.
%\YD {Delete:
%Additionally, since there are three vertex-disjoint (other than at $s$) paths from $s$ to $L$, there are at least three distinct end-vertices ending these
%external connectors on $L$.}
Note that $s$-connectors are always
internally vertex-disjoint from
$t$-connectors, as otherwise we would 
have a path from $s$ to $t$
in $G\setminus L \setminus (s,t)$. 

We call $P_{uv}$ and $P_{u'v'}$ the \emph{left and right sides} of $L$, respectively.
%, and assume, w.l.o.g., that there exists a $t$-connector ending on $P_{u'v'}$.
%\YD {Delete?: which we will call "the $t$-side of $L$";
%we also call $P_{uv}$ "the $s$-side of $L$"}.
There are several cases w.r.t.\ the sides of $L$
where external connectors end, which we examine below.
%\YD {I suggest canceling this list here, since it doubles the cases' titles.}
%(1) There exist $s$-connectors to both sides of $L$, and likewise w.r.t.\ $t$-connectors.
%(2) All $t$-connectors end on the $t$-side of $L$, and there is at least one $s$-connector ending there.
% on the same side (or symmetrically exchanging the role of $s$ and $t$ above). 
%(3) All $t$-connectors end on the $t$-side side of $L$ and all $s$-connectors end on its $s$-side.
In all these cases, we show a simple path from $s$ to $t$ 
containing $e_1$ and $e_2$ (except for
when $\{e_1,e_2, (s,t)\}$ is a 3-edge cut), from which we can construct the required cycle by adding $(s,t)$.

%\YDT {At most one of them contains $t$, and thus can be abridged to a connector from $t$.}
%Assume, w.l.o.g.,
%that there are two distinct such end-vertices $v_1,v_2$
%that end such disjoint paths 
%on $P_{uv}$.
%Now, we have the following cases,
%in which 

{\bf Case 1}: 
\emph{There exist both $s$- and $t$-connectors to
both sides of $L$.} 
%As there are at least three vertices where connectors end on $L$, 
Assume, w.l.o.g., that two of the (vertex-disjoint other than at $s$) connectors in the $s$-triple, $P_1$ and $P_2$, end on $P_{uv}$.
%; denote themby $v_1,v_2$.} 
Then, the requisite $s$ to $t$ path exists due to Lemma \ref{lem:path-splice}, as follows.
If $P_1$ and $P_2$ are $s$- and $t$-connectors, then $P_1,P_2$
constitute the Lemma conditions.
Otherwise, $P_1,P_2$ are both $s$-connectors. 
Since some $t$-connector, $P$, ending at $P_{uv}$ exists, and it is internally vertex-disjoint from $P_1$ and $P_2$, 
either $P$ and $P_1$ or $P$ and $P_2$ are vertex-disjoint,
which also suffices.
%we have an external connector from $s$ to $P_{uv}$ that is vertex-disjoint from an external connector from $t$ to $P_{uv}$, so the Lemma conditions hold.
%But suppose that the same vertex (say $v_1$) ends both the connector for $s$ and the connector for $t$. Then since $v_2$ also ends a connector for either $s$ or $t$, then the Lemma conditions hold here as well.

{\bf Case 2}: 
\emph{ All the
$t$-connectors end at $P_{u'v'}$,
and there is at least one $s$-connector, $P_s$,
ending there.} % $f(P_s)\in P_{u'v'}$. 
\emph{(The other variants are either $s$/$t$ or left/right symmetric or both, and thus can be treated similarly.)}
Then, since there are at least two vertex-disjoint (other than at $t$)  $t$-connectors, $P_1, P_2$, we have:
$f(P_1),f(P_2)\in P_{u'v'}$
and $f(P_1)\not = f(P_2)$.
The ending vertex of $P_s$ must thus be distinct from either $f(P_1)$, or $f(P_2)$, or both. Therefore, $P_s$ is vertex-disjoint from either $P_1$ or $P_2$ (or both).
Thus, the required simple path from $s$ to $t$ exists due to Lemma \ref{lem:path-splice}.

{\bf Case 3}: \emph{$($Neither case 1 nor case 2 occur, that is$)$ w.l.o.g., every $s$-connector ends at $P_{uv}$ and every $t$-connector ends at $P_{u'v'}$.}
Thus, $P_{uv}$ and $(s,t)$ separate $s$ from $P_{u'v'}$ and $t$, and likewise
$P_{u'v'}$ and $(s,t)$ separate $t$ from $P_{uv}$ and $s$. Therefore,
every $s$-connector is disjoint from every $t$-connector,
and both are internally vertex-disjoint from every bridge. Here too there are 2 sub-cases:

{\bf Case 3a}: \emph{No bridge exists.}
Then, $e_1, e_2, (s,t)$ form an edge cut of $G$, separating $G$ into its $s$-side and $t$-side, which is exclusion 2 of the theorem.
In this case,  
starting at $s$, we need to cross each of $e_1, e_2, (s,t)$ exactly once, so after any three crossings we end up on the  $t$-side of $G$ with no way back.

{\bf Case 3b}: \emph{At least one bridge exists.}
Let us show that here, the desired simple path from $s$ to $t$ can be constructed. 
Since $G$ is triconnected, every vertex $w\in L$ must have an incident edge $e$ not on $L$, 
as otherwise its two neighbors on $L$ form a 2-vertex cut.
We categorize such edges $e$ into:
\begin{enumerate}
    \item $e$ is on some external connector, thus called an \emph{external edge}.
    \item $e$ is on some bridge, thus called a \emph{bridge edge}.
    \item None of the above, in which case
    $e$ is called a  \emph{bypass edge}.
\end{enumerate}
We likewise call a vertex $w$ \emph{external, bridge}, or
\emph{bypass}, respectively, when $w$ has an incident edge of the respective type. 
Note that all vertices $w$ on  $L$ must be of at least one such type, but the types are non-exclusive: it is possible for $w$ to be bypass, and external, and bridge. 
Recall that at least two vertices on each path $P_{uv}$ and $P_{u'v'}$ must be external. 
Every external vertex $w$ on $P_{uv}$, resp., $P_{u'v'}$, is an end-vertex of some (not necessarily unique) external $s$-, resp., $t$-connector; % for either $s$ or (exclusive) $t$,
%which we denote by $P(w)$.
for each external vertex, we fix one such connector, denoting it by $P_{sw}$, resp., $P_{wt}$.
%\YDT {Delete?: For convenience, we orient each bridge, $B$, from $P_{uv}$ to $P_{u'v'}$, so that $b(B) \in P_{uv}$ and $f(B) \in P_{u'v'}$.}
For a bridge $B$, we denote by $v_B$ its end-vertex on $P_{uv}$ and by $v'_B$ its end-vertex on $P_{u'v'}$.

Now examine bridge end locations w.r.t. external vertex locations. For clarity,
we call the direction that is towards the $(v,v')$
edge on paths $P_{uv}, P_{u'v'}$ the {\em north}
direction, and that towards $(u,u')$ {\em south}.
Let $v_{N},v_{S}$ be the most northern (resp, most southern) external vertex on
$P_{uv}$,  and likewise $v'_{N},v'_{S}$ for vertices on $P_{u'v'}$. Since there are at least two external
vertices in each part of $L$, these vertices are all distinct.
For convenience, we will use $>$ to denote
"more northern than", which is defined only between
vertices on the same path $P_{uv}$ or $P_{u'v'}$.
There are now several cases, illustrated in Figure~\ref{fig:case3b}.

\begin{figure*}
    \centering
\resizebox{1.0\columnwidth}{!}{
\begin{tabular}{c||c||c}

\begin{tikzpicture}[minimum size=17mm,
  node distance=0.7cm and 2.5cm,
  >=stealth,
  bend angle=45,
  auto]
\tikzstyle{every node}=[font=\Huge]
\node (0) [circle, ultra thick, text centered, text width=0.5cm, draw=black] {$s$};
\node (1) [circle, ultra thick, text centered, text width=0.5cm, draw=black, right=of 0] {~};
\node (2) [circle, ultra thick, text centered, text width=0.5cm, draw=black, above=of 1] {~} edge [ultra thick, color=red] (1);
\node (3) [circle, ultra thick, text centered, text width=0.5cm, draw=black, below=of 1] {$\! v_S$};
\node (4) [circle, ultra thick, text centered, text width=0.5cm, draw=black, below=of 3] {~};
\node (5) [circle, ultra thick, text centered, text width=0.5cm, draw=black, below=of 4] {$u$} edge [color=red] (4);
\node (6) [circle, ultra thick, text centered, text width=0.5cm, draw=black, above=of 2]
{$\!\! v_B$} edge [color=red] (2);
\node (7) [circle, ultra thick, text centered, text width=0.5cm, draw=black, above=of 6] {$\!\! v_N$} edge [ultra thick, color=red] (0);
\node (8) [circle, ultra thick, text centered, text width=0.5cm, draw=black, above=of 7] {$v$} edge [color=red] (7);
\node (10) [circle, ultra thick, text centered, text width=0.5cm, draw=black, right=of 1] {$v'_{N}$};
\node (11) [circle, ultra thick, text centered, text width=0.5cm, draw=black, above=of 10] {~};
\node (12) [circle, ultra thick, text centered, text width=0.5cm, draw=black, above=of 11] {~};
\node (13) [circle, ultra thick, text centered, text width=0.5cm, draw=black, above=of 12] { $\!\! v'_B$} edge [color=red] node [anchor=east] {B} (6);
\node (14) [circle, ultra thick, text centered, text width=0.5cm, draw=black, below=of 10] {~};
\node (15) [circle, ultra thick, text centered, text width=0.5cm, draw=black, below=of 14] {$\!\! v'_S$};
\node (16) [circle, ultra thick, text centered, text width=0.5cm, draw=black, below=of 15] {$u'$} edge [color=red] (5) edge [color=red] (15) edge [color=red] node {$e_1$} (5);
\node (17) [circle, ultra thick, text centered, text width=0.5cm, draw=black, above=of 13] {$v'$} edge [color=red] (13) edge [color=red] node [anchor=south] {$e_2$} (8);
\node (20) [circle, ultra thick, text centered, text width=0.5cm, draw=black, right=of 10] {$t$}
edge [color=red] (15);

\draw [thick,-,>=stealth] (0) to (3);
\draw [thick,-,>=stealth,color=red] (1) to (3);
\draw [thick,-,>=stealth,color=red] (3) to (4);
\draw [thick,-,>=stealth] (6) to (7);
\draw [thick,-,>=stealth] (10) to (20);
\draw [thick,-,>=stealth] (10) to (11);
\draw [thick,-,>=stealth] (11) to (12);
\draw [thick,-,>=stealth] (12) to (13);
\draw [thick,-,>=stealth] (10) to (14);
\draw [thick,-,>=stealth] (14) to (15);
\draw [thick,-,>=stealth] (10) to (20);

\end{tikzpicture}

&

\begin{tikzpicture}[minimum size=17mm,
  node distance=0.7cm and 2.5cm,
  >=stealth,
  bend angle=45,
  auto]
\tikzstyle{every node}=[font=\Huge]
\node (0) [circle, ultra thick, text centered, text width=0.5cm, draw=black] {$s$};
\node (1) [circle, ultra thick, text centered, text width=0.5cm, draw=black, right=of 0] {~};
\node (2) [circle, ultra thick, text centered, text width=0.5cm, draw=black, above=of 1] {$\!\! v_N$} edge  (1);
\node (3) [circle, ultra thick, text centered, text width=0.5cm, draw=black, below=of 1] {$\!\! v_S$} edge [ultra thick,color=red]  (0);
\node (4) [circle, ultra thick, text centered, text width=0.5cm, draw=black, below=of 3] {$~$};
\node (5) [circle, ultra thick, text centered, text width=0.5cm, draw=black, below=of 4] {$u$} edge [color=red] (4);
\node (6) [circle, ultra thick, text centered, text width=0.5cm, draw=black, above=of 2] { $\!\! v_B$} edge  (2);
\node (7) [circle, ultra thick, text centered, text width=0.5cm, draw=black, above=of 6] {~} ;
\node (8) [circle, ultra thick, text centered, text width=0.5cm, draw=black, above=of 7] {$v$} edge [color=red] (7);
\node (10) [circle, ultra thick, text centered, text width=0.5cm, draw=black, right=of 1] {$\!\! v'_N$};
\node (11) [circle, ultra thick, text centered, text width=0.5cm, draw=black, above=of 10] {~};
\node (12) [circle, ultra thick, text centered, text width=0.5cm, draw=black, above=of 11] {~};
\node (13) [circle, ultra thick, text centered, text width=0.5cm, draw=black, above=of 12] {~} ;
\node (14) [circle, ultra thick, text centered, text width=0.5cm, draw=black, below=of 10] { $\!\! v'_B$}
edge [color=red] node [anchor=east] {B} (6);
\node (15) [circle, ultra thick, text centered, text width=0.5cm, draw=black, below=of 14] {$\!\! v'_S$};
\node (16) [circle, ultra thick, text centered, text width=0.5cm, draw=black, below=of 15] {$u'$} edge [color=red] (5) edge [color=red] (15) edge [color=red] node {$e_1$} (5);
\node (17) [circle, ultra thick, text centered, text width=0.5cm, draw=black, above=of 13] {$v'$} edge [color=red] (13) edge [color=red] node [anchor=south] {$e_2$} (8);
\node (20) [circle, ultra thick, text centered, text width=0.5cm, draw=black, right=of 10] {$t$}
edge (15);

\draw [ultra thick,-,>=stealth] (0) to (2);
\draw [thick,-,>=stealth] (1) to (3);
\draw [thick,-,>=stealth,color=red] (3) to (4);
\draw [thick,-,>=stealth,color=red] (6) to (7);
\draw [thick,-,>=stealth,color=red] (10) to (11);
\draw [thick,-,>=stealth,color=red] (11) to (12);
\draw [thick,-,>=stealth,color=red] (12) to (13);
\draw [thick,-,>=stealth] (10) to (14);
\draw [thick,-,>=stealth,color=red] (14) to (15);
\draw [thick,-,>=stealth,color=red] (10) to (20);

\end{tikzpicture}

&

\begin{tikzpicture}[minimum size=17mm,
  node distance=0.7cm and 2.5cm,
  >=stealth,
  bend angle=45,
  auto]
\tikzstyle{every node}=[font=\Huge]
\node (0) [circle, ultra thick, text centered, text width=0.5cm, draw=black] {$s$};
\node (1) [circle, ultra thick, text centered, text width=0.5cm, draw=black, right=of 0] {~};
\node (2) [circle, ultra thick, text centered, text width=0.5cm, draw=red, above=of 1] {$\!\! v_N$} edge [ultra thick] (1)
edge (0);
\node (3) [circle, ultra thick, text centered, text width=0.5cm, draw=black, below=of 1] {$\!\! v_S$};
\node (4) [circle, ultra thick, text centered, text width=0.5cm, draw=black, below=of 3] {~};
\node (5) [circle, ultra thick, text centered, text width=0.5cm, draw=black, below=of 4] {$u$} edge (4);
\node (6) [circle, ultra thick, text centered, text width=0.5cm, draw=black, above=of 2] { $\!\! v_B$} edge  (2);
\node (7) [circle, ultra thick, text centered, text width=0.5cm, draw=black, above=of 6] { $\!\! v_{B'}$}
edge node [anchor=west] {B'} (11);
\node (8) [circle, ultra thick, text centered, text width=0.5cm, draw=black, above=of 7] {$v$} edge (7);
\node (10) [circle, ultra thick, text centered, text width=0.5cm, draw=red, right=of 1] {$\!\! v'_N$};
\node (11) [circle, ultra thick, text centered, text width=0.5cm, draw=black, above=of 10] { $\!\! v'_{B'}$};
\node (12) [circle, ultra thick, text centered, text width=0.5cm, draw=black, above=of 11] {~};
\node (13) [circle, ultra thick, text centered, text width=0.5cm, draw=black, above=of 12] {$\!\! v'_B$} edge node [anchor=south] {B} (6);
\node (14) [circle, ultra thick, text centered, text width=0.5cm, draw=black, below=of 10] {~};
\node (15) [circle, ultra thick, text centered, text width=0.5cm, draw=black, below=of 14] {$\!\! v'_S$};
\node (16) [circle, ultra thick, text centered, text width=0.5cm, draw=black, below=of 15] {$u'$} edge  (5) edge  (15) edge  node {$e_1$} (5)
edge  (4);
\node (17) [circle, ultra thick, text centered, text width=0.5cm, draw=black, above=of 13] {$v'$} edge  (13) edge  node [anchor=south] {$e_2$} (8);
\node (20) [circle, ultra thick, text centered, text width=0.5cm, draw=black, right=of 10] {$t$}
edge  (15);

\draw [thick,-,>=stealth] (0) to (3);
\draw [thick,-,>=stealth] (1) to (3);
\draw [thick,-,>=stealth] (3) to (4);
\draw [thick,-,>=stealth] (6) to (7);
\draw [thick,-,>=stealth] (10) to (20);
\draw [thick,-,>=stealth] (10) to (11);
\draw [thick,-,>=stealth] (11) to (12);
\draw [thick,-,>=stealth] (12) to (13);
\draw [thick,-,>=stealth] (10) to (14);
\draw [thick,-,>=stealth] (14) to (15);
\draw [thick,-,>=stealth] (10) to (20);

\end{tikzpicture}

\end{tabular}
}
    \caption { Case 3b1 (left), Case 3b2  (center), impossibility in Case3b3 (right). North is up. The desired paths and  the 2-separator are in red.
 %   \YD {Let us remove the name $B''$ of the bridge in the bottom right.}
    }
\label{fig:case3b}
\end{figure*}
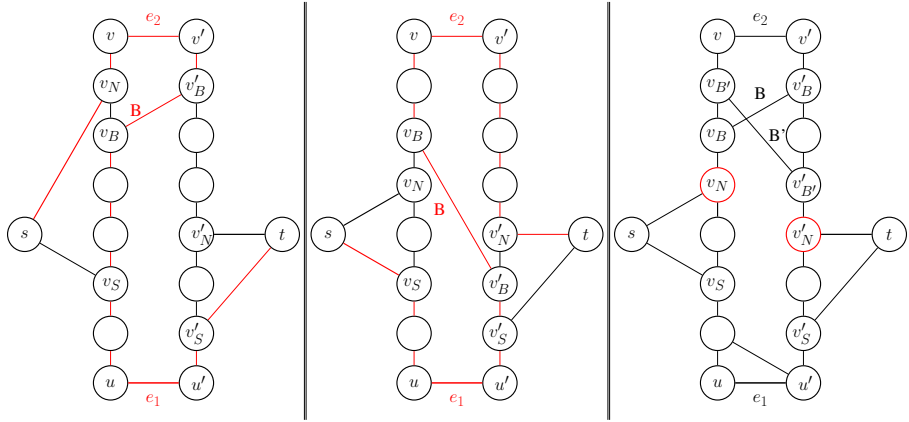

{\bf Case 3b1}: There is a bridge, $B$, such that $v_{N}>v_B$
and $v'_B>v'_S$ (Figure \ref{fig:case3b} (left)).
Then, the desired $s,t$-path is:
\[
(P_{sv_N},P_{uv}[v_{N}..v], (v,v'), P_{v'u'}[v'..v'_B], B, P_{vu}[v_B..u],
(u,u'), P_{u'v'}[u'..v'_S],P_{v'_St}).
\]

{\bf Case 3b2}: Symmetrically, there is a bridge, $B$, such that $v'_{N}>v'_B$
and $v_B>v_{S}$ (Figure \ref{fig:case3b} (center)). 
Then, the desired $s,t$-path is:
\[
(P_{sv_S},P_{uv}[v_S..u],(u,u'),P_{u'v'}[u'..v'_B], B, P_{vu}[v_B..v],(v,v'),P_{v'u'}[v'...v'_N],P_{v'_Nt}).
\]
Note that cases 3b1 and 3b2 are not necessarily mutually exclusive. 

{\bf Case 3b3}: Neither of the above two cases hold (Figure \ref{fig:case3b} (right)).
This implies that for every bridge $B$, either
$(v_B \geq v_{N})\wedge (v'_B \geq v'_N)$ (``northern bridge''),
or $(v_{S}\geq v_B) \wedge (v'_S \geq v'_B)$ (``southern bridge'').
Since at least one bridge must exist,
w.l.o.g., assume that a northern bridge exists
%the first disjunct w.r.t.\ some bridge 
(the other case is symmetric).
Let $B$ be a northern bridge with southernmost $v_B\geq v_N$, and $B'$ a northern bridge (possibly the same as $B$) with southernmost $v'_{B'} \geq v'_N$.

\paragraph{Necessity of bypasses} 
We show below that without appropriate bypasses, G would have a 2-vertex separator, such as $\{v_N, v'_N\}$ shown in Figure \ref{fig:case3b} (right).
Let us define notation for bypasses.
For a bypass  $I$, we denote by $N(I)$ and $S(I)$
its  north and south end-vertices, respectively.
%as the south end vertex of $I$.
We say that $I$ is a {\em bypass of vertex $x$} if $N(I)>x>S(I)$.
A bypass is called \emph{pure} if it has no internal elements in common with any bridge or external connector.

Let us show that either there exists a bypass of every $z$ with $v_B \geq z \geq v_N$, or 
there exists a bypass of every $z'$ with $v'_{B'} \geq z \geq v'_N$, or both. 
%\ES{Or possibly both. Do we need to say that?}
Assume the contrary.
Then, there exists a pair of vertices $z$, $v_B \geq z \geq v_N$, and $z'$, $v'_{B'} \geq z' \geq v'_N$,
with no bypass of either $z$ or $z'$.
Let us choose a vertex $w$ as follows.
If $z=v$ and $z'=v'$, then the end-vertices of $B$ are $v,v'$; since $(v,v')$ by definition is not a bridge, $B$ contains an internal vertex, which we denote by $w$.
Otherwise, let $w$ be $v$ if $z \neq v$, else $w=v' \neq z^\prime$.
We show now that $\{z,z^\prime\} $ separates $w$ from $\{s,t\}$, in contradiction to the triconnectivity of $G$, thus invalidating our assumption.

Assume, to the contrary, that $P$ is a simple path from $s$ to $w$ that includes neither $t$, nor $z$, nor $z'$.
Note that if $w \not \in L$, then $P$ could be extended along the right part of $B$ to $P_{u'v'}$.
Since $L$ can only be reached from $s$ through a vertex on $[v_N..v_S]$ and $z\geq v_N$, $P$ must have at least one vertex $x$ on $P_{uv}$ such that $z>x$ (we ruled out $z=x$ as $z \not \in P$).
Let $x^*$ be the last vertex on $P$ such that either $z> x^*$ or $z' > x^*$ and $y$ be the first vertex after $x^*$ on $P$ that is either on $L$ or on a northern bridge; let $P^*=P[x^*..y]$.
Assume that $x^* \in P_{uv}$. If $y \in P_{uv}$, then $P^*$ is a bypass of $z$, a contradiction.
If $y$ is on a northern bridge, then the concatenation of $P^*$ with the part of that bridge from $y$ to $P_{u'v'}$ is a northern bridge from $x^*$ with 
$v_B > z > x^*$. a contradiction to the definition of $B$.
The case $x^* \in P_{u'v'}$ is symmetric.

For similar reasons, a simple path from $t$ to $w$ not including $s,z,z'$ is also impossible, thereby
completing the proof of our statement. Thus, in the rest of our analysis of Case 3b3, we assume, w.l.o.g., that there exists
a bypass for every $z$ with $v_B \geq z \geq v_N$.

\paragraph{Purity of bypasses} Let ${\cal I}$ be the set of bypasses 
for all $z$ such that $v_B\geq z\geq v_N$. 
Let us show that all bypasses in ${\cal I}$ are pure. Otherwise, let $I\in \cal I$ be impure. 
Suppose that $I$ has an internal element in common with bridge $\tilde B$. 
If $\tilde B$ is a northern bridge,
%By the definition of case 3b3, either $(v_{\tilde B} \geq v_{N}) \wedge (v'_{\tilde B} \geq v'_N)$ or $(v_{S}\geq v_{\tilde B}) \wedge (v'_S \geq v'_{\tilde B})$.
%In the former case, 
$I$ can be spliced with $\tilde B$ creating a bridge $B''$ with $v_{B''}=S(I) < v_B$.
Note that since $v'_{B''} = v'_{\tilde B} \geq v'_N$, 
%we also have $v_{B''} \geq v_{N}$, by the definition of Case 3b3.
$B''$ is also a northern bridge.
This contradicts $v_B$ being the southernmost end-vertex at $P_{uv}$ of a northern bridge. %with $v_B \geq v_N$.
Otherwise, $\tilde B$ is a southern bridge. Then,
%In the latter case, 
$I$ can be spliced with $\tilde B$ creating a bridge $B''$ with $(v_{B''}=N(I) > v_N  >v_S) \wedge (v'_N > v'_S \geq v'_{\tilde B} =v'_{B''})$.
This accords with the definition of Case 3b2, not Case 3b3, a contradiction.
% with Case 3b3 under consideration.

Suppose now that $I$ has an internal element in common with either $s$- or $t$-connector $P$.
In the former case, $I$ can be spliced with $P$ creating
an $s$-connector $P_{sN(I)}$.
But $N(I) > v_N$, so $v_N$ is not
the northernmost external vertex on $P_{uv}$, a contradiction.
In the latter case, $I$ can be spliced with $P$ creating
a $t$-connector with an end-vertex on $P_{uv}$, contradicting the definition of Case 3.

%Then no bypass in ${\cal I}$
%has internal elements (vertices or edges) %in common with $B$, and likewise no internal
%elements in common with $P_{sv_N}$, which
%is shown as follows.
%Assume in contradiction that 
%some $I\in {\cal I}$ has internal elements
%in common with $B$.
%Then  $I$ can be spliced with the bridge $B$, creating a
%bridge $B''$ with $v_{B''}=S(I) < v_B$,
%so $v_B$ is not the southernmost bridge %vertex with $v_B \geq v_N$, a contradiction,
%or we have $v_N > v_{B''}$, contradicting the assumption underlying case 3b3, since 
%$v'_{B''}>v'_S$.
%Likewise, if $I\in {\cal I}$ has internal elements
%in common with $P_{sv_N}$, then $I$ can be spliced with $P_{sv_N}$, creating
%an external connector $P_{sN(I)}$.
%But $N(I) > v_N$, so $v_N$ is not
%the northernmost external vertex, a contradiction.

\paragraph{Bypass graphs} Let $v_{No} > v_B$ be the northernmost vertex on $P_{uv}$ in any bypass in ${\cal I}$,
and likewise $v_{So} < v_N$ be the southernmost such vertex. 
Now, define graph $G_{BP}$ as follows:
$G_{BP}$ consists of all vertices and edges in
all bypasses in ${\cal I}$, and all vertices
and edges in $P_{uv}$ between $v_{No}$ and $v_{So}$. 
By construction, $G_{BP}$ is (vertex) 2-connected. 
Thus, there exist two vertex-disjoint paths between any two sets of vertices in $G_{BP}$ of cardinality at least 2.
%including specifically 
Let us fix a pair of such paths for the sets $\{v_{No}, v_B\}$ and $\{ v_N, v_{So}\}$. 
We denote them by $P_{xy}$ with $xy$ indicating
the respective subscripts of $v$. 
%Note that because all bypasses in $\cal I$ have
%no common elements with $B$ or
%with $P_{sv_N}$ (except possibly at
%$v_B$ and $v_N$) then concatenating any
%simple path from $G_{BP}$ at these vertices
%with $B$ or $P_{sv_N}$ results in a simple path.

%are pure,
%these paths do not use any internal vertex of any bridge or external connector.
%\YD { Delete?: , allowing construction of the simple paths as defined below.}
We consider the two cases for the pair of paths as above.
In both, the desired simple path begins with $P_{sv_N}$ and ends with $(P_{vu}[v_{So}...u], (u,u'), P_{u'v'}[u'..v'_S],P_{v'_St})$.
Let us describe the middle part of the desired path.

\begin{figure*}
    \centering
\resizebox{1.0\columnwidth}{!}{
\begin{tabular}{c||c}
\begin{tikzpicture}[minimum size=17mm,
  node distance=0.7cm and 2.5cm,
  >=stealth,
  bend angle=45,
  auto]
  \tikzstyle{every node}=[font=\Huge]
\node (0) [circle, ultra thick, text centered, text width=0.5cm, draw=black] {$s$};
\node (1) [circle, ultra thick, text centered, text width=0.5cm, draw=black, right=of 0] {$v_N$}
edge [line width=4.0pt,color=red] 
node [anchor=south] {$P_{sv_N}$} (0);
\node (2) [circle, ultra thick, text centered, text width=0.5cm, draw=blue, above=of 1] {$~$} edge [line width=4.0pt, color=blue] (1);
\node (3) [circle, ultra thick, text centered, text width=0.5cm, draw=black, below=of 1] {$v_S$};
\node (4) [circle, ultra thick, text centered, text width=0.5cm, draw=black, below=of 3] { $\!\!v_{So}$}  ;
\node (5) [circle, ultra thick, text centered, text width=0.5cm, draw=black, below=of 4] {$u$} edge [line width=4.0pt,color=red] (4);
\node (6) [circle, ultra thick, text centered, text width=0.5cm, draw=black, above=of 2] { $\!\! v_B$} edge (2) edge [line width=4.0pt,color=blue,bend left] 
 node [anchor=west,yshift=-3.5cm,xshift=-1.1cm] { $P_{BSo}$} (4);
\node (7) [circle, ultra thick, text centered, text width=0.8cm, draw=black, above=of 6] {$\!\! v_{\! N\! o}$} edge [line width=4.0pt,color=blue, bend right] 
node [anchor=east, color=blue] {$P_{NNo}$} (2);
\node (8) [circle, ultra thick, text centered, text width=0.5cm, draw=black, above=of 7] {$v$} edge [line width=4.0pt,color=red] (7);
\node (10) [circle, ultra thick, text centered, text width=0.5cm, draw=black, right=of 1] {$v'_N$};
\node (11) [circle, ultra thick, text centered, text width=0.5cm, draw=black, above=of 10] {$~$};
\node (12) [circle, ultra thick, text centered, text width=0.5cm, draw=black, above=of 11] {$~$};
\node (13) [circle, ultra thick, text centered, text width=0.5cm, draw=black, above=of 12] {$\!\! v'_B$} edge [line width=4.0pt,color=red] node [anchor=east] {B} (6);
\node (14) [circle, ultra thick, text centered, text width=0.5cm, draw=black, below=of 10] {$~$};
\node (15) [circle, ultra thick, text centered, text width=0.5cm, draw=black, below=of 14] {$v'_S$};
\node (16) [circle, ultra thick, text centered, text width=0.5cm, draw=black, below=of 15] {$u'$} edge [line width=4.0pt,color=red] (5) edge [line width=4.0pt,color=red] (15) edge [line width=4.0pt,color=red] node {$e_1$} (5);
\node (17) [circle, ultra thick, text centered, text width=0.5cm, draw=black, above=of 13] {$v'$} edge [line width=4.0pt,color=red] (13) edge [line width=4.0pt,color=red] node [anchor=south] {$e_2$} (8);
\node (20) [circle, ultra thick, text centered, text width=0.5cm, draw=black, right=of 10] {$t$}
edge [line width=4.0pt,color=red] 
node [anchor=west,yshift=-0.5cm] {$P_{v'_St}$} (15);

\draw [thick,-,>=stealth] (0) to (3);
\draw [thick,-,>=stealth] (1) to (3);
\draw [thick,-,>=stealth] (3) to (4);
\draw [thick,-,>=stealth] (6) to (7);
\draw [thick,-,>=stealth] (10) to (20);
\draw [thick,-,>=stealth] (10) to (11);
\draw [thick,-,>=stealth] (11) to (12);
\draw [thick,-,>=stealth] (12) to (13);
%\draw [thick,-,>=stealth] (6) -- node [anchor=north] {B'}  (12);
\draw [thick,-,>=stealth] (10) to (14);
\draw [thick,-,>=stealth] (14) to (15);
\draw [thick,-,>=stealth] (10) to (20);

\end{tikzpicture}

&

\begin{tikzpicture}[minimum size=17mm,
  node distance=0.7cm and 2.5cm,
  >=stealth,
  bend angle=45,
  auto]
  \tikzstyle{every node}=[font=\Huge]
\node (0) [circle, ultra thick, text centered, text width=0.5cm, draw=black] {$s$};
\node (1) [circle, ultra thick, text centered, text width=0.5cm, draw=black, right=of 0] {$v_N$}
edge [line width=4.0pt,color=red] 
node [anchor=north, xshift=-0.2cm] {$P_{sv_N}$}
(0);
\node (2) [circle, ultra thick, text centered, text width=0.5cm, draw=blue, above=of 1] {$~$} edge [line width=4.0pt, color=blue]  (1);
\node (3) [circle, ultra thick, text centered, text width=0.5cm, draw=black, below=of 1] { $\!\!v_{\small So}$}
;
\node (4) [circle, ultra thick, text centered, text width=0.5cm, draw=black, below=of 3] {$v_S$};
\node (5) [circle, ultra thick, text centered, text width=0.5cm, draw=black, below=of 4] {$u$} edge [line width=4.0pt,color=red] (4);
\node (6) [circle, ultra thick, text centered, text width=0.5cm, draw=black, above=of 2] {$\!\! v_B$} edge [line width=4.0pt,line width=4.0pt, color=blue] node [color=blue] {$P_{NB}$} (2);
\node (7) [circle, ultra thick, text centered, text width=0.5cm, draw=black, above=of 6] { $\!\!\! v_{\! No}$} 
edge [line width=4.0pt,color=blue,bend right]  node [anchor=east,yshift=1cm] {$P_{NoSo}$} (3);
\node (8) [circle, ultra thick, text centered, text width=0.5cm, draw=black, above=of 7] {$v$} edge [line width=4.0pt,color=red] (7);
\node (10) [circle, ultra thick, text centered, text width=0.5cm, draw=black, right=of 1] {$v'_N$};
\node (11) [circle, ultra thick, text centered, text width=0.5cm, draw=black, above=of 10] {$~$};
\node (12) [circle, ultra thick, text centered, text width=0.5cm, draw=black, above=of 11] { $~$};
\node (13) [circle, ultra thick, text centered, text width=0.5cm, draw=black, above=of 12] {$\!\! v'_B$} edge [line width=4.0pt,color=red] node [anchor=east] {B} (6);
\node (14) [circle, ultra thick, text centered, text width=0.5cm, draw=black, below=of 10] {$~$};
\node (15) [circle, ultra thick, text centered, text width=0.5cm, draw=black, below=of 14] {$v'_S$};
\node (16) [circle, ultra thick, text centered, text width=0.5cm, draw=black, below=of 15] {$u'$} edge [color=red] (5) edge [line width=4.0pt,color=red] (15) edge [line width=4.0pt,color=red] node {$e_1$} (5);
\node (17) [circle, ultra thick, text centered, text width=0.5cm, draw=black, above=of 13] {$v'$} edge [line width=4.0pt,color=red] (13) edge [line width=4.0pt,color=red] node [anchor=south] {$e_2$} (8);
\node (20) [circle, ultra thick, text centered, text width=0.5cm, draw=black, right=of 10] {$t$}
edge [line width=4.0pt,color=red] 
node [anchor=west,yshift=-0.5cm] {$P_{v'_St}$} (15);

%\draw [thick,-,>=stealth,snake=snake] (0) to (3);
\draw [thick,-,>=stealth] (0) to (4);
\draw [thick,-,>=stealth] (1) to (3);
\draw [line width=4.0pt,color=red,-,>=stealth] (3) to (4);
\draw [thick,-,>=stealth] (6) to (7);
\draw [thick,-,>=stealth] (10) to (20);
\draw [thick,-,>=stealth] (10) to (11);
\draw [thick,-,>=stealth] (11) to (12);
\draw [thick,-,>=stealth] (12) to (13);
%\draw [thick,-,>=stealth] (6) -- node [anchor=north] {~~~B'}  (12);
\draw [thick,-,>=stealth] (10) to (14);
\draw [thick,-,>=stealth] (14) to (15);
\draw [thick,-,>=stealth] (10) to (20);

\end{tikzpicture}

\end{tabular}
}
    \caption{Case 3b3 I (left), Case 3b3 II (right). North is up. The desired path is in thick red and blue, where blue edges  denote the sub-paths in the bypass graph.}
    \label{fig:case3b3}
\end{figure*}
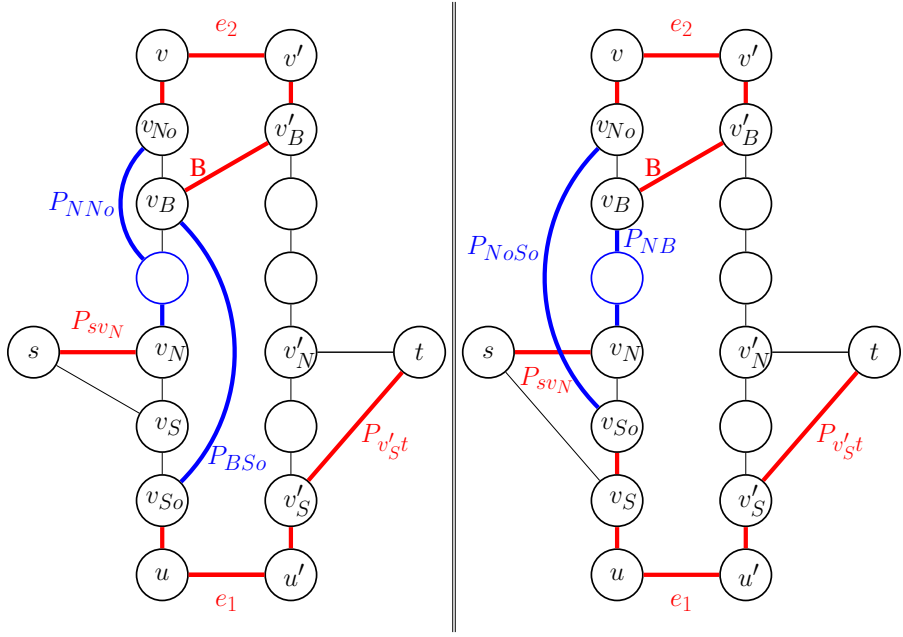

{\bf Case 3b3.I}: \emph{ There are  vertex-disjoint paths $P_{NNo}$} and $P_{BSo}$, see Figure 
\ref{fig:case3b3} (left).
Intuitively, here we have interleaved bypasses.
So, the required middle part bypasses the end-vertex $v_B$ of bridge $B$, goes to $(v,v')$ first, then crosses bridge $B$ from the right to the left, and then bypasses the external vertex $v_N$:
\[
(P_{NNo}, P_{uv}[v_{No}..v], (v,v'), P_{v'u'}[v'..v'_B], B,
P_{BSo}).
\]
%and then ends with: $(P_{vu}[v_{So}...u], (u,u'), P_{u'v'}[u'..v'_S],P(v'_S))$.

{\bf Case 3b3. II}: \emph{There are vertex-disjoint paths $P_{NoSo}$ and $P_{NB}$}, see Figure 
\ref{fig:case3b3} (right). 
Intuitively, this means
we have an overall bypass over the region
between external vertex $v_N$ and bridge
vertex $v_B$. So, the middle part crosses $B$ from the left to the right first, then takes $(v,'v)$, and bypasses
the bridge and external vertices:

\[
(P_{NB}, B, P_{u'v'}[v'_B..v'], (v',v), P_{vu}[v..v_{No}],
P_{NoSo}).
\]
%again ending with $(P_{vu}[v_{So}...u], (u,u'), P_{u'v'}[u'..v'_S],P(v'_S))$.

In both cases 3b3.I and 3b3.II, the constructed path is simple, since the bypasses we used are pure, and thus cannot (internal-vertex) intersect either $P_{sv_N}$ or $B$, 
and are completely disjoint from $P_{v'_St}$.
Note that Figure \ref{fig:case3b3} depicting case 3b3 (I and II) is a simplification with the bypass paths
shown as just one or two edges. In fact, these bypass paths can be quite convoluted, using edges of multiple bypasses interleaved with sections of $P_{uv}$.

As we have shown the desired $s,t$-path in all required cases, this concludes the proof of the lemma.
\end{proof}

Due to Lemma~\ref{lem:distinct}, the only remaining case of the \emph{if} part  of the proof of Theorem \ref{th:cycle}
that we have not considered is where some end-vertices of $e_1,e_2,e_3$ are not distinct. Let us analyze the possible alternatives, thereby
concluding the proof of Theorem \ref{th:cycle}:

\begin{enumerate}
    \item $e_1,e_2,e_3$ form a cycle: trivial.
    \item $e_1,e_2,e_3$ form a chain. W.l.o.g. let $e_1=(v_1,v_2),e_2=(v_2,v_3),e_3=(v_3,v_4)$. Let $G'=G\setminus \{v_2,v_3\}$, which is a connected graph since $G$ is 3-connected.
    Then, there is a path $P$ from $v_4$ to $v_1$ in $G'$, so $(e_1,e_2,e_3,P)$ is a simple cycle in $G$.
    \item Two edges have a common vertex, and one is disjoint. W.l.o.g., let $e_1=(v_1,v_2)$, $e_2=(v_2,v_3),e_3=(v_4,v_5)$, and $G'=G\setminus v_2$. Since $G$ is 3-connected, then $G'$ is 2-connected,
    and there exist vertex-disjoint paths $P_1, P_2$ from $\{v_4,v_5\}$ to $\{v_1,v_3\}$
    in $G'$.
    %Therefore, if  $P_1$ abuts $v_1$, then
    W.l.o.g., $P_1$ ends at $v_1$. Then,
    $(P_1,e_1,e_2,P_2,e_3)$ is a simple cycle in $G$. 
    %Otherwise, $(P_1,e_2,e_1,P_2,e_3)$ is a simple cycle in $G$.
    \item All edges share a common end-vertex. 
    This case is excluded by exception 1
    in the theorem statement.
\end{enumerate}

%******************************************
\section{CEP on Biconnected Graphs}\label{sec:biconnected}

In Section \ref{sec:basics}, we argued that 
in general undirected graphs, the CEP (and thus the PEP) can be easily reduced to biconnected graphs.
Above, we have determined necessary and sufficient conditions for the CEP in triconnected graphs.
In this section, we analyze the CEP on an arbitrary biconnected graph $G=(V,E)$, 
using its partition into triconnected components 
%As stated above, this can be done efficiently by generating the SPQR tree of $G$. --$>$
%using their partition into triconnected components 
provided by the SPQR tree ${\cal T} = {\cal T}(G)$.
For notation, definitions and properties of SPQR trees, refer to Section~\ref{sec:SPQR}.
We begin with defining the crucial concept of the component \emph{central} w.r.t. a given set of elements.

\begin{definition}[Central Component]
For any subset $S$ of elements of $G$, we call component $C$ of $\cal T$ \emph{central} w.r.t.\ $S$, if all elements $r_C(s), s\in S$, are distinct
$($that is, if
$\forall s_1, s_2 \in S: s_1\not = s_2 \Rightarrow r_C(s_1)\not = r_C(s_2))$.
\end{definition}

A central component does not necessarily exist in general; and if it exists, is not necessarily unique.
But for our purposes, it is sufficient that a central component always exists for sets of size 3.

\begin{lemma}[Central Component Lemma] \label{lem:central}
For every set $S$ of \emph{three} elements of $G$, there exists a component of $\cal T$ central w.r.t.\ $S$.
\end{lemma}

\begin{proof}
Denote the three given elements of $G$ by $x_1$, $x_2$, $x_3$. 
Let us choose three components $C_i$, $1 \le i \le 3$, such that $r_{C_i}(x_i)=x_i$, i.e. such that each $x_i$ is a \emph{real} element in $C_i$. 
%Observe that the central component w.r.t.\ $S$ in $\cal T$ can be recognized in the following structural way. 
%Consider three components $C(x_1)$, $C(x_2)$, $C(x_3)$ of $\cal T$. 
If at least two of these components coincide,
w.l.o.g., $C_1=C_2=C$,
then $C$ is central w.r.t.\ $S$. 
Indeed, by our assumption,
$x_1=r_C(x_1)$ and $x_2=r_C(x_2)$ are distinct non-virtual elements.
Additionally, either $r_C(x_3)=x_3$
or $r_C(x_3)$ is a virtual edge, which must in both cases be distinct from
$r_C(x_1)$ and $r_C(x_2)$.

Otherwise, %we have to assume that 
all three chosen components are distinct. If some $C_i$ lies on the path in $\cal T$ between the two other
components,
then the two elements $r_{C_i}(x_j)$, $j \neq i$, 
must be distinct: this is trivial if at least one of them is real, and
due to $C_i$ being on the path as above if both are represented by virtual edges.
Thus, $C_i$ is central w.r.t.\ $S$ here too. 

Otherwise, there exists a unique component, $C$, in tree $\cal T$ such that all $C_1,C_2,C_3$ are in distinct sub-trees hanging from $C$ in $\cal T$.
Let these sub-trees hang on structural edges $e_{\cal T}(e_i)$, respectively, where $e_i$ are distinct virtual edges of $C$.
Then, $C$ is central w.r.t.\ $S$ since each of the three $x_i$ is represented by 
%a distinct virtual edge 
$e_i$ in $C$, unless $x_i$ is a real element of $C$.
\end{proof}

    \remove{
\YD {Let us delete this paragraph, since the proof of the lemma is sufficiently algorithmic,(?)}
\ES{Maybe, but I think we are referring to this paragraph later (algorithms section) on so need to make sure this does not cause a problem later on.}
\YD {Check that the reference in line 4 of Section 5.2 is correct.}\ES{OK I think}
The following is another, algorithmic way to determine a component central w.r.t.\ $S$.
Let us begin from some component $C_1$ such that $r_{C_1}(x_1)=x_1$. If $x_2$ and $x_3$ are represented by distinct elements of $C_1$, 
then $C_1$ is central w.r.t.\ $S$. Otherwise, both $x_2$ and $x_3$ are represented by the same virtual edge, $e$, of $C_1$. 
Let $C$ be the other component that contains
virtual edge $e$.
Note that $e$ represents $x_1$ in $C$, but does not represent $x_2$ and $x_3$ in $C$. Once more, if $x_2$ and $x_3$ are represented by distinct elements of $C$, then $C$ is central w.r.t.\ $S$. Otherwise, we continue in the same way from $C$ as done in $C_1$ above.
This process is finite, since at each step, we move away from the originally chosen component $C$ in $\cal T$. The process stops at the component where $x_2$ and $x_3$ are represented by distinct elements, that is at a component central w.r.t.\ $S$.
}

We now state and prove a crucial proposition
on central components and must-include cycles.

\begin{proposition}
\label{p:central}
Let $S=\{x_1,x_2,x_3\}$ be a set of three elements of $G$, and $C$ be any component central w.r.t.\ $S$. Then, there exists a simple cycle containing all elements of $S$ in $G$ if and only if there exists a simple cycle containing all three elements representing them in $C$.
\end{proposition}

\begin{proof}

(\textbf{If})
Let $L'$ be a simple cycle in $C$
containing all $r_C(x_1),r_C(x_2),r_C(x_3)$.
If $L'$ contains no virtual edges, we are done since $L'$ is a cycle in $G$.
Otherwise, replace every virtual edge
$e\in L'$ by a path $P(e)$ as follows,
to create the simple cycle $L$ in $G$ as required.
Let virtual edge $e=(a,b)\in L'$.
%Let $(x,y)$ be any virtual edge in $L'$. 
If $e$ represents no element of $S$, then 
by Property~\ref{prop:biconnected},
there exists a simple path $P(e)$ from $a$ to $b$ in $G(C,e)$.
%\YD{ Delete: with all intermediate vertices and
%edges "represented by $(a,b)$ in $C$", that is such
%that for all $x\in P((a,b)), r_C(P)\in \{a, b, (a,b)\}$.}
(As a special case, path $P(e)$ may be the single edge $(a,b)$.)
Likewise, let $e$ represent element $x_i \in S$; then, $e$ represents no other element of $S$,
%there can be at most one such element,
since $C$ is a central component w.r.t. $S$.
By Property~\ref{prop:pathWvertex}, there exists a simple path from $a$ to $b$ containing $x_i$ in $G(C,e)$.
%\YD {Delete: with all intermediate represented by $(a,b)$ in $C$.}
By replacing every virtual edge $e$ in $L'$ by simple path $P(e)$, we obtain a cycle, $L$, in $G$ containing all elements of $S$.
Cycle $L$ is simple since graphs $G(C,e)$ as above are disjoint, except possibly for end-vertices of edges $e$, by Property~\ref{prop:sub-tree};
this exception cannot spoil simplicity of $L$, since $L'$ is simple.

(\textbf{Only if})
Let $L$ be a simple cycle in $G$ containing all three elements of $S$. 
%, and $C$ be a central component w.r.t.\ $S$.
If $L$ consists only of elements of $C$, it is a cycle in $C$, as required. 
Otherwise, cycle $L$ is sub-divided into the inclusion-maximal sub-paths whose edges are either all in $C$ or all in the same subgraph $G(C,(a,b))$, where $(a,b)$ is a virtual edge of $C$.
By Property~\ref{prop:sub-tree}, 
the end-vertices of each such sub-path in $G(C,(a,b))$ are $a$ and $b$. Therefore, replacing each such sub-path by virtual edge $(a,b)$ in $C$ results in a cycle, $L'$, in $C$. Cycle $L'$ is as required, since if $x_i$ was an element of $G(C,(a,b))$, then 
$r_C(x_i) = (a,b)$ is an edge of $L'$.
%\YD {DELETE:
%consider any inclusion-maximal sub-sequence of vertices $\tilde L$ in $L$ not in $C$. Since $C$ is a central component, $L$ must have at least 2 vertices in $C$, 
%\YD {I do not understand this "since". Can you explain it to me, Eyal?} %\ES{Not sure, I think you actually made this claim, so you will have to explain it... but it is needed as otherwise you cannot "splice"}
%and thus there exist distinct  vertices $a,b$ adjacent
%to vertices of $\tilde L$ on $L$. Then $(a,b)$ is a virtual
%edge in $C$.
%just before it and $y$ just after it in $L$ belong to $C$ and form a virtual edge $(x,y)$ in $C$.
%Moreover, if some element $z$ of $S$ is in
%$\tilde L$, then $(a,b)$ represents $z$ in $C$.
%If an edge $(a,b)$ in $L$ is not in $C$ but both $a$ and $b$ are in $C$, then $(a,b)$ is a virtual edge in $C$.
%By replacing in $L$ every maximal-inclusion vertex sub-sequence as above and every edge as above by the corresponding virtual edge in $C$, we obtain a simple cycle in $C$ containing all elements of $C$ representing the elements of $S$.}
\end{proof}

\emph{Remark\/}:
The reader may be curious about the relation of multiple central components to Proposition~\ref{p:central}.
In fact, absence of a cycle as required in a component, $C$, is a quite special situation.
This is impossible if $C$ is an S-node and happens in a P-node only if all all three its must-include elements are edges.
Since any R-node is triconnected, by Theorem~\ref{th:CEP}, %which implies that 
such a cycle can be absent in an R-node also only if all three must-include elements are edges.
%As stated in Lemma \ref{lem:unique}, in such a case,
By the following lemma, in the latter case, the central component should be unique.
%\YD {The place of Lemma 4.2 is changed, OK?}

%Note that for a given set $S$ of size 3, there may be more than one central component. However, in the following case we have a unique central component.

\begin{lemma}\label{lem:unique}
If all three elements of $S$ are represented by \emph{edges} in a component $C$ central w.r.t.\ $S$, then $C$ is a unique central component w.r.t.\ $S$.
\end{lemma}

\begin{proof}
Let distinct edges $e_1, e_2, e_3$ represent the %must-include 
elements $x_1, x_2, x_3 \in S$, respectively, in component $C$ central w.r.t.\ $S$. Consider another arbitrary component $C'$; %, which must have
necessarily, $C'\in B(C,e)$ for some virtual edge $e\in C$. 
Since $e_1, e_2, e_3$
are distinct edges, then only one of them can be the same as $e$, so suppose w.l.o.g. that
$e\not = e_1$ and $e\not = e_2$. This implies that
neither $x_1$ nor $x_2$ are in $B(C,e)$,
and therefore both elements $x_1, x_2$ are represented in $C'$ by the same virtual edge $e'$ (such that $C$ is in the sub-tree $B(C',e')$). This implies that $C'$ is not central w.r.t.\ $S$, as required.
%it is contained in the sub-tree of 
%$\cal T$ either represented by, w.l.o.g., virtual edge $e_3$, or by no one out of $e_1, e_2, e_3$. In either case, both 
\end{proof}

%\YD {DELETE:
%Let us summarize the remaining open issues for the PEP and CEP
%at this point. First,
%the PEP on graph $G$ is equivalent to the CEP on
%$G\cup \{(s,t)\}$, with the required elements in the cycle being $w_1,w_2,(s,t)$, i.e. two vertices and an edge.
%Thus it is sufficient address the CEP. Since there exists a cycle including three given elements in any graph $G$ if and only if there exists a cycle in a central component $C$ of ${\cal T}(G)$
%(Proposition \ref{p:central}),
%the CEP reduces to a 
%CEP of representing elements in $C$. Determining the latter is trivial in two cases: if $C$ is an $S$ node, the cycle obviously always exists. If $C$ is a $P$ node, the 
%cycle never exists if all three elements
%are represented by an edge in $C$, and always
%exists otherwise. The case where $C$ is an
%$R$ node is more complicated,
%and is handled by Theorem \ref{th:cycle}.
%We conclude this section by summarizing the results of the previous subsections as one theorem.}

We conclude this section by summarizing the results 
on must-include cycles from the
last two sections as one theorem.

\begin{theorem}\label{th:summary}
Let $G$ be a biconnected graph, and $x_1,x_2,x_3$ be three
distinct elements of $G$.
Let $C$ be any central component of ${\cal T}(G)$ w.r.t. $x_1,x_2, x_3$.
%with vertices $s,t$ such that $G^+=G\cup \{ (s,t)\}$ is biconnected, and $w_1,w_2$ be elements (vertices or edges) of $G$, all above elements distinct.
%$C$ be a central component of ${\cal T}(G)$ w.r.t. $x_1,x_2, x_3$.
%$(s,t),w_1,w_2$.
%Then, there exists a
Then, a simple  cycle in $G$ 
%includes $w_1$ and $w_2$ (equivalently, a simple cycle 
including $x_1,x_2, x_3$
%$(s,t),w_1,w_2$ in $G^+$)
exists if and only if 
either at least one of $r_C(x_1),r_C(x_2),r_C(x_3)$ is a vertex in $C$, or neither of the following conditions holds:
\begin{enumerate}
    \item Edges $r_C(x_1),r_C(x_2),r_C(x_3)$ %$r_C(w_1)$, $r_C(w_2)$, $r_C((s,t))$ 
    all have a common end-vertex in $C$.
    \item $C$ is an R-node and edges $r_C(x_1),r_C(x_2),r_C(x_3)$ %$r_C(w_1)$, $r_C(w_2)$, $r_C((s,t))$ 
    form a 3-edge cut of $C$.
\end{enumerate}
\end{theorem}

\begin{proof}
The theorem statement is consistent, since Lemma \ref{lem:central} guarantees existence of a central component $C$ w.r.t. $x_1,x_2,x_3$ in ${\cal T}(G)$. Additionally, if $C$ is not
unique, then by Lemma \ref{lem:unique},
$r_C(x_i)$ is a vertex for some
$3\geq i \geq 1$, in which case the cycle in $C$ always exists.
%\YD {I do not agree with the previous sentence here. It is not relevant to the proof of the theorem; its topic relates to the previous Remark (taken together with Lemma 4.2), and thus is appropriate to be placed there. Besides, to my mind, the statement of that sentence is not straightforward and thus should be proved.}
%\ES{I added it because a reviewer may say here "what if C is not unique? But if you think it is redundant, we can drop it.}
Due to Proposition \ref{p:central}, it is sufficient to prove that the conditions of the theorem are necessary and sufficient for existence of a simple cycle in $C$ that includes $r_C(x_1),r_C(x_2),r_C(x_3)$; observe that these three elements of $C$ are distinct since $C$ is central w.r.t. $x_1,x_2, x_3$.

Consider now all possible types
of node $C$. If $C$ is an $S$-node, the cycle as required trivially always exists. Accordingly, no three edges of an S-node (which is a cycle) share a common end-vertex, while exception 2 is not relevant.
If $C$ is a $P$-node, the cycle never exists if all three elements are represented by edges in $C$, which is covered by exception 1,
and always exists otherwise, which also accords with the theorem statement. 
In the case where $C$ is an $R$-node (which is a triconnected graph), % is more complicated.
the theorem follows immediately from %Lemma \ref{lem:central} (existence of a central component), %Proposition \ref{p:central} (reduction of problem to central component), 
Theorem \ref{th:CEP} (such a cycle always exists in a triconnected graph if not all elements are edges) and Theorem \ref{th:cycle} (the conditions for existence of a cycle with three must-include edges in triconnected graphs).
\end{proof}

%\noindent
%******************************

%\YD {Delete?: Note that exception 1 above holds for any three edges if the central component is a P-node, a fact used as a special case in the PEP algorithm in Section~\ref{sec:algorithms}.}
%The theorem also holds for the
%CEP, replacing $w_1,w_2$ by edges $e_1,e_2$ and with 
%$e_3=(s,t)$ being one of three required edges.

To illustrate the usage of Theorem~\ref{th:summary}, consider the graph shown in Figure \ref{fig:exampleSPQR} (top left).
In what follows, we use the reduction of the PEP to the CEP, the other parts of Figure \ref{fig:exampleSPQR} for illustration, and the proof of Proposition~\ref{p:central} to find the required must-include path, if it exists.
For simplicity, we choose  must-include path query examples where the source and destination vertices are already connected by an edge in the graph,
so that we can use the SPQR tree already shown in 
Figure \ref{fig:exampleSPQR}
(recall that the PEP query answer does not depend on presence/absence of such an edge).

Consider querying for a simple path from $w_4$ to $x$
that must include $u_1$ and $u_6$. Here, the equivalent 
CEP is
for elements $(w_4,x),u_1,u_6$,
for which the central component
is $C=R(u_1,u_2,u_4,u_5)$.
Cycle $L'=(u_1,u_2,u_5)$ in $C$ includes vertex $u_1$ and virtual edges $r_C((u_6))=(u_2,u_5) = e_1^{u_2,u_5}=Vir_5$ and $r_C((w_4,x))=(u_1,u_2) = e^{u_1,u_2}=Vir_6$.
By substituting edge $(u_1,u_2)$ with path $(u_1,x,w_4,t,u_3,u_2)$ and edge $(u_2,u_5)$ with path $(u_2,u_6,u_5)$, we get cycle $L=(u_1,x,w_4,t,u_3,u_2,u_6,u_5)$. By removing edge $(x,w_4)$ from $L$, we get simple path $(w_4,t,u_3$, $u_2,u_6,u_5,u_1,x)$, as required.

If queried for a must-include path from $w_4$ to $x$ through $u_1$ and $w_2$, the
equivalent
CEP query elements are
$(w_4,x),u_1,w_2$ and the
central component is $C=P(x,t)$. The representatives of $(w_4,x),u_1,w_2$ are three distinct virtual edges $e_i^{x,t}$, $1 \le i \le 3$. 
Since all these edges %connecting the same vertices $x$ and $t$ of $C$.
share both end-vertices,
no simple cycle including $(w_4,x),u_1,w_2$
exists, due to
exclusion 1 of Theorem~\ref{th:summary}. %, so no must-include cycle as required in the CEP problem induced by the query exist. 
Thus, the required must-include path also does not exist.

If queried for a must-include path from $u_2$ to $u_4$ via $w_4$ and $u_6$, the equivalent CEP
elements are $(u_2,u_4),w_4,u_6$
and the central component  is $C=R(u_1,u_2,u_4,u_5)$.
%Cycle $L'=(u_1,u_2,u_5,u_4)$ in $C$ goes through real edge $(u_1,u_4)$ and virtual edges $r_C(w_4)=(u_1,u_2)$ % = e^{u_1,u_2}=Vir_6$ 
%and $r_C((u_6))=(u_2,u_5)$. % = e_1^{u_2,u_5}=Vir_5$.
%By substituting the virtual edges as in the previous example, we get cycle $L=(u_1,x,w_4,t,u_3,u_2,u_6,u_5,u_4)$. By removing edge $(u_4,u_1)$ from $L$, we get path $(u_1,x,w_4,t,u_3,u_2,u_6,u_5,u_4)$  as required.
%If we change the source from $u_1$ to $u_2$, then the central component w.r.t. $(u_2,u_4),w_4,u_6$ remains the same. % is $C=R(u_1,u_2,u_4,u_5)$.
%As in the previous example, 
The triple of representatives $(u_2,u_4),(u_2,u_1),(u_2,u_5)$ of the must-include CEP elements in $C$ is covered by both exclusions 1 and 2 of Theorem~\ref{th:summary}, which implies that the answer to the query is negative.

%******************************************
\section{Algorithms for PEP/CEP and EPE}
\label{sec:algorithms}

We now turn to the algorithms for solving the CEP and EPE. %(and thus the PEP, see the reduction in Section~\ref{ss:basics}), and one for EPE.
Though both algorithms are stated to work for biconnected graphs,
the reductions provided in Section~\ref{ss:basics} imply similar results for general graphs.
To solve the PEP,
we use the CEP algorithm
and the reduction from PEP to CEP in Section~\ref{ss:problem statements}. 
%That is, given a PEP with $G, s, t, w_1,w_2$ solve
%the CEP on $G^+=G\cup (s,t)$
%with elements $(s,t), w_1, w_2$.

%*********************************
\subsection{Determining Must-Include Paths and Cycles}

We can now state Algorithm \ref{alg:CEP} for the CEP.
%(and the CEP as well, as a subproblem):
%given an undirected graph $G$, end-path vertices $s,t$ and intermediate vertices $w_1,w_2$, does there exist a simple path $P$ between $s$ and $t$ such that $w_1,w_2\in P$? 
Although an algorithm follows from the proofs
of Proposition \ref{p:central} and Lemma \ref{lem:central},
%(the algorithmic note after the latter proof \YD {To be edited}), 
here we state the algorithm explicitly and analyze its complexity.
%\YD {The function in the previous version of this algorithm did not solve the general CEP, just its case of one edge and two vertices.}

\begin{algorithm}[h]
	\caption{Determining Must-Include Cycle Existence (CEP)%, also solves the CEP.
	}
\hspace*{\algorithmicindent} 
\textbf{Input:} a biconnected graph $G$ and three of its elements $ x_1, x_2, x_3$ %, s, t, w_1, w_2$
	\begin{algorithmic}[1]
%	    \State Let $G^+=G\cup \{ (s,t)\} $
%	    \State {\bf return} \textproc{SolveCEP}($G^+, (s,t), w_1,w_2)$ 
%        \Function{SolveCEP}{$G, x_1, x_2, x_3$}
	    \State Compute $\cal T$, the SPQR tree of $G$
	   \State Compute representatives $r_C(x_1), r_C(x_2),r_C(x_3)$ for all components $C\in \cal T$ %\YD {choose T or $\cal T$ here and above}

  \State Find  a component $C$ central w.r.t. $x_1,x_2,x_3$
  \If {$C$ is an $S$ node or $r_C(x_i)$ is a vertex for at least one $3\geq i \geq1$}
      \State {\bf return} true
    \ElsIf {$C$ is a $P$ node}
         \State  {\bf return} false
    \EndIf
    \State /* Assert: $C$ is an $R$ node and all $r_C(x_i)$ are edges  */
    \If { $r_C(x_1), r_C(x_2), r_C(x_3)$ form an edge cut of $C$\\ 
    ~~~~~~~ or $r_C(x_1), r_C(x_2), r_C(x_3)$ all have a common end-vertex}
    \State {\bf return} false
    \Else
    \State {\bf return} true
    \EndIf
%   \EndFunction
    \end{algorithmic} %\vspace{-0.4cm}
    \label{alg:CEP}
\end{algorithm}

\begin{theorem}
Given any biconnected graph $G=(V,E)$,
%with vertices $s,t$  such that $G^+$ is biconnected,
Algorithm \ref{alg:CEP} always answers the CEP correctly, 
and has time complexity $O(|E|)$.
\end{theorem}

\begin{proof}
Correctness of Algorithm \ref{alg:CEP} follows from Theorem \ref{th:summary}.
The complexity of this algorithm is dominated by computing the SPQR tree of $G$, which is $O(|E|)$.
Computing the representatives can be done by traversing $\cal T$ from any of its vertices in $O(|E|)$ time as well.
It is clear from examining the rest of the steps of Algorithm \ref{alg:CEP}
that they all take time $O(|E|)$,
with some of the steps in fact only requiring constant time.
\end{proof}

Due to the equivalence between PEP and CEP stated in Section~\ref{sec:basics},
%for solving 
PEP$(G,s,t,$ $w_1,w_2)$ can be solved by executing Algorithm \ref{alg:CEP} for CEP$(G \cup (s,t),(s,t),w_1,w_2)$,
with the same time bound.
%and remove edge $(s,t)$ from the outputted cycle.

\subsubsection{Finding the Cycle}

It may be of interest to \emph{actually find} a simple cycle that includes the given elements $x_1, x_2, x_3$
%from $s$ through $w_1,w_2$ to $t$
when it exists. Naturally, we begin by 
constructing ${\cal T}(G)$ and finding the representatives of $x_1, x_2, x_3$ and
a central component $C$, as
in Algorithm \ref{alg:CEP}.
Then, consider any element $x_i$
that is not a real element in $C$.
Finding the portion of the required cycle
outside $C$ that contains $x_i$
entails finding a path between the ends of the virtual edge $r_C(x_i)$ outside $C$ that contains $x_i$. This can be done in linear time using the Ford-Fulkerson algorithm by finding
two vertex-disjoint (except at $x_i$, if $x_i$ is a  vertex) paths
between $x_i$ and the ends of $r_C(x_i)$
in the subgraph $G(C,r_C(x_i))$ (see Section~\ref{sec:basics}).

Finding the part of the cycle in $C$ that includes all the $r_C(x_i)$ (to be spliced with the paths as above)
is trivial for S-nodes and P-nodes. 
If $C$ is an R-node and at least one of $x_i$ is a real vertex of $C$, this can be done as described in Section~\ref{sec:PEPknown}.
The non-trivial task is finding the portion of the cycle as above in the central component  when $C$ is an R-node and all $r_C(x_i)$ are edges.
Here, there are  many cases to
handle, so we specify this part of the algorithm only implicitly. 
Observe that the proof of Theorem \ref{th:cycle},
including its required lemmas, 
specifies the required cycle in a constructive way, except for finding paths as in Menger's theorem; the latter is covered by Section~\ref{ss:basics}.
Thus, by following the cases
of the proof, we (implicitly) get an $O(|E|)$
algorithm for finding a simple $x_1,x_2,x_3$-cycle.
By the reduction in Section~\ref{sec:basics}, this also implies an $O(|E|)$
algorithm for finding a simple 
$s,\{w_1,w_2\},t$-path.

%\YD {DELETE:
%is in most cases constructive, and specifies the required path from $s$ to $t$ (that contains the edges representing $w_1$ and $w_2$) directly.
%The only parts in the proof that are not
%constructive are on finding the paths 
%between vertex sets, from vertices
%to the cycle $L$ used in creating the
%path, and the connectors of various types
%used in the proof (bridges, external connectors, and bypasses).
%In some cases, we need a specific number of paths,
%such as when determining the paths in Proposition
%\ref{prop:tri+cycle}, or when we need to find %the bypass paths in case 3b3 of Lemma \ref{lem:distinct}.
%It is well known that such paths can be found in time $O(|E|)$ by using the Ford-Fulkerson algorithm \cite{10.5555/1942094}.
%and the flow decomposition (\YD {To be checked in Even's book} folklore, see, e.g., in \cite{DDK}).
%In other cases, we need to determine
%path existence and can use BFS, including:
%finding out whether the condition
%for Lemma \ref{lem:distinct} requiring
%no $s,t$-path in $G\setminus L$ holds,
%finding the northmost and southmost
%$s$-connector and $t$-connector in case 3 of 
%Lemma \ref{lem:distinct},
%finding the southmost and northmost bridges.
%Thus, by following the cases
%of the proof we (implicitly) get an $O(|E|)$
%algorithm for finding a simple $s,\{w_1,w_2\},t$-path.
%}

%**************************************
\subsection{Enumerating Exclusion Pairs}

Consider now the EPE problem. 
Recall the notation $G^+=G\cup (s,t)$.
The idea here is to visit every component $C'$ that could potentially be central in ${\cal T} (G^+)$ w.r.t. $(s,t)$ and some two vertices not in $C'$, 
and for each such component emit all
possible vertex pairs 
excluded from being on a simple $s,t$-path by Theorem \ref{th:summary}.

Algorithm \ref{alg:EPE} for the EPE
starts with the component $C$
which has $(s,t)$ as a real edge, 
assigning $C$ to be
the root of ${\cal T}(G^+)$, and visits
components $C'$ recursively
using procedure \textproc{Traverse}.
In each such $C'$,
we consider all virtual edges other than
$r_{C'}((s,t))$ as potentially generating exclusions; 
we use $E^{\mbox{virt}}(C^-)$ to denote the set of virtual edges in $C^-= C' \setminus r_{C'}((s,t))$.
%\YD {There are $\{$ and $\}$ here, in line 7 of Algorithm 2, in the third line before Theorem 5.2 and in Theorem 5.2. This does not accord the notation in line 1 of Algorithm 2 and, as far as I recall, in many places in the rest of our paper after $\setminus$ and $\cup$. Let us try being consistent.}\ES{Do we want to say that we omit the $\{ ~~\}$ for singleton sets (and then can drop them always), or just be precise and always use the $\{ ~~\}$, and add them everywhere? What is your preference? But also note that strictly speaking $G\cup \{ (s,t)\}$
%is also problematic, as $G$ is not a set.
%This is shorthand for: $(V, E\cup\{ (s,t)\})$
%where $G=(V,E)$. But perhaps this is so standard that we don't have to say it?}
%\YD {Thank you for pointing out that also $G\cup \{ (s,t)\}$ is problematic. So, I'd certainly prefer omitting $\{$ and $\}$ everywhere for singleton sets.}
Recall that $r_{C''}(C')$, the representation of $C'$ in its child $C''$ is the virtual edge in $C''$  corresponding to the structural edge between $C'$ and $C''$.
This allows us to cheaply supply $r_{C''}((s,t))=r_{C''}(C')$ as an additional argument in
the recursive call at the end
of procedure \textproc{Traverse}.
%allows us to cheaply supply $r_{C''}((s,t))$ as an additional argument in the recursive call.}
For convenience, we use $s',t'$ to denote the ends of $r_{C'}((s,t))$ in $C'$.

\begin{algorithm}[h]
	\caption{Exclusion-Pairs Enumeration (EPE)}
	\hspace*{\algorithmicindent} \textbf{Input:} a biconnected graph $G$ and two of its vertices $s, t$
	\begin{algorithmic}[1]
	    \State Let $G^+=G\cup (s,t)$
	    \State Compute ${\cal T=T}(G^+)$, the SPQR tree of $G^+$
    \State Root $\cal T$ at $C$, the unique component of $\cal T$
       where     $r_C((s,t))=(s,t)$ is a real edge
    \State Call \textproc{Traverse}($C,(s,t)$)
    \State {\bf exit}
    \Procedure{Traverse}{$C',(s',t')$}
       \State Let $C^- = C'\setminus  (s',t') $
       \If {$C'$ is a $P$ node}
        \If{$|E^{\mbox{virt}}(C^-)| \ge 2$}
            \State \textproc{EmitPairs}($C',E^{\mbox{virt}}(C^-)$)
        \EndIf
      \ElsIf{$C'$ is a non-leaf $R$ node}
%      \State Let $Cuts=$Find3EdgeCuts$(C',(s',t')))$
%      \For {all $E_{EX} \in Cuts$

     \For {$u\in \{s',t'\}$}
      \State Let $E_{EX}= \{ (u,v)|(u,v) \in E^{\mbox{virt}}(C^-) \}$
        \If{$|E_{EX}| \ge 2$}
          \State \textproc{EmitPairs}($C',E_{EX}$) \label{EP:ends}
        \EndIf
    \EndFor
     \State Let $Cuts = $ \textproc{Find2EdgeCuts}$(C^-,(s',t'))$
     \For {all $E_{EX} \in Cuts$ 
     s.t. $E_{EX}\subseteq E^{\mbox{virt}}(C^-)$
       \\ ~~~~~~~~~~~~~~~~~~~~~~and
      the edges in $E_{EX}$ are not both incident on either $s'$ or $t'$      }
         \State \textproc{EmitPairs}($C', E_{EX}$) \label{EP:cut}
      \EndFor
  \EndIf
  \For {all $C''\in \mbox{children}(C')$}
      \State Call \textproc{Traverse}($C'',r_{C''}(C')$)
  \EndFor
\EndProcedure
    \end{algorithmic}
    \label{alg:EPE}
\end{algorithm}

In each component $C'$, given some set $E'$ 
of virtual edges in $C'$, we use function \textproc{EmitPairs}($C',E'$) to
emit exclusion pairs of vertices from $G$, as follows.
For a virtual edge $e$ in $C'$, 
denote by $V(C',e)$ the set of 
vertices (from $G$) in $G(C',e)$ 
excluding the end-vertices of $e$. 
%\YD {nevertheless, the word "appearing" should be explained in the following text from Section 2: vertices and edges appearing in T'.}  STILL TO DO!!!!! **** !!!
Then, for every pair of distinct virtual edges $\{ e_1,e_2\}\subseteq E'$, the exclusion vertex pair $\{v_1,v_2\}$ is emitted for every
pair of vertices $v_1 \in V(C',e_1)$ and $v_2 \in V(C',e_2)$.
\textproc{EmitPairs}($C',E'$) can work either implicitly,
 i.e., just output pointers to $C'$ and
 to the set of virtual edges $E'$,
 or explicitly, i.e., output each exclusion
 vertex pair $\{ v_1, v_2\}$
 explicitly.

\textproc{EmitPairs} is called exactly as
indicated by Theorem \ref{th:summary}.
The S-nodes do not cause exclusion pairs,
because they are cycles and
thus there always trivially exists a cycle containing any
set of elements therein, so the algorithm does not
call \textproc{EmitPairs} for S-nodes.
If $C'$ is a P-node, all pairs of virtual edges in $C^-$ cause exclusions, which are emitted.

If $C'$ is an R-node, exclusions are
emitted (in line \ref{EP:ends}) for every set of cardinality at least 2 of virtual edges of $C^-$  incident on either end of $r_{C'}((s,t))$, as well as (in line \ref{EP:cut}) 
for every pair of virtual edges forming a 3-edge cut of $C'$ together with edge $(s',t')$,
unless these edges are both incident on either $s'$ or on $t'$, in which case we do not need to emit them again (already emitted in line \ref{EP:ends}).
%\YDT {Note that for any graph $C'$ with edge $(s',t')$, the 3-edge cuts separating $s'$ and $t'$ are exactly the 2-edge cuts separating $s'$ and $t'$ in graph $C^- = C' \setminus \{(s',t')\}$ with the extra edge $(s',t')$.}
Note that $\tilde{E}$ is a 2-edge cut separating
$s'$ and $t'$ in $C^- = C'\setminus (s', t')$ if and only if
$\tilde{E}\cup (s',t')$ is a 3-edge cut separating
$s'$ and $t'$ in $C'$.
%\YD {I like your version, Eyal. Just let us avoid E' (used above with another meaning. Say, $\tilde E$. Besides, maybe recall $C^- = C' \setminus \{(s',t')\}$.}
We show how to find all such
2-edge cuts effectively in Section~\ref{sec:EdgeCuts},
thereby describing function \textproc{Find2EdgeCuts}.

\begin{theorem}
Given any graph $G=(V,E)$ with vertices $s,t$ such that $G^+ = G \cup (s,t)$ is biconnected, Algorithm \ref{alg:EPE} always emits the correct exclusion pairs. 
The algorithm has time complexity $O(|E|)$ for implicit output and  $O(\max\{|E|,N\}) = O(|V|^2)$ for explicit output, where $N$ is the number of emitted exclusion pairs.
\end{theorem}

\begin{proof} (outline)
Correctness of Algorithm \ref{alg:EPE} is shown as follows. First, we claim that
whenever \textproc{Traverse}$(C',(s',t'))$ is called, we have $(s',t')=r_{C'}((s,t))$,
which is shown by trivial induction.
Next, let $w_1, w_2$ be an exclusion pair w.r.t. $s,t$. 
Then, by Theorem \ref{th:summary}, there exists a component $C$ of type P or type R
central w.r.t. $w_1,w_2,(s,t)$.
In the former case (type P), 
$r_C(w_1), r_C(w_2),r_C((s,t))$
are distinct virtual edges (except for $r_C((s,t))$, which may be a real edge) in $C$.
%\YD {Delete: (otherwise a cycle $w_1,w_2,(s,t)$ exists)}. 
All P nodes are visited and all 
exclusion pairs w.r.t.
such edge sets are emitted.
In the latter
case (type R), $r_C(w_1), r_C(w_2),$ $r_C((s,t))$
are all virtual edges (except for $r_C((s,t))$, which may be a real edge), such
that either they all meet at one end of $r_C((s,t))$, or they form a 3-edge cut. Since 
%\YD {I do not understand the expression: all these cases also specifically emit an exclusion}\ES{Is this better: \textproc{Traverse} explicitly emits exclusions in all these cases},
\textproc{Traverse} explicitly emits exclusions in all these cases,
the algorithm emits all exclusion pairs. 
The converse follows since the above are the only exclusions emitted by the algorithm.

The total time for all \textproc{Find2EdgeCuts} calls is $O(|E|)$ (see Section \ref{sec:EdgeCuts}).
By examination of the pre-order traversal of the components of $\cal T$, the $O(|E|)$ complexity for \emph{implicit output} follows
from each edge in components of $\cal T$ being examined at most a constant number of times
and the $O(|E|)$ size of the entire SPQR tree data structure \cite{SPQRtrees} (see 
%Theorem~\ref{pr:SPQR-size-time} (from 
Section~\ref{sec:SPQR}).
%\YD {Do we need this addition, Eyal? Or the entire style of the paper is relying on those properties without mentioning them?}
%\ES{In Section 4 which you have yet to examine
%these properties are used in several locations.}
%\YD {I meant just a few words added by me, now distinguished in blue, see above.}

For \emph{explicit output},
the required list of vertex pairs can be generated in $O(\max\{|E|,N\})$ using the following lemma:

\begin{lemma}
The number of nodes in the sub-tree of ${\cal T}$ hanging on any virtual edge $e$ in any component $C_0$ is at most twice the number of internal
(that is, excluding the end-vertices of $e$)
vertices of $G$ in $V(C_0,e)$. \end{lemma}

%\YD {I do not understand: is $C$ the same (unique) $C$ as in Algorithm 2 or it is an arbitrary component. In the latter case, it is worth to choose another notion instead of $C$. In any case, we should resolve this ambiguity explicitly.} \ES{You are right, it is not the same C. What should we call it (note we already have C' in the proof. Should we do C' and C'' or is that too many primes?)}

\begin{proof} Consider the scan in ${\cal T}$ starting from virtual edge $e$ in $C_0$. 
Each descendant component $C^\prime$ contains at least one vertex of $G$ not counted
before visiting $C^\prime$, except for P-nodes. Observe that any P-node $C^\prime$ has 
least one child and has no children of type P. Hence, we 
can choose anyone of its children and charge it for $C^\prime$. The lemma follows.
\end{proof}

Due to the lemma,
assuming the list of vertices of $G$ in each component $C'$ is explicitly stored, 
the complexity of creating any list $L_{C'e}=V(C',e)$ via the subtree scan is thus linear in  the list size, $O(|V(C',e)|)$. 
Therefore, even using the straightforward scheme of creating the lists of vertices $L_{C'e}$ for all $e\in E'$
whenever \textproc{EmitPairs}($C',E'$) is called, and then
simply iterating over the pairs of lists to produce the output, takes time linear in the number of vertex pairs emitted.
 \end{proof}

\textit{Remark\/}: Note that Algorithm \ref{alg:EPE} can easily be modified
to emit all \emph{excluded pairs of edges}, rather than vertices. 
Simply omit the conditions checking whether the edges in
$E_{EX}$ are virtual, allowing real edges  to be emitted as well.
The semantics of \textproc{EmitPairs} is then
changed to emit the element pairs using the set of {\em edges} in each $G(C',e)$
instead of the set of vertices therein.
%\YD {I did not find where notion $BC(T,C')$ was defined. What is it?}\ES{It is a typo, of course... sorry. But also we have a new notation that is clearer so updated above) Then I rephrased everything below, is this better?}
%\YD {Now, quite OK with me.}
Additionally, when computing the set of edges for
the 2nd argument ($E_{EX}$) of \textproc{EmitPairs},
real edges are also added: whenever
Algorithm \ref{alg:EPE} indicates $E^{\mbox{virt}}(C^-)$ this is replaced by
\emph{all} edges of $C^-$.
Then, \textproc{EmitPairs} interprets a real edge
$e$ in its $E_{EX}$ argument as an exclusion
involving $e$ itself, rather than the reference to the edges in $G(C',e)$.
%whenever $u,v$ is a
%real edge in $E_{EX}$ then
%\YD {I cannot decipher the rest of the sentence (though do understand what is meant there).}
%$(u,v)$ itself is
%the one participating in the exclusion
%with the other edges in $E_{EX}$.
Generalizing to all excluded pairs of elements (vertices and edges) can be done similarly.

\subsection {Finding the 2-Edge Cuts}
\label{sec:EdgeCuts}

We now turn to the only remaining part of our solution to EPE:
finding  all 2-edge cuts between $s'$ and $t'$ in $C^-=C\setminus (s', t') $, where $C$ is a  central component in $\cal T$ of type R found by Algorithm~\ref{alg:EPE} (see the call to procedure \textproc{Find2EdgeCuts} therein).
In a general graph, 
it is not so trivial to list all $s,t$-mincuts, and
there might be exponentially many of them.
%the edge %\YDT {two terminals} 
%min-cuts efficiently, \YD {I'd like to restore "e.g.," here, since this is not the unique problem; the other one is a need to create a recursive program for listing all the cuts.} because there might be exponentially many such cuts. 
However, our case is quite specific, allowing ''leaving behind'' both edges of the last found 2-edge cut before continuing the search for the next 2-edge cut, as discussed below; this allows for a simple and efficient algorithm.
In what follows, we use the notation and knowledge on  basics of graph connectivity and related algorithms presented in Section~\ref{ss:basics}. % and \ref{ss:basics-alg}.}

\begin{figure*}
    \centering
\resizebox{1.0\columnwidth}{!}{
\begin{tabular}{c||c}
\begin{tikzpicture}[minimum size=13mm,
  node distance=0.5cm and 2.5cm,
  >=stealth,
  bend angle=45,
  auto]
  \tikzstyle{every node}=[font=\Huge]
\node (0) [circle, ultra thick, text centered, text width=0.5cm, draw=black] {$s$};
\node (1) [circle, ultra thick, text centered, text width=0.5cm, draw=black, right=of 0] {}
edge [ultra thick, color=green] (0);
\node (2) [circle, ultra thick, text centered, text width=0.5cm, draw=black, above=of 1] {} edge [ultra thick] (1)
edge [ultra thick, color=green] (0);

\node (10) [circle, ultra thick, text centered, text width=0.5cm, draw=black, right=of 1] {}
edge [ultra thick, color=red] (1);
\node (11) [circle, ultra thick, text centered, text width=0.5cm, draw=black, above=of 10] {}
edge [ultra thick] (10)
edge [ultra thick, color=red] (2);
\node (12) [circle, ultra thick, text centered, text width=0.5cm, draw=black, right=of 10] {}
edge [ultra thick] (10)
edge [ultra thick] (11)
;

\node (17) [circle, ultra thick, text centered, text width=0.5cm, draw=black, right=of 12] {}
edge [ultra thick, color=blue] (12);
\node (18) [circle, ultra thick, text centered, text width=0.5cm, draw=black, above=of 17] {} edge [ultra thick] (17)
edge [ultra thick, color=blue] (11);
\node (20)
[circle, ultra thick, text centered, text width=0.5cm, draw=black, right=of 17]
{$t$}
edge [ultra thick] (17)
edge [ultra thick] (18)
edge [dashed,bend left] (0)
;
\node (21)
[circle, ultra thick, text centered, text width=0.5cm, draw=black, above=of 18]
{$~$}
edge [ultra thick] (18)
edge [ultra thick] (20)
edge [ultra thick,bend right] (17)
;

\node (30) [above=of 0] {$G(V_0)$};

\node (31) [above=of 2] {$G(V_1)$};

\node (32) [above=of 11] {$~~~~~~~~~~~~~G(V_2)$};

\node (33) [above=of 21] {$~~~~~~~~~~~~~G(V_3)$};

\end{tikzpicture}

&

\begin{tikzpicture}[minimum size=13mm,
  node distance=0.5cm and 2.5cm,
  >=stealth,
  bend angle=45,
  auto]
  \tikzstyle{every node}=[font=\Huge]
\node (0) [circle, ultra thick, text centered, text width=0.5cm, draw=black] {$s$};
\node (1) [circle, ultra thick, text centered, text width=0.5cm, draw=red, right=of 0] {}
edge [ultra thick, color=green] (0);
\node (2) [circle, ultra thick, text centered, text width=0.5cm, draw=black, above=of 1] {} edge [ultra thick] (1)
edge [ultra thick, color=green] (0);

\node (10) [circle, ultra thick, text centered, text width=0.5cm, draw=black, right=of 1] {}
edge [ultra thick, color=blue] (1);
\node (11) [circle, ultra thick, text centered, text width=0.5cm, draw=black, above=of 10] {}
edge [ultra thick] (10);
\node (12) [circle, ultra thick, text centered, text width=0.5cm, draw=black, right=of 10] {}
edge [ultra thick] (10)
edge [ultra thick] (11)
;
\node (13) [circle, ultra thick, text centered, text width=0.5cm, draw=black, right=of 11] {}
edge [ultra thick] (10)
edge [ultra thick] (11)
edge [ultra thick] (12)
;
\node (30) [below=of 10,xshift=2cm] {$G(V_i)$};

\node (17) [circle, ultra thick, text centered, text width=0.5cm, draw=red, right=of 12] {}
edge [ultra thick, color=blue] (12);
\node (18) [circle, ultra thick, text centered, text width=0.5cm, draw=black, above=of 17] {} edge [ultra thick] (17)
edge [ultra thick, color=blue, bend right] (2);
\node (20)
[circle, ultra thick, text centered, text width=0.5cm, draw=black, right=of 17]
{$t$}
edge [ultra thick] (17)
edge [ultra thick] (18)
edge [dashed, bend left] (0)
;
\node (21)
[circle, ultra thick, text centered, text width=0.5cm, draw=black, above=of 18]
{$~$}
edge [ultra thick] (18)
edge [ultra thick] (20)
edge [ultra thick,bend right] (17)
;
\end{tikzpicture}

\end{tabular}
}
    \caption{Three 2-edge cuts sequence in $G^-$ (different colors, left). Overlap of 2-edge cuts in $G^-$ (blue) implies a 2-edge cut in $G$ and, thus, a 2-vertex cut in $G$ (red vertices, right).}
    \label{f:2-cuts}
\end{figure*}
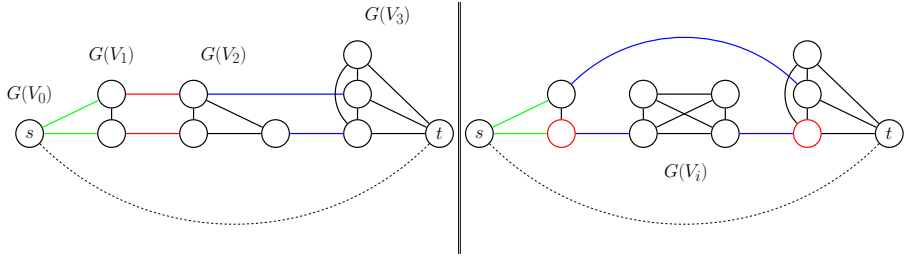

    \remove{
\YD {nontrivial --$>$ following}
\YD {There are words "not so trivial" in line 4}
part: when any central component
$C$ being visited is an R-node, 
we need to find all possible
2-edge cuts between $s'$ and $t'$ in 
$C^-=C\setminus \{ (s', t')\} $.
This is done by
function \textproc{Find2EdgeCuts} (Algorithm \ref{alg:cuts}).
% \textit{Finding 3-Edge Cuts between $s$ and $t$} works. [[ In this text, we will write ``cut'' meaning ``edge cut.(???)]] 
In general, it is not so trivial to find all the edge min-cuts efficiently. However, our case is quite specific, allowing ''leaving behind'' both edges of the last found 2-edge cut before continuing the search for the next 2-edge cut. For illustration, Figure~\ref{f:2-cuts} (left) shows the situation that we use in the routine. Let $E_i$ be an $s',t'$ cut,
denote by $V^{s'}_i$ the set of vertices in side of the cut that contains $s'$,
and by $V^{t'}_i$ the vertices on the $t'$ side.
Formally, we have:
}

Let us analyze the structure of 3-edge $s,t$-cuts %separating two given vertices $s$ and $t$ (called 3-edge $s,t$-cuts) 
in an arbitrary triconnected graph
%, changing a bit the notation for simplicity.
%Consider a triconnected graph 
$G=(V,E)$, $|V| \ge 4$, such that $(s,t) \in E$.
%\YD {See the e-mail sent to you now.}
%Assume, to the contrary, that $\{(u_1,v_1), (u_2,v_2)\}$ is such a 2-cut disconnecting $u_1,v_1$ from $u_2,v_2$.
%If all $u_1,v_1,u_2,v_2$ are distinct, then two vertices $u_1$ and $v_2$ disconnect $v_1$ from $u_2$, a contradiction to triconnectivity of $G$.
%Otherwise, assume, w.l.o.g., that $u_1=v_1$ and $u_2 \neq v_2$. Since $|V| \ge 4$, either $\{u_1\}$ or $\{u_2,v_2\}$ are vertex cuts, a similar contradiction.
%Likewise, there is no 1-edge cut in $G$, since either one or the other end-vertex of a single edge disconnecting $G$ forms a vertex cut alone.
%only if that minimum cut cardinality is exactly 3. 
%In what follows, we are interested in edge cuts of either $G$ or $G^-$ that disconnect $s$ and $t$, called $s,t$-cuts. 
Recall that any (minimal) edge cut divides $G$ into exactly two connected induced subgraphs, whose vertex sets form a 2-partition of $V$.
For any 
% We assign to 
$s,t$-cut $E'$, we denote the corresponding 2-partition of $V$ by $(V^s(E'), V^t(E'))$, $s \in V^s(E')$, $t \in V^t(E')$. 

By Theorem~\ref{l:edge-vertex cuts}, there is no 1- or 2-edge cut in $G$. 
Henceforth in our analysis, we assume that there exist %$q \ge 1$ 
3-edge $s,t$-cuts in $G$;
hence, those cuts are minimum $s,t$-cuts in $G$. 
By Theorem \ref{th:crossing}, these edge cuts do not cross each other.
%(see Section~\ref{ss:basics}).
This means that for any two such edge cuts $E'$ and $E''$, 
either $V^s(E') \cap V^t(E'')$ or $V^s(E'') \cap V^t(E')$ is empty;
equivalently, either $V^s(E') \subset V^s(E'')$ or $V^s(E'') \subset V^s(E')$.
Therefore, if there are $q$ 3-edge $s,t$-cuts in $G$, then we can order them as $E_i$, $1 \le i \le q$, so that for any $1 \le i <j \le q$, $V^s(E_i) \subset V^s(E_j)$
(a nested sequence of 2-partitions).
%, where $V^s_k$ denotes $V^s_k$
We denote: $V_0 = V^s(E_1)$, $V_q = V^t(E_q)$, and, for any $1 \le i \le q-1$, $V_i = V^s(E_{i+1}) \setminus V^s(E_i) \neq \emptyset$.

Denote $G^- = (V,E \setminus (s,t))$.
Obviously, each 3-edge $s,t$-cut $E'$ in $G$ contains edge $(s,t)$ and, thus, naturally defines the 2-edge $s,t$-cut $E' \setminus (s,t)$ of $G^-$ dividing $V$ in the same way.
%Let us now flip to graph $G^-$
That is, for every $i$, there are exactly two edges between $V^s(E_i)$ and $V^t(E_i)$ %, in addition to $(s,t)$, 
in $G^-$. 
Let us show that both those edges are between $V_{i-1}$ and $V_i$, for all $i$.
Assume, to the contrary, that there is an edge between $V^s(E_i)$ and $V^t(E_{i+1})$, ``jumping over'' $V_i$
and thus common to $E_i$ and $E_{i+1}$ (as well as, possibly, to other $E_j$).
A simple count shows that there is at most one edge between $V^s(E_i)$ and $V_i$ and at most one edge between $V_i$ and $V^t(E_{i+1})$. That is, there are overall at most two edges connecting $V_i$ to the rest of $G$, a contradiction to the absence of 1- and 2-edge cuts in $G$.
For illustration, see Figure~\ref{f:2-cuts} (right). 
Summarizing, graph $G^-$ has a very special structure: a sequence of $q+1$ induced subgraphs $G(V_i)$ with exactly two edges between any two subsequent subgraphs, and no other edges between these subgraphs.
See Figure~\ref{f:2-cuts} (left) for an example.

We now use the flow techniques:
construct flow network $N$ from graph $G^-$ by replacing each edge $(u,v)$ therein by two directed edges $(u,v)$ and $(v,u)$ of capacity 1 each and by assigning $s$ as its source and $t$ as its sink.
Recall that by Theorem~\ref{th:max-flow min-cut},
%the max-flow min-cut theorem of Ford and Fulkerson \cite{10.5555/1942094}, 
the size of any maximum flow in $N$ equals the size of any minimum (edge) $s,t$-cut. 
Thus, there is no 3-edge $s,t$-cut in $G$ if and only if there is no 2-edge $s,t$-cut in $G^-$ and if and only if the max-flow size in $N$ is at least 3.

We consider the case where the max-flow size in $N$ is 2 and the $q$ cuts $E_i \setminus (s,t)$ are the minimum 2-edge $s,t$-cuts (``min-cuts'') in $N$.
Let $f$ be any maximal flow of size 2 in $N$, and $N_f$ be the residual network w.r.t. $f$.
Any edge $(u,v)$ in $N$ is absent in $N_f$ if and only if it is saturated: $f(u,v)=1$.
Note that both edges of any cut $E_i \setminus (s,t)$ are saturated by flow $f$ in the direction from $V_{i-1}$ to $V_i$ in $N$, so that their two inverse edges (directed from $V_i$ to $V_{i-1}$) are the only edges connecting $V_{i-1}$ and $V_i$ in $N_f$.
Picard and Queyranne \cite {PQ} proved that all min-cuts together divide the vertices of $N_f$ into the vertex sets of its strongly connected components. 
That is, in our case, the strongly connected components of $N_f$ are the induced graphs $G(V_i)$.

Recall that each iteration of the Ford-Fulkerson algorithm %\cite{10.5555/1942094} 
executes the vertex scan (say, BFS, to be specific) from $s$ in the current residual network. If $t$ is labeled, then the flow is incremented. In our case, the flow size increases by exactly one at each iteration, so if $t$ is labeled by the third BFS, then the max-flow size is at least 3, and thus, there is no 2-egde $s,t$-cut in $G^-$. 
Otherwise, the flow $f$ constructed after the first two iterations is maximum. Thus, the BFS scan 
%built after exactly two iterations. 
%At the third iteration, BFS is 
executed in $N_f$ labels exactly the vertices of $V_0$ and then stops, since the orientation of the edges of $E_1$, the only edges connecting $V_0$ to its outside in $N_f$, prevents BFS from advancing. We thus are able to detect the edge cut $E_1$ as the two edges whose one end-vertex is labeled and the other is not.
If after that we label artificially the end-vertices of edges in $E_1$ external to $V_0$, i.e. those in $V_1$, and let BFS continue, then the BFS would label exactly the vertices of $V_1$, stopping similarly. Then, we can find the edge cut $E_2$, and so on up to finding $E_q$; thereafter, BFS labels all vertices of $V_q$, including $t$. 

See procedure \textproc{Find2EdgeCuts} (Algorithm~\ref{alg:cuts}) implementing this scheme, where the notation is changed to that compatible with the main algorithm solving EPE: %$G \rightarrow C$, 
$G^- \rightarrow C^-$, $s \rightarrow s'$, $t \rightarrow t'$.
%; its correctness is based above.
Some technical trick is added to the usual BFS scan for detecting the edges of each next min-cut. % by the way.
Instead of scanning only the edges of $N_f$ out-going from the currently labeled vertex $v$, we scan {\em all} the edges of $C^-$ (equivalently, the edges of $N$) incident on $v$. 
Such an edge $(v,v')$ can only be in the next min-cut if it has flow 1 assigned (and thus is absent in $N_f$) and if $v'$ is  currently  unlabeled. In each such case, we include that edge in the temporary list $Cut$. 
In the post-processing after the labeling process stops, we exclude from $Cut$ all edges with both end-vertices labeled, thus retaining only the two edges going from the current (labeled) strongly connected component to the next (unlabeled) one.

%\YD {Maybe, Enqueue(QUEUE,$s'$) in line 8?}

%$<<$ Time analysis $>>$
The linear time bound $O(|E|)$ of Algorithm~\ref{alg:cuts} running on $G$ %\textproc{Find2EdgeCuts} 
follows from the  $O(|E|)$  time of each iteration of the Ford-Fulkerson algorithm and of BFS.
%Two Ford-Fulkerson algorithm iterations are run initially, and then one more specialized iteration where
All operations deviating from the standard BFS are either constant-time or done at most twice per edge.
That is, the runtime of procedure \textproc{Find2EdgeCuts} on $C^-$ is proportional to the number edges therein, which is $O(|E|)$
total in all components
(see Section~\ref{sec:SPQR}).
Hence, the total time of all runs of \textproc{Find2EdgeCuts} in Algorithm~\ref{alg:EPE} on $G$ is $O(|E|)$.

\begin{algorithm}[t!]
	\caption{\textit{Finding 2-Edge Cuts between $s'$ and $t'$} }\label{alg:cuts}
	\hspace*{\algorithmicindent} 
	\textbf{Input:} graph $C^-$, such that
	$C^- \cup  (s',t')$ is 3-connected %, its vertices $s, t$
	 \\
	\hspace*{\algorithmicindent} 
	\textbf{Output:} list of edge pairs $Cuts$
	\begin{algorithmic}[1]
		\Function{Find2EdgeCuts}{$C^-,(s',t')$}
%	    \State Let $C^-= C \setminus \{(s',t')\} $
%	    \State $Cuts$ $\leftarrow \emptyset$
	    \State Construct flow network $N$ from $C^-$: each edge of $C^-$ 
	    \State ~~~generates two anti-parallel edges in $N$ of capacity 1 each,
	    \State ~~~source=$s'$, target=$t'$
	    \State Execute two iterations of the Ford-Fulkerson algorithm on $N$
%	    \State ~~~ resulting in flow $f$ of size 2 and residual network $N_f$
	    \State $Cuts$ $\leftarrow \emptyset$, QUEUE $\leftarrow \emptyset$, Cut $\leftarrow \emptyset$
	    \State Unlabel all vertices in $C^-$
	    \State Label $s'$ and Enqueue(QUEUE,$\{ s'\}$)
	    %\YD{Delete: /* Set $j\leftarrow 0$} */
%
    \While {vertex $t'$ is not labeled}
          \If {Empty(QUEUE)}
          \label{l:empty}
             \State Remove from Cut all $(u,v)$
                where $u,v$ are both labeled
            \State Add Cut to $Cuts$
            \State  Let $(u',u''),(v',v'')$ be the edges in Cut s.t. $u',v'$ are labeled
%             \State Remove $(u',u'')$ and $(v',v'')$ from $N$
             \State Label $u''$ and $v''$ and Enqueue(QUEUE,$\{ u'',v''\}$)
             \State Cut $\leftarrow \emptyset$  
%	        \State Unlabel all vertices in $C^-$
             %~~~~~ /* \YD{Delete: Set $j
 %            \leftarrow j+1$ } */
%
          \Else
             \State $v\leftarrow$ Dequeue(QUEUE)
             \For{all $(v,v')\in C^-$ 
                s.t. $v'$ is not labeled}
                \If{ $(v,v')$ has flow=1 assigned} %in $C^-$ 
                   \State Add $(v,v')$ to Cut
                \Else 
                   \State Label $v'$ and Enqueue(QUEUE,$\{ v' \} $) \label{l:addvertex}
                 \EndIf
             \EndFor
          \EndIf
    \EndWhile
    \State  {\bf return} $Cuts$
    \EndFunction
    \end{algorithmic} %\vspace{-0.4cm}
\end{algorithm}

%\noindent
%*************************************

%Let us contract each strongly connected component of the residual network $N_f$ w.r.t. a maximal flow $f$ into a (super-)node, retaining the edges going between them; denote the resulting (multi-) graph by $\tilde G$; we denote by 
%Let $\tilde G$ be the (super-)graph of all strongly connected components of the residual network w.r.t. a maximal flow . 
%Let us call an edge $s,t$-cut $E'$ \emph{directed} if all its edges 
%Then, the minimum $s,t$-cuts in $N$ are induced exactly by the ``directed'' edge cuts in $\tilde G$, that is the edge 

%********************************************************
\bibliography{LSP}

\providecommand\pagehyphen{-}
\begin{thebibliography}{10}

\bibitem{DBLP:journals/algorithmica/BattistaT96}
Giuseppe~Di Battista and Roberto Tamassia.
\newblock On-line maintenance of triconnected components with {SPQR}-trees.
\newblock {\em Algorithmica}, 15(4):302--318, 1996.

\bibitem{berczi_et_al:LIPIcs:2017:7824}
Kristof Berczi and Yusuke Kobayashi.
\newblock {The Directed Disjoint Shortest Paths Problem}.
\newblock In Kirk Pruhs and Christian Sohler, editors, {\em 25th Annual
  European Symposium on Algorithms (ESA 2017)}, volume~87 of {\em Leibniz
  International Proceedings in Informatics (LIPIcs)}, pages 13:1--13:13,
  Dagstuhl, Germany, 2017. Schloss Dagstuhl--Leibniz-Zentrum fuer Informatik.

\bibitem{10.5555/1614191}
Thomas~H. Cormen, Charles~E. Leiserson, Ronald~L. Rivest, and Clifford Stein.
\newblock {\em Introduction to Algorithms, Third Edition}.
\newblock The MIT Press, 3rd edition, 2009.

\bibitem{ShimonyEtAl2022SOCS}
G.~Dahan, I.~Tabib, S.~E. Shimony, and A.~Felner.
\newblock Generalized longest path problems.
\newblock In {\em {Symposium on Combinatorial Search}}, pages 56--64, 2022.

\bibitem{DKL}
E.~Dinic, Alexander Karzanov, and M.~Lomonosov.
\newblock The system of minimum edge cuts in a graph.
\newblock {\em In book: Issledovaniya po Diskretnoǐ Optimizatsii (Engl. title:
  Studies in Discrete Optimizations), A.A. Fridman, ed., Nauka, Moscow,
  290-306, in Russian,}, 01 1976.

\bibitem{Dirac1960}
Gabriel~Andrew Dirac.
\newblock In abstrakten graphen vorhandene vollständige 4-graphen und ihre
  unterteilungen.
\newblock {\em Mathematische Nachrichten}, 22(1-2):61--85, 1960.

\bibitem{EILAMTZOREFF1998113}
Tali Eilam-Tzoreff.
\newblock The disjoint shortest paths problem.
\newblock {\em Discrete Applied Mathematics}, 85(2):113--138, 1998.

\bibitem{even2011graph}
S.~Even.
\newblock {\em Graph Algorithms}.
\newblock Cambridge University Press, 2011.

\bibitem{10.5555/1942094}
D.~R. Ford and D.~R. Fulkerson.
\newblock {\em Flows in Networks}.
\newblock Princeton University Press, USA, 2010.

\bibitem{garey1979computers}
M.~R. Garey and D.~S. Johnson.
\newblock {\em Computers and Intractability: A Guide to the Theory of
  NP-Completeness (Series of Books in the Mathematical Sciences)}.
\newblock W. H. Freeman, first edition edition, 1979.

\bibitem{SPQRtrees}
Carsten Gutwenger and Petra Mutzel.
\newblock A linear time implementation of {SPQR}-trees.
\newblock In {\em Proceedings of the 8th International Symposium on Graph
  Drawing}, 9 2000.

\bibitem{Harary1969}
F.~Harary.
\newblock {\em Graph Theory}.
\newblock Addison-Wesley, Reading, MA, 1969.

\bibitem{doi:10.1137/0202012}
J.~E. Hopcroft and R.~E. Tarjan.
\newblock Dividing a graph into triconnected components.
\newblock {\em SIAM Journal on Computing}, 2(3):135--158, 1973.

\bibitem{10.5555/1051910}
Jon Kleinberg and Eva Tardos.
\newblock {\em Algorithm Design}.
\newblock Addison-Wesley Longman Publishing Co., Inc., USA, 2005.

\bibitem{10.1215/S0012-7094-37-00336-3}
Saunders~Mac Lane.
\newblock {A structural characterization of planar combinatorial graphs}.
\newblock {\em Duke Mathematical Journal}, 3(3):460 -- 472, 1937.

\bibitem{k-connected}
Misha Lavrov.
\newblock k-connected graph, there exists a cycle that contains any 2 edges and
  any k-2 vertices.
\newblock
  \url{https://math.stackexchange.com/questions/3599622/k-connected-graph-there-exists-a-cycle-that-contains-any-2-edges-and-any-k-2?noredirect=1\&lq=1},
  2020.

\bibitem{10.1145/1061425.1061430}
James~F. Lynch.
\newblock The equivalence of theorem proving and the interconnection problem.
\newblock {\em ACM Sigda Newsletter}, 5(3):31–36, sep 1975.

\bibitem{k-connect}
Nick Matteo.
\newblock G a k-connected graph. show that a set of k-2 vertices and a set of
  two edges lie on a common cycle.
\newblock
  \url{https://math.stackexchange.com/questions/3600673/g-a-k-connected-graph-show-that-a-set-of-k-2-vertices-and-a-set-of-two-edges-li?rq=1},
  2020.

\bibitem{Menger1927}
Karl Menger.
\newblock Zur allgemeinen kurventheorie.
\newblock {\em Fundamenta Mathematicae}, 10(1):96--115, 1927.

\bibitem{PQ}
Jean-Claude Picard and Maurice Queyranne.
\newblock On the structure of all minimum cuts in a network and applications.
\newblock {\em Mathematical Programming}, 22(121), 1982.

\bibitem{DisjointTrees}
B.~A. Reed, N.~Robertson, A.~Schrijver, and P.~D. Seymour.
\newblock Finding disjoint trees in planar graphs in linear time.
\newblock {\em Contemporary Mathematics}, 14, 1993.

\bibitem{DBLP:journals/jct/RobertsonS95b}
Neil Robertson and Paul~D. Seymour.
\newblock Graph minors .xiii. the disjoint paths problem.
\newblock {\em J. Comb. Theory, Ser. {B}}, 63(1):65--110, 1995.

\bibitem{DBLP:journals/networks/Seymour80}
Paul~D. Seymour.
\newblock Four-terminus flows.
\newblock {\em Networks}, 10(1):79--86, 1980.

\bibitem{5366808}
Hars Vardhan, Shreejith Billenahalli, Wanjun Huang, Miguel Razo, Arularasi
  Sivasankaran, Limin Tang, Paolo Monti, Marco Tacca, and Andrea Fumagalli.
\newblock Finding a simple path with multiple must-include nodes.
\newblock In {\em 2009 IEEE International Symposium on Modeling, Analysis \&
  Simulation of Computer and Telecommunication Systems}, pages 1--3, 2009.

\end{thebibliography}

\end{document}